# Properties and applications of transversal operators

Daniel J. Greenhoe

*Abstract*: This paper presents some properties and applications of "transversal operators". Two transversal operators are presented: a "translation" operator $\mathbf{T}$ and a "dilation" operator $\mathbf{D}$. Such operators are used in common analysis systems including Fourier series analysis, Fourier analysis, Gabor analysis, multiresolution analysis (MRA), and wavelet analysis. Like the unitary Fourier transform operator $\bar{\mathbf{F}}$, the transversal operators $\mathbf{T}$ and $\mathbf{D}$ are *unitary*. Demonstrations of the usefulness of these three unitary operators are found in the proofs of results found in some common analytic systems including MRA analysis and wavelet analysis.

# Contents

# 1  Background: operators on linear spaces

## 1.1  Star-algebras

**Definition 1.1**  [1] Let $(\mathbb{F}, +, \cdot)$ be a field. Let $X$ be a set and let $+$ be an operator in $X^{X^2}$ and $\otimes$ be an operator in $X^{\mathbb{F} \times X}$. The structure $\boldsymbol{L} \triangleq \left( X, +, \cdot, (\mathbb{F}, \dot{+}, \dot{\times}) \right)$ is a **linear space** over the field $(\mathbb{F}, +, \cdot)$ if

| | | | | | |
|---|---|---|---|---|---|
| 1. | $\exists \mathbb{0} \in X$ | such that | $x + \mathbb{0} = x$ | $\forall x \in X$ | ($+$ *identity*) |
| 2. | $\exists y \in X$ | such that | $x + y = \mathbb{0}$ | $\forall x \in X$ | ($+$ *inverse*) |
| 3. | | | $(x + y) + z = x + (y + z)$ | $\forall x,y,z \in X$ | ($+$ is *associative*) |
| 4. | | | $x + y = y + x$ | $\forall x,y \in X$ | ($+$ is *commutative*) |
| 5. | | | $1 \cdot x = x$ | $\forall x \in X$ | ($\cdot$ *identity*) |
| 6. | | | $\alpha \cdot (\beta \cdot x) = (\alpha \cdot \beta) \cdot x$ | $\forall \alpha,\beta \in S$ and $x \in X$ | ($\cdot$ *associates* with $\cdot$) |
| 7. | | | $\alpha \cdot (x + y) = (\alpha \cdot x) + (\alpha \cdot y)$ | $\forall \alpha \in S$ and $x,y \in X$ | ($\cdot$ *distributes* over $+$) |
| 8. | | | $(\alpha + \beta) \cdot x = (\alpha \cdot x) + (\beta \cdot x)$ | $\forall \alpha,\beta \in S$ and $x \in X$ | ($\cdot$ *pseudo-distributes* over $+$) |

The set $X$ is called the **underlying set**. The elements of $X$ are called **vectors**. The elements of $\mathbb{F}$ are called **scalars**. A linear space is also called a **vector space**. If $\mathbb{F} \triangleq \mathbb{R}$, then $\boldsymbol{L}$ is a **real linear space**. If $\mathbb{F} \triangleq \mathbb{C}$, then $\boldsymbol{L}$ is a **complex linear space**.

All *linear space*s are equipped with an operation by which vectors in the spaces can be added together. Linear spaces also have an operation that allows a scalar and a vector to be "multiplied" together. But linear spaces in general have no operation that allows two vectors to be multiplied together. A linear space together with such an operator is an **algebra**.[2]

---

[1]   [84], pages 40–41, ⟨Definition 2.1 and following remarks⟩,  [56], page 41, ☞ [58], pages 1–2, ☞ [103], ⟨Chapter IX⟩, ☞ [104], pages 119–120, ▣ [12], pages 134–135

**Definition 1.2** [3] Let $A$ be an *algebra*.

An algebra $A$ is **unital** if $\exists u \in A$ such that $ux = xu = x$ $\forall x \in A$

**Definition 1.3** [4] Let $A$ be an *algebra*.

The pair $(A, *)$ is a $*$-**algebra** if

| | | | | | |
|---|---|---|---|---|---|
| 1. | $(x+y)^*$ | $=$ | $x^* + y^*$ | $\forall_{x,y \in A}$ | *(distributive)* and |
| 2. | $(\alpha x)^*$ | $=$ | $\bar{\alpha} x^*$ | $\forall_{x \in A,\ \alpha \in \mathbb{C}}$ | *(conjugate linear)* and |
| 3. | $(xy)^*$ | $=$ | $y^* x^*$ | $\forall_{x,y \in A}$ | *(antiautomorphic)* and |
| 4. | $x^{**}$ | $=$ | $x$ | $\forall_{x \in A}$ | *(involutory)* |

The operator $*$ is called an **involution** on the algebra $A$.

**Proposition 1.4** [5] *Let $(A, *)$ be an* UNITAL $*$-ALGEBRA.

$$x \text{ is invertible} \implies \begin{cases} 1. & x^* \text{ is INVERTIBLE} & \forall_{x \in A} \quad and \\ 2. & (x^*)^{-1} = (x^{-1})^* & \forall_{x \in A} \end{cases}$$

**Definition 1.5** [6] Let $(A, \|\cdot\|)$ be a $*$-*algebra* (Definition 1.3 page 4).

An element $x \in A$ is **hermitian** or **self-adjoint** if $x^* = x$.

An element $x \in A$ is **normal** if $xx^* = x^*x$.

An element $x \in A$ is a **projection** if $xx = x$ *(involutory)* and $x^* = x$ *(hermitian)*.

**Theorem 1.6** [7] *Let $(A, \|\cdot\|)$ be a $*$-ALGEBRA (Definition 1.3 page 4).*

$$\underbrace{x = x^* \text{ and } y = y^*}_{x \text{ and } y \text{ are hermitian}} \implies \begin{cases} x + y = (x+y)^* & \text{\small($x + y$ is self adjoint)} \\ x^* = (x^*)^* & \text{\small($x^*$ is self adjoint)} \\ \underbrace{xy = (xy)^*}_{(xy) \text{ is hermitian}} \iff \underbrace{xy = yx}_{commutative} \end{cases}$$

**Definition 1.7** (Hermitian components) [8] Let $(A, \|\cdot\|)$ be a $*$-*algebra* (Definition 1.3 page 4).

The **real part** of $x$ is defined as $\Re x \triangleq \dfrac{1}{2}\left(x + x^*\right)$

The **imaginary part** of $x$ is defined as $\Im x \triangleq \dfrac{1}{2i}\left(x - x^*\right)$

---

[2] It has been estimated that as of 2005, there are "somewhere between 50,000 and 200,000" algebras. [62], page v

[3] [38], page 1

[4] [109], page 178, [46], page 241

[5] [38], page 5

[6] [109], page 178, [46], page 242

[7] [97], page 429

[8] [97], page 430, [109], page 179, [46], page 242

**Theorem 1.8** [9] *Let* $(A, *)$ *be a* ∗-ALGEBRA *(Definition 1.3 page 4).*

$$\Re x = (\Re x)^* \qquad \forall x \in A \qquad (\Re x \text{ is hermitian})$$
$$\Im x = (\Im x)^* \qquad \forall x \in A \qquad (\Im x \text{ is hermitian})$$

**Theorem 1.9** (Hermitian representation) [10] *Let* $(A, *)$ *be a* ∗-ALGEBRA *(Definition 1.3 page 4).*

$$a = x + iy \qquad \Longleftrightarrow \qquad x = \Re a \quad and \quad y = \Im a$$

**Definition 1.10** [11] Let $A$ be an algebra.

The pair $(A, \|\cdot\|)$ is a **normed algebra** if

$$\|xy\| \le \|x\| \, \|y\| \qquad \forall x, y \in A \qquad \text{(multiplicative condition)}$$

A normed algebra $(A, \|\cdot\|)$ is a **Banach algebra** if $(A, \|\cdot\|)$ is also a Banach space.

**Proposition 1.11**

$$(A, \|\cdot\|) \text{ is a normed algebra} \qquad \Longrightarrow \qquad \text{multiplication is } \textbf{continuous} \text{ in } (A, \|\cdot\|)$$

**Definition 1.12** [12]

The triple $(A, \|\cdot\|, *)$ is a $C^*$ **algebra** if

1. $(A, \|\cdot\|)$      is a Banach algebra   and
2. $(A, *)$      is a ∗-algebra    and
3. $\|x^* x\| = \|x\|^2 \quad \forall x \in A.$

**Theorem 1.13** [13] *Let* $A$ *be an algebra.*

$$(A, \|\cdot\|, *) \text{ is a } C^* \textbf{ algebra} \qquad \Longrightarrow \qquad \|x^*\| = \|x\|$$

## 1.2 Operators on linear spaces

### 1.2.1 Operator Algebra

An operator is simply a function that maps from a linear space to another linear space (or to the same linear space).

**Definition 1.14** [14] A function $\mathbf{A}$ in $Y^X$ is an **operator** in $Y^X$ if $X$ and $Y$ are both *linear spaces*.

---

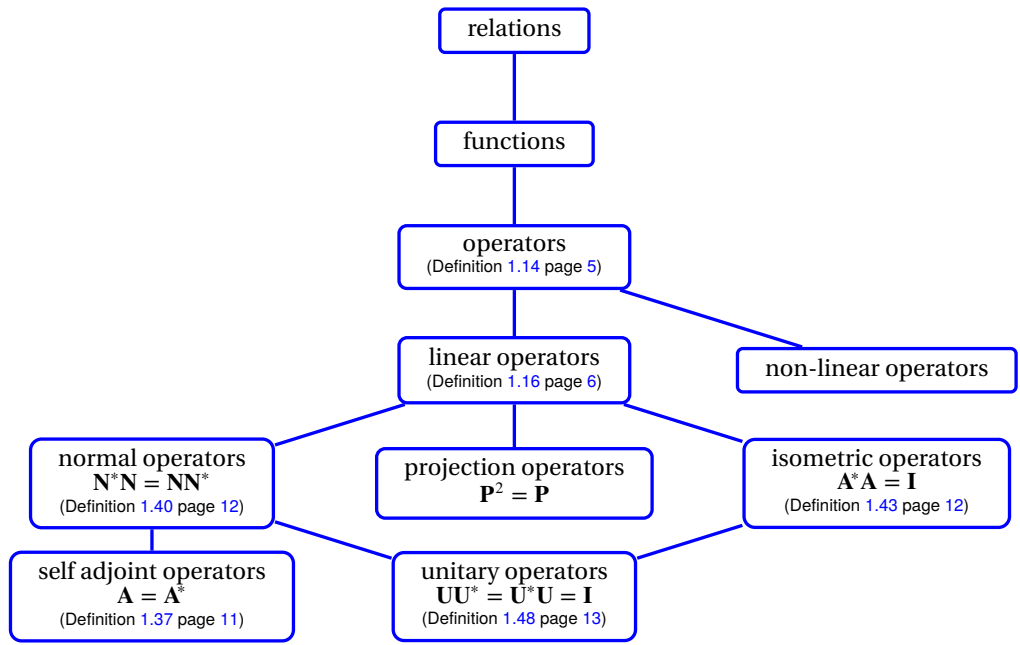

Figure 1: Some operator types

Two operators $\mathbf{A}$ and $\mathbf{B}$ in $\boldsymbol{Y}^{\boldsymbol{X}}$ are **equal** if $\mathbf{A}\boldsymbol{x} = \mathbf{B}\boldsymbol{x}$ for all $\boldsymbol{x} \in \boldsymbol{X}$. The inverse relation of an operator $\mathbf{A}$ in $\boldsymbol{Y}^{\boldsymbol{X}}$ always exists as a *relation* in $2^{\boldsymbol{X}\boldsymbol{Y}}$, but may not always be a *function* (may not always be an operator) in $\boldsymbol{Y}^{\boldsymbol{X}}$.

The operator $\mathbf{I} \in \boldsymbol{X}^{\boldsymbol{X}}$ is the *identity* operator if $\mathbf{I}\boldsymbol{x} = \mathbf{I}$ for all $\boldsymbol{x} \in \boldsymbol{X}$.

**Definition 1.15**  [15]  Let $\boldsymbol{X}^{\boldsymbol{X}}$ be the set of all operators with from a linear space $\boldsymbol{X}$ to $\boldsymbol{X}$. Let $\mathbf{I}$ be an operator in $\boldsymbol{X}^{\boldsymbol{X}}$. Let $\circledast(\boldsymbol{X})$ be the *identity element* in $\boldsymbol{X}^{\boldsymbol{X}}$.

$\qquad$ $\mathbf{I}$ is the **identity operator** in $\boldsymbol{X}^{\boldsymbol{X}}$ if $\qquad$ $\mathbf{I} = \circledast(\boldsymbol{X})$.

### 1.2.2   Linear operators

**Definition 1.16**  [16]  Let $\boldsymbol{X} \triangleq \big(X, +, \cdot, (\mathbb{F}, \dotplus, \dot\times)\big)$ and $\boldsymbol{Y} \triangleq \big(Y, +, \cdot, (\mathbb{F}, \dotplus, \dot\times)\big)$ be linear spaces.

---

[15]  🕮  [97], page 411
[16]  🕮  [84], page 55, 🕮 [3], page 224, 🕮 [66], page 6, 🕮 [126], page 33

An operator $\mathbf{L} \in Y^X$ is **linear** if

    1.    $\mathbf{L}(x+y)$  $=$  $\mathbf{L}x + \mathbf{L}y$     $\forall x,y \in X$       (*additive*)      and

    2.    $\mathbf{L}(\alpha x)$  $=$  $\alpha \mathbf{L}x$      $\forall x \in X, \quad \forall \alpha \in \mathbb{F}$     (*homogeneous*).

The set of all linear operators from $X$ to $Y$ is denoted $\mathcal{L}(X, Y)$ such that

$$\mathcal{L}(X, Y) \triangleq \left\{ \mathbf{L} \in Y^X \, \middle| \, \mathbf{L} \text{ is linear} \right\} \qquad .$$

**Theorem 1.17**   [17]   *Let* $\mathbf{L}$ *be an operator from a linear space* $X$ *to a linear space* $Y$, *both over a field* $\mathbb{F}$.

$$\mathbf{L} \text{ is LINEAR} \quad \Longrightarrow \quad \begin{cases} 1. & \mathbf{L}\mathbb{0} & = & \mathbb{0} \\ 2. & \mathbf{L}(-x) & = & -(\mathbf{L}x) & \forall x \in X \\ 3. & \mathbf{L}(x - y) & = & \mathbf{L}x - \mathbf{L}y & \forall x,y \in X \\ 4. & \mathbf{L}\left(\displaystyle\sum_{n=1}^{N} \alpha_n x_n\right) & = & \displaystyle\sum_{n=1}^{N} \alpha_n (\mathbf{L}x_n) & x_n \in X, \alpha_n \in \mathbb{F} \end{cases}$$

**Theorem 1.18**   [18]   *Let* $\mathcal{L}(X, Y)$ *be the set of all linear operators from a linear space* $X$ *to a linear space* $Y$. *Let* $\mathcal{N}(\mathbf{L})$ *be the* NULL SPACE *of an operator* $\mathbf{L}$ *in* $Y^X$ *and* $\mathcal{I}(\mathbf{L})$ *the* IMAGE SET *of* $\mathbf{L}$ *in* $Y^X$.

     $\mathcal{L}(X, Y)$    *is a linear space*                 (*space of linear transforms*)

     $\mathcal{N}(\mathbf{L})$    *is a linear subspace of* $X$     $\forall \mathbf{L} \in Y^X$

     $\mathcal{I}(\mathbf{L})$    *is a linear subspace of* $Y$     $\forall \mathbf{L} \in Y^X$

**Example 1.19**   [19]   Let $C([a, b], \mathbb{R})$ be the set of all *continuous* functions from the closed real interval $[a, b]$ to $\mathbb{R}$.

     $C([a, b], \mathbb{R})$ is a linear space.

**Theorem 1.20**   [20]   *Let* $\mathcal{L}(X, Y)$ *be the set of linear operators from a linear space* $X$ *to a linear space* $Y$. *Let* $\mathcal{N}(\mathbf{L})$ *be the* NULL SPACE *of a linear operator* $\mathbf{L} \in \mathcal{L}(X, Y)$.

     $\mathbf{L}x = \mathbf{L}y$     $\Longleftrightarrow$     $x - y \in \mathcal{N}(\mathbf{L})$

     $\mathbf{L}$ *is* INJECTIVE     $\Longleftrightarrow$     $\mathcal{N}(\mathbf{L}) = \{\mathbb{0}\}$

**Theorem 1.21**   [21]   *Let* $W$, $X$, $Y$, *and* $Z$ *be linear spaces over a field* $\mathbb{F}$.

    1.   $\mathbf{L}(\mathbf{M}\mathbf{N})$  $=$  $(\mathbf{L}\mathbf{M})\mathbf{N}$     $\forall \mathbf{L} \in \mathcal{L}(Z,W), \mathbf{M} \in \mathcal{L}(Y,Z), \mathbf{N} \in \mathcal{L}(X,Y)$    (ASSOCIATIVE)

    2.   $\mathbf{L}(\mathbf{M} \mathbin{\dot{+}} \mathbf{N})$  $=$  $(\mathbf{L}\mathbf{M}) \mathbin{\dot{+}} (\mathbf{L}\mathbf{N})$     $\forall \mathbf{L} \in \mathcal{L}(Y,Z), \mathbf{M} \in \mathcal{L}(X,Y), \mathbf{N} \in \mathcal{L}(X,Y)$    (LEFT DISTRIBUTIVE)

    3.   $(\mathbf{L} \mathbin{\dot{+}} \mathbf{M})\mathbf{N}$  $=$  $(\mathbf{L}\mathbf{N}) \mathbin{\dot{+}} (\mathbf{M}\mathbf{N})$     $\forall \mathbf{L} \in \mathcal{L}(Y,Z), \mathbf{M} \in \mathcal{L}(Y,Z), \mathbf{N} \in \mathcal{L}(X,Y)$    (RIGHT DISTRIBUTIVE)

    4.   $\alpha(\mathbf{L}\mathbf{M})$  $=$  $(\alpha\mathbf{L})\mathbf{M} = \mathbf{L}(\alpha\mathbf{M})$     $\forall \mathbf{L} \in \mathcal{L}(Y,Z), \mathbf{M} \in \mathcal{L}(X,Y), \alpha \in \mathbb{F}$    (HOMOGENEOUS)

---

## 1.3 Operators on Normed linear spaces

### 1.3.1 Operator norm

**Definition 1.22** [22] Let $\mathcal{L}(\boldsymbol{X}, \boldsymbol{Y})$ be the space of linear operators over normed linear spaces $\boldsymbol{X}$ and $\boldsymbol{Y}$. [23]

The **operator norm** $\|\|\cdot\|\|$ is defined as
$$\|\|\mathbf{A}\|\| \triangleq \sup_{x \in \boldsymbol{X}} \{\, \|\mathbf{A}\boldsymbol{x}\| \mid \|\boldsymbol{x}\| \le 1 \,\} \qquad \forall \mathbf{A} \in \mathcal{L}(\boldsymbol{X}, \boldsymbol{Y})$$

The pair $(\mathcal{L}(\boldsymbol{X}, \boldsymbol{Y}), \|\|\cdot\|\|)$ is the **normed space of linear operators** on $(\boldsymbol{X}, \boldsymbol{Y})$.

Proposition 1.23 (next) shows that the functional defined in Definition 1.22 (previous) is a *norm*.

**Proposition 1.23** [24] *Let $(\mathcal{L}(\boldsymbol{X}, \boldsymbol{Y}), \|\|\cdot\|\|)$ be the normed space of linear operators over the normed linear spaces $\boldsymbol{X} \triangleq (X, +, \cdot, (\mathbb{F}, \dotplus, \dottimes), \|\cdot\|)$ and $\boldsymbol{Y} \triangleq (Y, +, \cdot, (\mathbb{F}, \dotplus, \dottimes), \|\cdot\|)$.*

*The functional $\|\|\cdot\|\|$ is a **norm** on $\mathcal{L}(\boldsymbol{X}, \boldsymbol{Y})$. In particular,*

| | | | | | |
|---|---|---|---|---|---|
| *1.* | $\|\|\mathbf{A}\|\|$ | $\ge$ | $0$ | $\forall \mathbf{A} \in \mathcal{L}(\boldsymbol{X},\boldsymbol{Y})$ | (NON-NEGATIVE) *and* |
| *2.* | $\|\|\mathbf{A}\|\|$ | $=$ | $0 \iff \mathbf{A} \stackrel{\circ}{=} \mathbb{0}$ | $\forall \mathbf{A} \in \mathcal{L}(\boldsymbol{X},\boldsymbol{Y})$ | (NONDEGENERATE) *and* |
| *3.* | $\|\|\alpha \mathbf{A}\|\|$ | $=$ | $|\alpha| \, \|\|\mathbf{A}\|\|$ | $\forall \mathbf{A} \in \mathcal{L}(\boldsymbol{X},\boldsymbol{Y}), \alpha \in \mathbb{F}$ | (HOMOGENEOUS) *and* |
| *4.* | $\|\|\mathbf{A} \dotplus \mathbf{B}\|\|$ | $\le$ | $\|\|\mathbf{A}\|\| + \|\|\mathbf{B}\|\|$ | $\forall \mathbf{A} \in \mathcal{L}(\boldsymbol{X},\boldsymbol{Y})$ | (SUBADDITIVE). |

*Moreover, $(\mathcal{L}(\boldsymbol{X}, \boldsymbol{Y}), \|\|\cdot\|\|)$ is a **normed linear space**.*

**Lemma 1.24** [25] *Let $(\mathcal{L}(\boldsymbol{X}, \boldsymbol{Y}), \|\|\cdot\|\|)$ be the normed space of linear operators over normed linear spaces $\boldsymbol{X} \triangleq (X, +, \cdot, (\mathbb{F}, \dotplus, \dottimes), \|\cdot\|)$ and $\boldsymbol{Y} \triangleq (Y, +, \cdot, (\mathbb{F}, \dotplus, \dottimes), \|\cdot\|)$.*
$$\|\|\mathbf{L}\|\| = \sup_{x} \{\, \|\mathbf{L}\boldsymbol{x}\| \mid \|\boldsymbol{x}\| = 1 \,\} \qquad \forall \boldsymbol{x} \in \mathcal{L}(\boldsymbol{X}, \boldsymbol{Y})$$

**Proposition 1.25** [26] *Let $\mathbf{I}$ be the identity operator in the normed space of linear operators $(\mathcal{L}(\boldsymbol{X}, \boldsymbol{X}), \|\|\cdot\|\|)$.*
$$\|\|\mathbf{I}\|\| = 1$$

**Theorem 1.26** [27] *Let $(\mathcal{L}(\boldsymbol{X}, \boldsymbol{Y}), \|\|\cdot\|\|)$ be the normed space of linear operators over normed linear spaces $\boldsymbol{X}$ and $\boldsymbol{Y}$.*

| | | | |
|---|---|---|---|
| $\|\mathbf{L}\boldsymbol{x}\|$ | $\le$ | $\|\|\mathbf{L}\|\| \, \|\boldsymbol{x}\|$ | $\forall \mathbf{L} \in \mathcal{L}(\boldsymbol{X}, \boldsymbol{Y}), \, \boldsymbol{x} \in \boldsymbol{X}$ |
| $\|\|\mathbf{KL}\|\|$ | $\le$ | $\|\|\mathbf{K}\|\| \, \|\|\mathbf{L}\|\|$ | $\forall \mathbf{K}, \mathbf{L} \in \mathcal{L}(\boldsymbol{X}, \boldsymbol{Y})$ |

---

### 1.3.2 Bounded linear operators

**Definition 1.27** [28]   Let $(\mathcal{L}(X, Y), \|\|\cdot\|\|)$ be a normed space of linear operators.
     An operator $\mathbf{B}$ is **bounded** if $\|\|\mathbf{B}\|\| < \infty$.
     The quantity $\mathcal{B}(X, Y)$ is the set of all **bounded linear operators** on $(X, Y)$ such that
     $\mathcal{B}(X, Y) \triangleq \{\mathbf{L} \in \mathcal{L}(X, Y) \mid \|\|\mathbf{L}\|\| < \infty\}$.

**Theorem 1.28** [29] *Let $(\mathcal{L}(X, Y), \|\|\cdot\|\|)$ be the set of linear operators over normed linear spaces*
$X \triangleq \left(X, +, \cdot, (\mathbb{F}, \dotplus, \dottimes), \|\cdot\|\right)$ *and* $Y \triangleq \left(Y, +, \cdot, (\mathbb{F}, \dotplus, \dottimes), \|\cdot\|\right)$.
     *The following conditions are all* EQUIVALENT:

| | | | |
|---|---|---|---|
| 1. | $\mathbf{L}$ *is continuous at* A SINGLE POINT $x_0 \in X$ | $\forall \mathbf{L} \in \mathcal{L}(X, Y)$ | $\Longleftrightarrow$ |
| 2. | $\mathbf{L}$ *is* CONTINUOUS *(at every point* $x \in X$) | $\forall \mathbf{L} \in \mathcal{L}(X, Y)$ | $\Longleftrightarrow$ |
| 3. | $\|\|\mathbf{L}\|\| < \infty$ (*$\mathbf{L}$ is* BOUNDED) | $\forall \mathbf{L} \in \mathcal{L}(X, Y)$ | $\Longleftrightarrow$ |
| 4. | $\exists M \in \mathbb{R}$   such that   $\|\mathbf{L}x\| \leq M \|x\|$ | $\forall \mathbf{L} \in \mathcal{L}(X, Y), x \in X$ | |

### 1.3.3 Adjoints on normed linear spaces

**Definition 1.29**   Let $\mathcal{B}(X, Y)$ be the space of bounded linear operators on normed linear spaces $X$ and $Y$. Let $X^*$ be the *topological dual space* of $X$.
     $\mathbf{B}^*$ is the **adjoint** of an operator $\mathbf{B} \in \mathcal{B}(X, Y)$ if
     $f(\mathbf{B}x) = \left[\mathbf{B}^* f\right](x)$      $\forall f \in X^*, x \in X$

**Theorem 1.30** [30] *Let $\mathcal{B}(X, Y)$ be the space of bounded linear operators on normed linear spaces $X$ and $Y$.*

| | | | |
|---|---|---|---|
| $(\mathbf{A} \dotplus \mathbf{B})^*$ | $=$ | $\mathbf{A}^* \dotplus \mathbf{B}^*$ | $\forall \mathbf{A}, \mathbf{B} \in \mathcal{B}(X, Y)$ |
| $(\lambda \mathbf{A})^*$ | $=$ | $\lambda \mathbf{A}^*$ | $\forall \mathbf{A}, \mathbf{B} \in \mathcal{B}(X, Y)$ |
| $(\mathbf{A}\mathbf{B})^*$ | $=$ | $\mathbf{B}^* \mathbf{A}^*$ | $\forall \mathbf{A}, \mathbf{B} \in \mathcal{B}(X, Y)$ |

**Theorem 1.31** [31] *Let $\mathcal{B}(X, Y)$ be the space of bounded linear operators on normed linear spaces $X$ and $Y$. Let $\mathbf{B}^*$ be the adjoint of an operator $\mathbf{B}$.*
     $\|\|\mathbf{B}\|\| = \|\|\mathbf{B}^*\|\|$      $\forall \mathbf{B} \in \mathcal{B}(X, Y)$

---

## 1.4   Operators on Inner-product spaces

**Theorem 1.32** [32] *Let* $\mathbf{A}, \mathbf{B} \in \mathcal{B}(X, X)$ *be bounded linear operators on an inner-product space* $X \triangleq \big( X, +, \cdot, (\mathbb{F}, \dotplus, \dottimes), \langle \triangle \,|\, \triangledown \rangle \big)$.

$$\begin{aligned}
\langle \mathbf{B}x \mid x \rangle &= 0 &\forall x \in X &\quad\Longleftrightarrow\quad& \mathbf{B}x &= \mathbb{0} &\forall x \in X\\
\langle \mathbf{A}x \mid x \rangle &= \langle \mathbf{B}x \mid x \rangle &\forall x \in X &\quad\Longleftrightarrow\quad& \mathbf{A} &= \mathbf{B}
\end{aligned}$$

A fundamental concept of operators on inner-product spaces is the *operator adjoint* (Proposition 1.33 page 10). The adjoint of an operator is a kind of generalization of the conjugate of a complex number in that

- ✎ Both are *star-algebras* (Theorem 1.35 page 10).
- ✎ Both support decomposition into "real" and "imaginary" parts (Theorem 1.9 page 5).

Structurally, the operator adjoint provides a convenient symmetric relationship between the *range space* and *null space* of an operator (Theorem 1.36 page 11).

**Proposition 1.33** [33] *Let* $\mathcal{B}(H, H)$ *be the space of bounded linear operators on a Hilbert space* $H$.[34]

> *An operator* $\mathbf{B}^*$ *is the* ADJOINT *of* $\mathbf{B} \in \mathcal{B}(H, H)$ *if*
> $$\langle \mathbf{B}x \mid y \rangle = \langle x \mid \mathbf{B}^* y \rangle \qquad \forall x, y \in H.$$

**Example 1.34**   (Matrix algebra: $\mathbf{A}^* = \mathbf{A}^H$)   In matrix algebra,
- ✎ The inner-product operation $\langle x \mid y \rangle$ is represented by $y^H x$.
- ✎ The linear operator is represented as a matrix $A$.
- ✎ The operation of $A$ on vector $x$ is represented as $Ax$.
- ✎ The adjoint of matrix $A$ is the Hermitian matrix $A^H$.

Structures that satisfy the four conditions of the next theorem are known as *\*-algebras* ("*star-algebras*", Definition 1.3 page 4). Other structures which are \*-algebras include the field of complex numbers $\mathbb{C}$ and any ring of complex square $n \times n$ matrices.[35]

**Theorem 1.35**   (operator star-algebra) [36] *Let* $H$ *be a Hilbert space with operators* $\mathbf{A}, \mathbf{B} \in \mathcal{B}(H, H)$ *and with adjoints* $\mathbf{A}^*, \mathbf{B}^* \in \mathcal{B}(H, H)$. *Let* $\bar{\alpha}$ *be the complex conjugate of some* $\alpha \in \mathbb{C}$.

---

[32] ✎ [112], page 310, ⟨Theorem 12.7, Corollary⟩ ✎ [53], page 217, ⟨Theorem C.9⟩
[33] ✎ [97], page 220, ✎ [112], page 311, ✎ [49], page 182, ✎ [100], page 49, ✎ [126], page 41
[34] *bounded operator*:    Definition 1.27 (page 9)
      *adjoint*:                Definition 1.29 (page 9)
[35] ✎ [114], page 1
[36] ✎ [59], pages 39–40, ✎ [112], page 311, ✎ [53], pages 219–220, ⟨Theorem C.10⟩

*The pair $(H, *)$ is a $*$-ALGEBRA (STAR-ALGEBRA). In particular,*

| | | | | | | |
|---|---|---|---|---|---|---|
| 1. | $(\mathbf{A} \dotplus \mathbf{B})^* $ | $=$ | $\mathbf{A}^* + \mathbf{B}^*$ | $\forall \mathbf{A}, \mathbf{B} \in H$ | (DISTRIBUTIVE) | *and* |
| 2. | $(\alpha \mathbf{A})^*$ | $=$ | $\bar{\alpha} \mathbf{A}^*$ | $\forall \mathbf{A}, \mathbf{B} \in H$ | (CONJUGATE LINEAR) | *and* |
| 3. | $(\mathbf{AB})^*$ | $=$ | $\mathbf{B}^* \mathbf{A}^*$ | $\forall \mathbf{A}, \mathbf{B} \in H$ | (ANTIAUTOMORPHIC) | *and* |
| 4. | $\mathbf{A}^{**}$ | $=$ | $\mathbf{A}$ | $\forall \mathbf{A}, \mathbf{B} \in H$ | (INVOLUTORY) | |

**Theorem 1.36** [37] *Let $Y^X$ be the set of all operators from a linear space $X$ to a linear space $Y$. Let $\mathcal{N}(\mathbf{L})$ be the NULL SPACE of an operator $\mathbf{L}$ in $Y^X$ and $\mathcal{I}(\mathbf{L})$ the IMAGE SET of $\mathbf{L}$ in $Y^X$.*

$$\mathcal{N}(\mathbf{A}) = \mathcal{I}(\mathbf{A}^*)^\perp$$
$$\mathcal{N}(\mathbf{A}^*) = \mathcal{I}(\mathbf{A})^\perp$$

## 1.5 Special Classes of Operators

### 1.5.1 Self Adjoint Operators

**Definition 1.37** [38] Let $\mathbf{B} \in \mathcal{B}(H, H)$ be a bounded operator with adjoint $\mathbf{B}^*$ on a Hilbert space $H$.

The operator $\mathbf{B}$ is said to be **self-adjoint** or **hermitian** if $\mathbf{B} \triangleq \mathbf{B}^*$.

**Example 1.38** (Autocorrelation operator) Let $\mathsf{x}(t)$ be a random process with autocorrelation

$$R_{\mathsf{xx}}(t, u) \triangleq \underbrace{\mathrm{E}[\mathsf{x}(t)\mathsf{x}^*(u)]}_{\text{expectation}}.$$

Let an autocorrelation operator $\mathbf{R}$ be defined as $[\mathbf{R}\mathsf{f}](t) \triangleq \int_\mathbb{R} \underbrace{R_{\mathsf{xx}}(t, u)}_{\text{kernel}} \mathsf{f}(u) \; \mathrm{d}u$.

$$\mathbf{R} = \mathbf{R}^* \qquad \text{(The autocorrelation operator } \mathbf{R} \text{ is } \textit{self-adjoint)}$$

**Theorem 1.39** [39] *Let $\mathbf{S} : H \to H$ be an operator over a Hilbert space $H$ with eigenvalues $\{\lambda_n\}$ and eigenfunctions $\{\psi_n\}$ such that $\mathbf{S}\psi_n = \lambda_n \psi_n$ and let $\|x\| \triangleq \sqrt{\langle x \mid x \rangle}$.*

$$\underbrace{\mathbf{S} = \mathbf{S}^*}_{\mathbf{S} \text{ is SELF ADJOINT}} \implies \begin{cases} 1. & \langle \mathbf{S}x \mid x \rangle \in \mathbb{R} & \textit{(the hermitian quadratic form of } \mathbf{S} \textit{ is real)} \\ 2. & \lambda_n \in \mathbb{R} & \textit{(eigenvalues of } \mathbf{S} \textit{ are real)} \\ 3. & \lambda_n \neq \lambda_m \implies \langle \psi_n \mid \psi_m \rangle = 0 & \textit{(eigenfunctions associated with distinct eigenvalues are orthogonal)} \end{cases}$$

---

### 1.5.2   Normal Operators

**Definition 1.40**  [40]   Let $\mathcal{B}(\boldsymbol{X}, \boldsymbol{Y})$ be the space of bounded linear operators on normed linear spaces $\boldsymbol{X}$ and $\boldsymbol{Y}$. Let $\mathbf{N}^*$ be the adjoint of an operator $\mathbf{N} \in \mathcal{B}(\boldsymbol{X}, \boldsymbol{Y})$.

$\quad$ $\mathbf{N}$ is **normal** if $\quad \mathbf{N}^*\mathbf{N} = \mathbf{N}\mathbf{N}^*$.

**Theorem 1.41**  [41]  *Let $\mathcal{B}(\boldsymbol{H}, \boldsymbol{H})$ be the space of bounded linear operators on a Hilbert space $\boldsymbol{H}$. Let $\mathcal{N}(\mathbf{N})$ be the* NULL SPACE *of an operator $\mathbf{N}$ in $\mathcal{B}(\boldsymbol{H}, \boldsymbol{H})$ and $\mathcal{I}(\mathbf{N})$ the* IMAGE SET *of $\mathbf{N}$ in $\mathcal{B}(\boldsymbol{H}, \boldsymbol{H})$.*

$$\underbrace{\mathbf{N}^*\mathbf{N} = \mathbf{N}\mathbf{N}^*}_{\mathbf{N}\ is\ normal} \quad\Longleftrightarrow\quad \|\mathbf{N}^*\boldsymbol{x}\| = \|\mathbf{N}\boldsymbol{x}\| \qquad \forall \boldsymbol{x} \in \boldsymbol{H}$$

**Theorem 1.42**  [42]  *Let $\mathcal{B}(\boldsymbol{H}, \boldsymbol{H})$ be the space of bounded linear operators on a Hilbert space $\boldsymbol{H}$. Let $\mathcal{N}(\mathbf{N})$ be the* NULL SPACE *of an operator $\mathbf{N}$ in $\mathcal{B}(\boldsymbol{H}, \boldsymbol{H})$ and $\mathcal{I}(\mathbf{N})$ the* IMAGE SET *of $\mathbf{N}$ in $\mathcal{B}(\boldsymbol{H}, \boldsymbol{H})$.*

$$\underbrace{\mathbf{N}^*\mathbf{N} = \mathbf{N}\mathbf{N}^*}_{\mathbf{N}\ is\ normal} \quad\Longrightarrow\quad \underbrace{\mathcal{N}(\mathbf{N}^*) = \mathcal{N}(\mathbf{N})}_{\mathbf{N}\ and\ \mathbf{N}^*\ have\ the\ same\ null\ space}$$

### 1.5.3   Isometric operators

An operator on a pair of normed linear spaces is *isometric* (next definition) if it is an *isometry*.

**Definition 1.43**   Let $\left( X, +, \cdot, (\mathbb{F}, \dotplus, \dot\times), \|\cdot\| \right)$ and $\left( Y, +, \cdot, (\mathbb{F}, \dotplus, \dot\times), \|\cdot\| \right)$ be *normed linear spaces* (Definition 3.6 page 33).

$\quad$ An operator $\mathbf{M} \in \mathcal{L}(\boldsymbol{X}, \boldsymbol{Y})$ is **isometric** if

$\quad$ $\|\mathbf{M}\boldsymbol{x}\| = \|\boldsymbol{x}\| \qquad \forall x \in X.$

**Theorem 1.44**  [43]  *Let $\left( X, +, \cdot, (\mathbb{F}, \dotplus, \dot\times), \|\cdot\| \right)$ and $\left( Y, +, \cdot, (\mathbb{F}, \dotplus, \dot\times), \|\cdot\| \right)$ be* NORMED LINEAR SPACES. *Let $\mathbf{M}$ be a linear operator in $\mathcal{L}(\boldsymbol{X}, \boldsymbol{Y})$.*

$$\underbrace{\|\mathbf{M}\boldsymbol{x}\| = \|\boldsymbol{x}\| \quad \forall x \in X}_{isometric\ in\ length} \quad\Longleftrightarrow\quad \underbrace{\|\mathbf{M}\boldsymbol{x} - \mathbf{M}\boldsymbol{y}\| = \|\boldsymbol{x} - \boldsymbol{y}\| \quad \forall x, y \in X}_{isometric\ in\ distance}$$

---

Isometric operators have already been defined (Definition 1.43 page 12) in the more general normed linear spaces, while Theorem 1.44 (page 12) demonstrated that in a normed linear space $X$, $\|Mx\| = \|x\| \iff \|Mx - My\| = \|x - y\|$ for all $x, y \in X$. Here in the more specialized inner-product spaces, Theorem 1.45 (next) demonstrates two additional equivalent properties.

**Theorem 1.45** [44] *Let $\mathcal{B}(X, X)$ be the space of bounded linear operators on a normed linear space $X \triangleq \left( X, +, \cdot, \langle \mathbb{F}, \dotplus, \dottimes \rangle, \|\cdot\| \right)$. Let $N$ be a bounded linear operator in $\mathcal{L}(X, X)$, and $I$ the identity operator in $\mathcal{L}(X, X)$. Let $\|x\| \triangleq \sqrt{\langle x \mid x \rangle}$.*
      *The following conditions are all **equivalent**:*

|   |   |   |   |   |   |
|---|---|---|---|---|---|
| 1. | $M^*M$ | $=$ | $I$ | | $\iff$ |
| 2. | $\langle Mx \mid My \rangle$ | $=$ | $\langle x \mid y \rangle$ | $\forall x,y \in X$  (M is surjective) | $\iff$ |
| 3. | $\|Mx - My\|$ | $=$ | $\|x - y\|$ | $\forall x,y \in X$  (isometric in distance) | $\iff$ |
| 4. | $\|Mx\|$ | $=$ | $\|x\|$ | $\forall x \in X$  (isometric in length) | |

**Theorem 1.46** [45] *Let $\mathcal{B}(X, Y)$ be the space of bounded linear operators on normed linear spaces $X$ and $Y$. Let $M$ be a bounded linear operator in $\mathcal{B}(X, Y)$, and $I$ the identity operator in $\mathcal{L}(X, X)$. Let $\Lambda$ be the set of eigenvalues of $M$. Let $\|x\| \triangleq \sqrt{\langle x \mid x \rangle}$.*

$$\underbrace{M^*M = I}_{M \text{ is isometric}} \implies \begin{cases} \|\|M\|\| & = & 1 & \text{\scriptsize (UNIT LENGTH)} \quad and \\ |\lambda| & = & 1 & \forall \lambda \in \Lambda \end{cases}$$

**Example 1.47** (One sided shift operator) [46] Let $X$ be the set of all sequences with range $\mathbb{W} \; (0, 1, 2, \ldots)$ and shift operators defined as

|   |   |   |   |   |
|---|---|---|---|---|
| 1. | $S_r \left( x_0, x_1, x_2, \ldots \right)$ | $\triangleq$ | $\left( 0, x_0, x_1, x_2, \ldots \right)$ | (right shift operator) |
| 2. | $S_l \left( x_0, x_1, x_2, \ldots \right)$ | $\triangleq$ | $\left( x_1, x_2, x_3, \ldots \right)$ | (left shift operator) |

   1. $S_r$ is an isometric operator.
   2. $S_r^* = S_l$

### 1.5.4 Unitary operators

**Definition 1.48** [47] Let $\mathcal{B}(X, Y)$ be the space of bounded linear operators on normed linear spaces $X$ and $Y$. Let $U$ be a bounded linear operator in $\mathcal{B}(X, Y)$, and $I$ the identity operator in $\mathcal{B}(X, X)$.
      The operator $U$ is **unitary** if   $U^*U \;=\; UU^* = I$.

---

**Proposition 1.49**  *Let $\mathcal{B}(X, Y)$ be the space of bounded linear operators on normed linear spaces $X$ and $Y$. Let $\mathbf{U}$ and $\mathbf{V}$ be bounded linear operators in $\mathcal{B}(X, Y)$.*

$$\left.\begin{array}{l}\mathbf{U} \text{ is } \text{UNITARY} \quad and \\ \mathbf{V} \text{ is } \text{UNITARY}\end{array}\right\} \quad \Longrightarrow \quad (\mathbf{UV}) \text{ is } \text{UNITARY}.$$

**Theorem 1.50**  [48]  *Let $\mathcal{B}(H, H)$ be the space of bounded linear operators on a Hilbert space $H$. Let $\mathbf{U}$ be a bounded linear operator in $\mathcal{B}(H, H)$, and $\mathcal{I}(\mathbf{U})$ the IMAGE SET of $\mathbf{U}$.*
   *The following conditions are **equivalent**:*

| | | | | |
|---|---|---|---|---|
| 1. | $\mathbf{UU}^* = \mathbf{U}^*\mathbf{U} = \mathbf{I}$ | | *(unitary)* | $\Longleftrightarrow$ |
| 2. | $\langle \mathbf{U}x \mid \mathbf{U}y\rangle = \langle \mathbf{U}^*x \mid \mathbf{U}^*y\rangle = \langle x \mid y\rangle$ | *and*   $\mathcal{I}(\mathbf{U}) = X$ | *(surjective)* | $\Longleftrightarrow$ |
| 3. | $\|\mathbf{U}x - \mathbf{U}y\| = \|\mathbf{U}^*x - \mathbf{U}^*y\| = \|x - y\|$ | *and*   $\mathcal{I}(\mathbf{U}) = X$ | *(isometric in distance)* | $\Longleftrightarrow$ |
| 4. | $\|\mathbf{U}x\| = \|x\|$ | *and*   $\mathcal{I}(\mathbf{U}) = X$ | *(isometric in length)* | |

**Theorem 1.51**  [49]  *Let $\mathcal{B}(H, H)$ be the space of bounded linear operators on a Hilbert space $H$. Let $\mathbf{U}$ be a bounded linear operator in $\mathcal{B}(H, H)$, $\mathcal{N}(\mathbf{U})$ the NULL SPACE of $\mathbf{U}$, and $\mathcal{I}(\mathbf{U})$ the* IMAGE SET *of $\mathbf{U}$.*

$$\underbrace{\mathbf{UU}^* = \mathbf{U}^*\mathbf{U} = \mathbf{I}}_{\mathbf{U} \text{ is unitary}} \Longrightarrow \left\{ \begin{array}{rcll} \mathbf{U}^{-1} & = & \mathbf{U}^* & \\ \mathcal{I}(\mathbf{U}) & = & \mathcal{I}(\mathbf{U}^*) & = \quad X \\ \mathcal{N}(\mathbf{U}) & = & \mathcal{N}(\mathbf{U}^*) & = \quad \{0\} \\ \|\|\mathbf{U}\|\| & = & \|\|\mathbf{U}^*\|\| & = \quad 1 \quad \text{(UNIT LENGTH)} \end{array} \right.$$

✎PROOF:

(1)   Note that $\mathbf{U}$, $\mathbf{U}^*$, and $\mathbf{U}^{-1}$ are all both **isometric** and **normal**:

$$\begin{array}{rcccll} \mathbf{U}^*\mathbf{U} & & & = & \mathbf{I} & \Longrightarrow \quad \mathbf{U} \text{ is isometric} \\ \mathbf{UU}^* & = & \mathbf{U}^*\mathbf{U} & = & \mathbf{I} & \Longrightarrow \quad \mathbf{U}^* \text{ is isometric} \\ \mathbf{U}^{-1} & = & \mathbf{U}^* & & & \Longrightarrow \quad \mathbf{U}^{-1} \text{ is isometric} \end{array}$$

$$\begin{array}{rcccl} \mathbf{U}^*\mathbf{U} & = & \mathbf{UU}^* & = & \Longrightarrow \quad \mathbf{U} \text{ is normal} \\ \mathbf{UU}^* & = & \mathbf{U}^*\mathbf{U} & = & \Longrightarrow \quad \mathbf{U}^* \text{ is normal} \\ \mathbf{U}^{-1} & = & \mathbf{U}^* & & \Longrightarrow \quad \mathbf{U}^{-1} \text{ is normal} \end{array}$$

(2)   Proof that $\mathbf{U}^*\mathbf{U} = \mathbf{UU}^* = \mathbf{I} \Longrightarrow \mathcal{I}(\mathbf{U}) = \mathcal{I}(\mathbf{U}^*) = H$: by Theorem 1.50 page 14.

(3)   Proof that $\mathbf{U}^*\mathbf{U} = \mathbf{UU}^* = \mathbf{I} \Longrightarrow \mathcal{N}(\mathbf{U}) = \mathcal{N}(\mathbf{U}^*) = \mathcal{N}(\mathbf{U}^{-1})$:

$$\begin{array}{ll} \mathcal{N}(\mathbf{U}^*) = \mathcal{N}(\mathbf{U}) & \text{because } \mathbf{U} \text{ and } \mathbf{U}^* \text{ are both normal and by Theorem 1.41 page 12} \\ = \mathcal{I}(\mathbf{U})^\perp & \text{by Theorem 1.36 page 11} \\ = X^\perp & \text{by above result} \\ = \{0\} & \end{array}$$

---

(4) Proof that $\mathbf{U}^*\mathbf{U} = \mathbf{U}\mathbf{U}^* = \mathbf{I} \implies \left\|\mathbf{U}^{-1}\right\| = \|\mathbf{U}^*\| = \|\mathbf{U}\| = 1$:

Because $\mathbf{U}$, $\mathbf{U}^*$, and $\mathbf{U}^{-1}$ are all isometric and by Theorem 1.46 page 13.

# 2   Background: harmonic analysis

## 2.1   Families of functions

This paper is largely set in the space of *Lebesgue square-integrable functions* $L_{\mathbb{R}}^2$ (Definition 2.2 page 15). The space $L_{\mathbb{R}}^2$ is a subspace of the space $\mathbb{R}^{\mathbb{R}}$, the set of all functions with *domain* $\mathbb{R}$ (the set of real numbers) and *range* $\mathbb{R}$. The space $\mathbb{R}^{\mathbb{R}}$ is a subspace of the space $\mathbb{C}^{\mathbb{C}}$, the set of all functions with *domain* $\mathbb{C}$ (the set of complex numbers) and *range* $\mathbb{C}$. That is, $L_{\mathbb{R}}^2 \subseteq \mathbb{R}^{\mathbb{R}} \subseteq \mathbb{C}^{\mathbb{C}}$. In general, the notation $Y^X$ represents the set of all functions with domain $X$ and range $Y$ (Definition 2.1 page 15). Although this notation may seem curious, note that for finite $X$ and finite $Y$, the number of functions (elements) in $Y^X$ is $\left|Y^X\right| = |Y|^{|X|}$.

**Definition 2.1**   Let $X$ and $Y$ be sets.
The space $Y^X$ represents the set of all functions with *domain $X$* and *range $Y$* such that
$$Y^X \triangleq \{\mathsf{f}(x) \,|\, \mathsf{f}(x) : X \to Y\}$$

**Definition 2.2**   Let $\mathbb{R}$ be the set of real numbers, $\mathscr{B}$ the set of *Borel sets* on $\mathbb{R}$, and $\mu$ the standard *Borel measure* on $\mathscr{B}$. Let $\mathbb{R}^{\mathbb{R}}$ be as in Definition 2.1 page 15.
The space of **Lebesgue square-integrable functions** $L_{(\mathbb{R}, \mathscr{B}, \mu)}^2$ (or $L_{\mathbb{R}}^2$) is defined as
$$L_{\mathbb{R}}^2 \triangleq L_{(\mathbb{R}, \mathscr{B}, \mu)}^2 \triangleq \left\{ \mathsf{f} \in \mathbb{R}^{\mathbb{R}} \,\left|\, \left(\int_{\mathbb{R}} |\mathsf{f}|^2\right)^{\frac{1}{2}} \mathrm{d}\mu < \infty \right. \right\}.$$
The **standard inner product** $\langle \triangle \,|\, \triangledown \rangle$ on $L_{\mathbb{R}}^2$ is defined as
$$\langle \mathsf{f}(x) \,|\, \mathsf{g}(x) \rangle \triangleq \int_{\mathbb{R}} \mathsf{f}(x)\mathsf{g}^*(x) \, \mathrm{d}x.$$
The **standard norm** $\|\cdot\|$ on $L_{\mathbb{R}}^2$ is defined as $\|\mathsf{f}(x)\| \triangleq \langle \mathsf{f}(x) \,|\, \mathsf{f}(x) \rangle^{\frac{1}{2}}$

**Definition 2.3** [50]   Let $X$ be a set.
The **indicator function** $\mathbb{1} \in \{0,1\}^{2^X}$ is defined as
$$\mathbb{1}_A(x) = \begin{cases} 1 & \text{for } x \in A & \forall_{x \in X, \ A \in 2^X} \\ 0 & \text{for } x \notin A & \forall_{x \in X, \ A \in 2^X} \end{cases}$$
The indicator function $\mathbb{1}$ is also called the **characteristic function**.

---

## 2.2 Trigonometric functions

### 2.2.1 Definitions

**Lemma 2.4** [51] *Let* $\mathcal{C}$ *be the* SPACE OF ALL CONTINUOUSLY DIFFERENTIABLE REAL FUNCTIONS *and* $\frac{d}{dx} \in \mathcal{C}^{\mathcal{C}}$ *the differentiation operator.* $\frac{d}{dx}^2 f + f = 0 \iff$

$$
\left\{
\begin{aligned}
f(x) &= \underbrace{[f](0) \sum_{n=0}^{\infty} (-1)^n \frac{x^{2n}}{(2n)!}}_{\textit{even terms}} + \underbrace{\left[\frac{d}{dx} f\right](0) \sum_{n=0}^{\infty} (-1)^n \frac{x^{2n+1}}{(2n+1)!}}_{\textit{odd terms}} \\
&= \left( f(0) + \left[\frac{d}{dx} f\right](0) x \right) - \left( \frac{f(0)}{2!} x^2 + \frac{\left[\frac{d}{dx} f\right](0)}{3!} x^3 \right) + \left( \frac{f(0)}{4!} x^4 + \frac{\left[\frac{d}{dx} f\right](0)}{5!} x^5 \right) \cdots
\end{aligned}
\right\}
$$

**Definition 2.5** [52] Let $\mathcal{C}$ be the *space of all continuously differentiable real functions* and $\frac{d}{dx} \in \mathcal{C}^{\mathcal{C}}$ the differentiation operator.

The **cosine** function $\cos(x)$ is the function $f \in \mathcal{C}$ that satisfies the following conditions:

$$\underbrace{\frac{d}{dx}^2 f + f = 0}_{\text{2nd order homogeneous differential equation}} \qquad \underbrace{f(0) = 1}_{\text{1st initial condition}} \qquad \underbrace{\left[\frac{d}{dx} f\right](0) = 0}_{\text{2nd initial condition}}$$

The **sine** function $\sin(x)$ is the function $g \in \mathcal{C}$ that satisfies the following conditions:

$$\underbrace{\frac{d}{dx}^2 g + g = 0}_{\text{2nd order homogeneous differential equation}} \qquad \underbrace{g(0) = 0}_{\text{1st initial condition}} \qquad \underbrace{\left[\frac{d}{dx} g\right](0) = 1}_{\text{2nd initial condition}}$$

**Theorem 2.6** [53]

$$
\begin{aligned}
\cos(x) &= \sum_{n=0}^{\infty} (-1)^n \frac{x^{2n}}{(2n)!} &= 1 - \frac{x^2}{2} + \frac{x^4}{4!} - \frac{x^6}{6!} + \cdots & \qquad \forall x \in \mathbb{R} \\
\sin(x) &= \sum_{n=0}^{\infty} (-1)^n \frac{x^{2n+1}}{(2n+1)!} &= x - \frac{x^3}{3!} + \frac{x^5}{5!} - \frac{x^7}{7!} + \cdots & \qquad \forall x \in \mathbb{R}
\end{aligned}
$$

**Proposition 2.7** [54] *Let* $\mathcal{C}$ *be the space of all continuously differentiable real functions and*

---

[51] ☞ [110], page 156, ☞ [92]

[52] ☞ [110], page 157, ☞ [36], pages 228–229

[53] ☞ [110], page 157

[54] ☞ [110], page 157. The general solution for the *non-homogeneous* equation $\frac{d}{dx}^2 f(x) + f(x) = g(x)$ with initial conditions $f(a) = 1$ and $f'(a) = \rho$ is

$$f(x) = \cos(x) + \rho \sin(x) + \int_a^x g(y) \sin(x - y) \, dy.$$

This type of equation is called a *Volterra integral equation of the second type.* References: ☞ [37], page 371, ☞ [92]. Volterra equation references: ☞ [105], page 99, ☞ [87], ☞ [88].

$\frac{d}{dx} \in C^C$ the differentiation operator. Let $f'(0) \triangleq \left[\frac{d}{dx}f\right](0)$.

$$\underbrace{\frac{d^2}{dx}f + f = 0}_{\text{2ND ORDER HOMOGENEOUS DIFFERENTIAL EQUATION}} \iff f(x) = f(0)\cos(x) + f'(0)\sin(x) \quad \forall f \in C, \forall x \in \mathbb{R}$$

**Theorem 2.8** [55] *Let* $\frac{d}{dx} \in C^C$ *be the differentiation operator.*

$$\begin{aligned}\frac{d}{dx}\cos(x) &= -\sin(x) &\forall x \in \mathbb{R} \\ \frac{d}{dx}\sin(x) &= \cos(x) &\forall x \in \mathbb{R}\end{aligned}$$

### 2.2.2 The complex exponential

**Definition 2.9** The function $f \in \mathbb{C}^C$ is the **exponential function** $\exp(ix) \triangleq f(x)$ if

1. $\frac{d^2}{dx}f + f = 0$   (second order homogeneous differential equation)   and
2. $f(0) = 1$   (first initial condition)   and
3. $\left[\frac{d}{dx}f\right](0) = i$   (second initial condition)

**Theorem 2.10** (Euler's identity) [56]

$e^{ix} = \cos(x) + i\sin(x) \quad \forall x \in \mathbb{R}$

**Corollary 2.11**

$e^{ix} = \displaystyle\sum_{n \in \mathbb{W}} \frac{(ix)^n}{n!} \quad \forall x \in \mathbb{R}$

**Corollary 2.12** [57]

$e^{i\pi} + 1 = 0$

The exponential has two properties that makes it extremely special:

☞ The exponential is an eigenvalue of any LTI operator (Theorem 2.13 page 17).

☞ The exponential generates a continuous point spectrum for the differential operator.

**Theorem 2.13** [58] *Let* $\mathbf{L}$ *be an operator with kernel* $h(t, \omega)$ *and*

$$\hat{h}(s) \triangleq \left\langle h(t, \omega) \mid e^{st} \right\rangle \quad \text{(LAPLACE TRANSFORM)}.$$

$\left.\begin{array}{l} 1. \ \mathbf{L} \ \textit{is linear and} \\ 2. \ \mathbf{L} \ \textit{is time-invariant} \end{array}\right\} \implies \mathbf{L}e^{st} = \underbrace{\hat{h}^*(-s)}_{\textit{eigenvalue}} \underbrace{e^{st}}_{\textit{eigenvector}}$

---

### 2.2.3 Trigonometric Identities

**Corollary 2.14** (Euler formulas) [59]

$$\cos(x) = \Re\left(e^{ix}\right) = \frac{e^{ix} + e^{-ix}}{2} \qquad \forall x \in \mathbb{R}$$

$$\sin(x) = \Im\left(e^{ix}\right) = \frac{e^{ix} - e^{-ix}}{2i} \qquad \forall x \in \mathbb{R}$$

**Theorem 2.15** [60]

$$e^{(\alpha+\beta)} = e^{\alpha}\,e^{\beta} \qquad \forall \alpha, \beta \in \mathbb{C}$$

## 2.3 Fourier Series

The *Fourier Series* expansion of a periodic function is simply a complex trigonometric polynomial. In the special case that the periodic function is even, then the Fourier Series expansion is a cosine polynomial.

**Definition 2.16** [61] The **Fourier Series operator** $\hat{\mathbf{F}} : L_{\mathbb{R}}^2 \to \ell_{\mathbb{R}}^2$ is defined as

$$\left[\hat{\mathbf{F}}f\right](n) \triangleq \frac{1}{\sqrt{\tau}} \int_0^\tau f(x) e^{-i\frac{2\pi}{\tau}nx}\, dx \qquad \forall f \in \left\{ f \in L_{\mathbb{R}}^2 \,\middle|\, f \text{ is periodic with period } \tau \right\}$$

**Theorem 2.17** *Let $\hat{\mathbf{F}}$ be the Fourier Series operator.*
*The **inverse Fourier Series** operator $\hat{\mathbf{F}}^{-1}$ is given by*

$$\left[\hat{\mathbf{F}}^{-1}\left(\tilde{x}_n\right)_{n\in\mathbb{Z}}\right](x) \triangleq \frac{1}{\sqrt{\tau}} \sum_{n\in\mathbb{Z}} \tilde{x}_n e^{i\frac{2\pi}{\tau}nx} \qquad \forall (\tilde{x}_n) \in \ell_{\mathbb{R}}^2$$

**Theorem 2.18** *The **Fourier Series adjoint** operator $\hat{\mathbf{F}}^*$ is given by*

$$\hat{\mathbf{F}}^* = \hat{\mathbf{F}}^{-1}$$

---

[59]  [33],  [18], page 12
[60]  [111], page 1
[61]  [78], page 3

✎ PROOF:

$$
\begin{aligned}
\left\langle \hat{\mathbf{F}}\mathsf{x}(x) \;\middle|\; \bar{\mathsf{y}}(n) \right\rangle_{\mathbb{Z}} &= \left\langle \frac{1}{\sqrt{\tau}}\int_0^\tau \mathsf{x}(x)e^{-i\frac{2\pi}{\tau}nx}\,\mathrm{d}x \;\middle|\; \bar{\mathsf{y}}(n) \right\rangle_{\mathbb{Z}} && \text{by definition of } \hat{\mathbf{F}} \text{ Definition 2.16 page 18} \\
&= \frac{1}{\sqrt{\tau}}\int_0^\tau \mathsf{x}(x)\left\langle e^{-i\frac{2\pi}{\tau}nx} \;\middle|\; \bar{\mathsf{y}}(n) \right\rangle_{\mathbb{Z}}\,\mathrm{d}x && \text{by additivity property of } \langle \triangle \mid \triangledown \rangle \\
&= \int_0^\tau \mathsf{x}(x)\frac{1}{\sqrt{\tau}}\left\langle \bar{\mathsf{y}}(n) \;\middle|\; e^{-i\frac{2\pi}{\tau}nx} \right\rangle_{\mathbb{Z}}^{*}\,\mathrm{d}x && \text{by property of } \langle \triangle \mid \triangledown \rangle \\
&= \int_0^\tau \mathsf{x}(x)\left[\hat{\mathbf{F}}^{-1}\bar{\mathsf{y}}(n)\right]^{*}\,\mathrm{d}x && \text{by definition of } \hat{\mathbf{F}}^{-1} \text{ page 18} \\
&= \left\langle \mathsf{x}(x) \;\middle|\; \underbrace{\hat{\mathbf{F}}^{-1}}_{\hat{\mathbf{F}}^{*}}\bar{\mathsf{y}}(n) \right\rangle_{\mathbb{R}}
\end{aligned}
$$

✐

The Fourier Series operator has several nice properties:

- ✎ $\hat{\mathbf{F}}$ is *unitary* [62] (Corollary 2.19 page 19).
- ✎ Because $\hat{\mathbf{F}}$ is unitary, it automatically has several other nice properties such as being *isometric,* and satisfying *Parseval's equation,* satisfying *Plancheral's formula,* and more (Corollary 2.20 page 19).

**Corollary 2.19** *Let* $\mathbf{I}$ *be the identity operator and let* $\hat{\mathbf{F}}$ *be the Fourier Series operator with adjoint* $\hat{\mathbf{F}}^{*}$.

$$\hat{\mathbf{F}}\hat{\mathbf{F}}^{*} = \hat{\mathbf{F}}^{*}\hat{\mathbf{F}} = \mathbf{I} \qquad \textit{($\hat{F}$ is unitary…and thus also normal and isometric)}$$

✎ PROOF:    This follows directly from the fact that $\hat{\mathbf{F}}^{*} = \hat{\mathbf{F}}^{-1}$ (Theorem 2.18 (page 18)).    ✐

**Corollary 2.20** *Let* $\hat{\mathbf{F}}$ *be the Fourier series operator,* $\hat{\mathbf{F}}^{*}$ *be its adjoint, and* $\hat{\mathbf{F}}^{-1}$ *be its inverse.*

$$
\begin{aligned}
\mathcal{R}(\hat{\mathbf{F}}) &= \mathcal{R}(\hat{\mathbf{F}}^{-1}) &&= \boldsymbol{L}_{\mathbb{R}}^2 \\
\left\|\left|\hat{\mathbf{F}}\right|\right\| &= \left\|\left|\hat{\mathbf{F}}^{-1}\right|\right\| &&= 1 && \textsc{(unitary)} \\
\left\langle \hat{\mathbf{F}}\boldsymbol{x} \;\middle|\; \hat{\mathbf{F}}\boldsymbol{y} \right\rangle &= \left\langle \hat{\mathbf{F}}^{-1}\boldsymbol{x} \;\middle|\; \hat{\mathbf{F}}^{-1}\boldsymbol{y} \right\rangle &&= \left\langle \boldsymbol{x} \;\middle|\; \boldsymbol{y} \right\rangle && \textsc{(Parseval's equation)} \\
\left\|\hat{\mathbf{F}}\boldsymbol{x}\right\| &= \left\|\hat{\mathbf{F}}^{-1}\boldsymbol{x}\right\| &&= \|\boldsymbol{x}\| && \textsc{(Plancherel's formula)} \\
\left\|\hat{\mathbf{F}}\boldsymbol{x} - \hat{\mathbf{F}}\boldsymbol{y}\right\| &= \left\|\hat{\mathbf{F}}^{-1}\boldsymbol{x} - \hat{\mathbf{F}}^{-1}\boldsymbol{y}\right\| &&= \|\boldsymbol{x} - \boldsymbol{y}\| && \textsc{(isometric)}
\end{aligned}
$$

✎ PROOF:    These results follow directly from the fact that $\hat{\mathbf{F}}$ is unitary (Corollary 2.19 page 19) and from the properties of unitary operators (Theorem 1.51 page 14).    ✐

**Theorem 2.21**    *The set*
$$\left\{ \frac{1}{\sqrt{\tau}}e^{i\frac{2\pi}{\tau}nx} \;\middle|_{n\in\mathbb{Z}} \right\}$$
*is an* ORTHONORMAL BASIS *for all functions* $\mathsf{f}(x)$ *with support in* $[0,\ \tau]$.

---

[62] *unitary operators*: Definition 1.48 page 13

## 2.4   Fourier Transform

### 2.4.1   Properties

**Definition 2.22** [63] The **Fourier Transform** operator $\tilde{\mathbf{F}}$ is defined as[64]

$$\left[\tilde{\mathbf{F}}\mathsf{f}\right](\omega) \triangleq \frac{1}{\sqrt{2\pi}} \int_{\mathbb{R}} \mathsf{f}(x)\,e^{-i\omega x}\;\mathrm{d}x \qquad \forall \mathsf{f} \in L^2_{(\mathbb{R},\mathscr{B},\mu)}$$

This definition of the Fourier Transform is also called the **unitary Fourier Transform**.

**Remark 2.23** (**Fourier transform scaling factor**) [65] If the Fourier transform operator $\tilde{\mathbf{F}}$ and inverse Fourier transform operator $\tilde{\mathbf{F}}^{-1}$ are defined as

$$\tilde{\mathbf{F}}\mathsf{f}(x) \;\triangleq\; A\int_{\mathbb{R}} \mathsf{f}(x)e^{-i\omega x}\,\mathrm{d}x \quad \text{and} \quad \tilde{\mathbf{F}}^{-1}\tilde{\mathsf{f}}(\omega) \;\triangleq\; B\int_{\mathbb{R}} \mathsf{F}(\omega)e^{i\omega x}\,\mathrm{d}\omega \quad,$$

then $A$ and $B$ can be any constants as long as $AB = \frac{1}{2\pi}$. The Fourier transform is often defined with the scaling factor $A$ set equal to 1 such that $\left[\tilde{\mathbf{F}}\mathsf{f}(x)\right](\omega) = \int_{\mathbb{R}} \mathsf{f}(x)\,e^{-i\omega x}\,\mathrm{d}x$. In this case, the inverse Fourier transform operator $\tilde{\mathbf{F}}^{-1}$ is either defined as

✎   $\left[\tilde{\mathbf{F}}^{-1}\mathsf{f}(x)\right](f) \triangleq \int_{\mathbb{R}} \mathsf{f}(x)e^{i2\pi fx}\,\mathrm{d}x$   (using oscillatory frequency free variable $f$) or

✎   $\left[\tilde{\mathbf{F}}^{-1}\mathsf{f}(x)\right](\omega) \triangleq \frac{1}{2\pi}\int_{\mathbb{R}} \mathsf{f}(x)e^{i\omega x}\,\mathrm{d}x$   (using angular frequency free variable $\omega$).

In short, the $2\pi$ has to show up somewhere, either in the argument of the exponential ($e^{-i2\pi ft}$) or in front of the integral ($\frac{1}{2\pi}\int \cdots$). One could argue that it is unnecessary to burden the exponential argument with the $2\pi$ factor ($e^{-i2\pi ft}$), and thus could further argue in favor of using the angular frequency variable $\omega$ thus giving the inverse operator definition $\left[\tilde{\mathbf{F}}^{-1}\mathsf{f}(x)\right](\omega) \triangleq \frac{1}{2\pi}\int_{\mathbb{R}} \mathsf{f}(x)\,e^{-i\omega x}\,\mathrm{d}x$. But this causes a new problem. In this case, the Fourier operator $\tilde{\mathbf{F}}$ is not *unitary* (see Theorem 2.25 page 20)—in particular, $\tilde{\mathbf{F}}\tilde{\mathbf{F}}^* \neq \mathbf{I}$, where $\tilde{\mathbf{F}}^*$ is the *adjoint* of $\tilde{\mathbf{F}}$; but rather, $\tilde{\mathbf{F}}\left(\frac{1}{2\pi}\tilde{\mathbf{F}}^*\right) = \left(\frac{1}{2\pi}\tilde{\mathbf{F}}^*\right)\tilde{\mathbf{F}} = \mathbf{I}$. But if we define the operators $\tilde{\mathbf{F}}$ and $\tilde{\mathbf{F}}^{-1}$ to both have the scaling factor $\frac{1}{\sqrt{2\pi}}$, then $\tilde{\mathbf{F}}$ and $\tilde{\mathbf{F}}^{-1}$ are inverses *and* $\tilde{\mathbf{F}}$ is *unitary*—that is, $\tilde{\mathbf{F}}\tilde{\mathbf{F}}^* = \tilde{\mathbf{F}}^*\tilde{\mathbf{F}} = \mathbf{I}$.

**Theorem 2.24** (Inverse Fourier transform) [66] *Let $\tilde{\mathbf{F}}$ be the Fourier Transform operator* (Definition 2.22 page 20). *The inverse $\tilde{\mathbf{F}}^{-1}$ of $\tilde{\mathbf{F}}$ is*

$$\left[\tilde{\mathbf{F}}^{-1}\tilde{\mathsf{f}}\right](x) \triangleq \frac{1}{\sqrt{2\pi}} \int_{\mathbb{R}} \tilde{\mathsf{f}}(\omega)e^{i\omega x}\,\mathrm{d}\omega \qquad \forall \tilde{\mathsf{f}} \in L^2_{(\mathbb{R},\mathscr{B},\mu)}$$

**Theorem 2.25** *Let $\tilde{\mathbf{F}}$ be the Fourier Transform operator with inverse $\tilde{\mathbf{F}}^{-1}$ and adjoint $\tilde{\mathbf{F}}^*$.*

$$\tilde{\mathbf{F}}^* = \tilde{\mathbf{F}}^{-1}$$

---

✎ PROOF:

$$\left\langle \tilde{\mathbf{F}}f \mid g \right\rangle = \left\langle \frac{1}{\sqrt{2\pi}} \int_{\mathbb{R}} f(x)\, e^{-i\omega x}\; \mathrm{d}x \;\middle|\; g(\omega) \right\rangle \qquad \text{by definition of } \tilde{\mathbf{F}} \text{ page } 20$$

$$= \frac{1}{\sqrt{2\pi}} \int_{\mathbb{R}} f(x) \left\langle e^{-i\omega x} \mid g(\omega) \right\rangle\; \mathrm{d}x \qquad \text{by } \textit{additive property} \text{ of } \langle \triangle \mid \triangledown \rangle$$

$$= \int_{\mathbb{R}} f(x) \frac{1}{\sqrt{2\pi}} \left\langle g(\omega) \mid e^{-i\omega x} \right\rangle^{*}\; \mathrm{d}x \qquad \text{by } \textit{conjugate symmetric property} \text{ of } \langle \triangle \mid \triangledown \rangle$$

$$= \left\langle f(x) \mid \frac{1}{\sqrt{2\pi}} \left\langle g(\omega) \mid e^{-i\omega x} \right\rangle \right\rangle \qquad \text{by definition of } \langle \triangle \mid \triangledown \rangle$$

$$= \left\langle f \mid \underbrace{\tilde{\mathbf{F}}^{-1}}_{\tilde{\mathbf{F}}^{*}}g \right\rangle \qquad \text{by Theorem } 2.24 \text{ page } 20$$

✐

The Fourier Transform operator has several nice properties:

- ✐ $\tilde{\mathbf{F}}$ is *unitary* (Definition 1.48 page 13) (Corollary 2.26—next corollary).
- ✐ Because $\tilde{\mathbf{F}}$ is unitary, it automatically has several other nice properties (Theorem 2.27 page 21).

**Corollary 2.26**   *Let $\mathbf{I}$ be the identity operator and let $\tilde{\mathbf{F}}$ be the Fourier Transform operator with adjoint $\tilde{\mathbf{F}}^{*}$ and inverse $\tilde{\mathbf{F}}^{-1}$.*

$$\underbrace{\tilde{\mathbf{F}}\tilde{\mathbf{F}}^{*} = \tilde{\mathbf{F}}^{*}\tilde{\mathbf{F}} = \mathbf{I}}_{\tilde{\mathbf{F}}^{*}=\tilde{\mathbf{F}}^{-1}} \qquad (\tilde{\mathbf{F}} \text{ is unitary})$$

✎ PROOF:    This follows directly from the fact that $\tilde{\mathbf{F}}^{*} = \tilde{\mathbf{F}}^{-1}$ (Theorem 2.25 page 20).    ✐

**Theorem 2.27**   *Let $\tilde{\mathbf{F}}$ be the Fourier transform operator with adjoint $\tilde{\mathbf{F}}^{*}$ and inverse $\tilde{\mathbf{F}}$. Let $\|\|\cdot\|\|$ be the operator norm with respect to the vector norm $\|\cdot\|$ with respect to the Hilbert space $\left(\mathbb{C}^{\mathbb{R}}, \langle \triangle \mid \triangledown \rangle\right)$. Let $\mathcal{R}(\mathbf{A})$ be the range of an operator $\mathbf{A}$.*

$$\mathcal{R}(\mathbf{F}\tau) = \mathcal{R}(\tilde{\mathbf{F}}^{-1}) = L_{\mathbb{R}}^{2}$$

$$\|\|\tilde{\mathbf{F}}\|\| = \|\|\tilde{\mathbf{F}}^{-1}\|\| = 1 \qquad \text{(UNITARY)}$$

$$\left\langle \tilde{\mathbf{F}}f \mid \tilde{\mathbf{F}}g \right\rangle = \left\langle \tilde{\mathbf{F}}^{-1}f \mid \tilde{\mathbf{F}}^{-1}g \right\rangle = \langle f \mid g \rangle \qquad \text{(PARSEVAL'S EQUATION)}$$

$$\|\tilde{\mathbf{F}}f\| = \|\tilde{\mathbf{F}}^{-1}f\| = \|f\| \qquad \text{(PLANCHEREL'S FORMULA)}$$

$$\|\tilde{\mathbf{F}}f - \tilde{\mathbf{F}}g\| = \|\tilde{\mathbf{F}}^{-1}f - \tilde{\mathbf{F}}^{-1}g\| = \|f - g\| \qquad \text{(ISOMETRIC)}$$

✎ PROOF:    These results follow directly from the fact that $\tilde{\mathbf{F}}$ is unitary (Corollary 2.26 page 21) and from the properties of unitary operators (Theorem 1.51 page 14).    ✐

**Theorem 2.28**  (Shift relations)  *Let $\tilde{\mathbf{F}}$ be the Fourier transform operator.*

$$\tilde{\mathbf{F}}[\mathsf{f}(x-u)](\omega) = e^{-i\omega u}\left[\tilde{\mathbf{F}}\mathsf{f}(x)\right](\omega)$$
$$\left[\tilde{\mathbf{F}}\left(e^{ivx}\mathsf{g}(x)\right)\right](\omega) = \left[\tilde{\mathbf{F}}\mathsf{g}(x)\right](\omega-v)$$

**Theorem 2.29**  (Complex conjugate)  *Let $\tilde{\mathbf{F}}$ be the Fourier Transform operator and $*$ represent the complex conjugate operation on the set of complex numbers.*

$$\tilde{\mathbf{F}}\mathsf{f}^*(-x) = \left[\tilde{\mathbf{F}}\mathsf{f}(x)\right]^* \qquad \forall \mathsf{f}\in L^2_{(\mathbb{R},\mathscr{B},\mu)}$$

### 2.4.2  Convolution

**Definition 2.30** [67]

The **convolution operation** is defined as

$$\left[\mathsf{f}\star\mathsf{g}\right](x) \triangleq \mathsf{f}(x)\star\mathsf{g}(x) \triangleq \int_{u\in\mathbb{R}}\mathsf{f}(u)\mathsf{g}(x-u)\,\mathrm{d}u \qquad \forall \mathsf{f},\mathsf{g}\in L^2_{(\mathbb{R},\mathscr{B},\mu)}$$

Theorem 2.31 (next) demonstrates that multiplication in the "time domain" is equivalent to convolution in the "frequency domain" and vice-versa.

**Theorem 2.31**  (convolution theorem)  *Let $\tilde{\mathbf{F}}$ be the Fourier Transform operator and $\star$ the convolution operator.*

$$\underbrace{\tilde{\mathbf{F}}\left[\mathsf{f}(x)\star\mathsf{g}(x)\right](\omega)}_{\text{convolution in "time domain"}} = \underbrace{\sqrt{2\pi}\left[\tilde{\mathbf{F}}\mathsf{f}\right](\omega)\left[\tilde{\mathbf{F}}\mathsf{g}\right](\omega)}_{\text{multiplication in "frequency domain"}} \qquad \forall \mathsf{f},\mathsf{g}\in L^2_{(\mathbb{R},\mathscr{B},\mu)}$$

$$\underbrace{\tilde{\mathbf{F}}\left[\mathsf{f}(x)\mathsf{g}(x)\right](\omega)}_{\text{multiplication in "time domain"}} = \underbrace{\frac{1}{\sqrt{2\pi}}\left[\tilde{\mathbf{F}}\mathsf{f}\right](\omega)\star\left[\tilde{\mathbf{F}}\mathsf{g}\right](\omega)}_{\text{convolution in "frequency domain"}} \qquad \forall \mathsf{f},\mathsf{g}\in L^2_{(\mathbb{R},\mathscr{B},\mu)}.$$

---

[67] ☞ [9], page 6

✎Proof:

$$\tilde{\mathbf{F}}\Big[\mathsf{f}(x) \star \mathsf{g}(x)\Big](\omega) = \tilde{\mathbf{F}}\left[\int_{u\in\mathbb{R}} \mathsf{f}(u)\mathsf{g}(x-u)\,\mathrm{d}u\right](\omega) \qquad\qquad \text{by def. of } \star \text{ (Definition 2.30 page 22)}$$

$$= \int_{u\in\mathbb{R}} \mathsf{f}(u)\Big[\tilde{\mathbf{F}}\mathsf{g}(x-u)\Big](\omega)\,\mathrm{d}u$$

$$= \int_{u\in\mathbb{R}} \mathsf{f}(u)e^{-i\omega u}\Big[\tilde{\mathbf{F}}\mathsf{g}(x)\Big](\omega)\,\mathrm{d}u \qquad\qquad \text{by Theorem 2.28 page 22}$$

$$= \sqrt{2\pi}\underbrace{\left(\frac{1}{\sqrt{2\pi}}\int_{u\in\mathbb{R}} \mathsf{f}(u)e^{-i\omega u}\,\mathrm{d}u\right)}_{[\tilde{\mathbf{F}}\mathsf{f}](\omega)}\Big[\tilde{\mathbf{F}}\mathsf{g}\Big](\omega)$$

$$= \sqrt{2\pi}\Big[\tilde{\mathbf{F}}\mathsf{f}\Big](\omega)\,\Big[\tilde{\mathbf{F}}\mathsf{g}\Big](\omega) \qquad\qquad \text{by def. of } \tilde{\mathbf{F}} \text{ (Definition 2.22 page 20)}$$

$$\tilde{\mathbf{F}}[\mathsf{f}(x)\mathsf{g}(x)](\omega) = \tilde{\mathbf{F}}\Big[\big(\tilde{\mathbf{F}}^{-1}\tilde{\mathbf{F}}\mathsf{f}(x)\big)\,\mathsf{g}(x)\Big](\omega) \qquad\qquad \text{by def. of operator inverse (page 6)}$$

$$= \tilde{\mathbf{F}}\left[\left(\tfrac{1}{\sqrt{2\pi}}\int_{v\in\mathbb{R}} \big[\tilde{\mathbf{F}}\mathsf{f}(x)\big](v)e^{ivx}\,\mathrm{d}v\right)\mathsf{g}(x)\right](\omega) \qquad \text{by Theorem 2.24 page 20}$$

$$= \tfrac{1}{\sqrt{2\pi}}\int_{v\in\mathbb{R}} \big[\tilde{\mathbf{F}}\mathsf{f}(x)\big](v)\big[\tilde{\mathbf{F}}\big(e^{ivx}\,\mathsf{g}(x)\big)\big](\omega,v)\,\mathrm{d}v$$

$$= \tfrac{1}{\sqrt{2\pi}}\int_{v\in\mathbb{R}} \big[\tilde{\mathbf{F}}\mathsf{f}(x)\big](v)\big[\tilde{\mathbf{F}}\mathsf{g}(x)\big](\omega-v)\,\mathrm{d}v \qquad \text{by Theorem 2.28 page 22}$$

$$= \tfrac{1}{\sqrt{2\pi}}\big[\tilde{\mathbf{F}}\mathsf{f}\big](\omega) \star \big[\tilde{\mathbf{F}}\mathsf{g}\big](\omega) \qquad\qquad \text{by def. of } \star \text{ (Definition 2.30 page 22)}$$

✐

## 2.5 Operations on sequences

### 2.5.1 Z-transform

**Definition 2.32** [68] Let $X^Y$ be the set of all functions from a set $Y$ to a set $X$. Let $\mathbb{Z}$ be the set of integers.

A function $\mathsf{f}$ in $X^Y$ is a **sequence** over $X$ if $Y = \mathbb{Z}$.

A sequence may be denoted in the form $(x_n)_{n\in\mathbb{Z}}$ or simply as $(x_n)$.

**Definition 2.33** [69] Let $(\mathbb{F}, +, \cdot)$ be a field.

The **space of all absolutely square summable sequences** $\ell_{\mathbb{F}}^2$ over $\mathbb{F}$ is defined as

$$\ell_{\mathbb{F}}^2 \triangleq \left\{ (x_n)_{n\in\mathbb{Z}} \,\middle|\, \sum_{n\in\mathbb{Z}} |x_n|^2 < \infty \right\}$$

---

[68] 📖 [19], page 1, 📖 [129], page 23, ⟨Definition 2.1⟩, 📖 [75], page 31
[69] 📖 [85], page 347, ⟨Example 5.K⟩

The space $\ell_{\mathbb{R}}^2$ is an example of a *separable Hilbert space*. In fact, $\ell_{\mathbb{R}}^2$ is the *only* separable Hilbert space in the sense that all separable Hilbert spaces are isomorphically equivalent. For example, $\ell_{\mathbb{R}}^2$ is isomorphic to $L_{\mathbb{R}}^2$, the *space of all absolutely square Lebesgue integrable functions*.

**Definition 2.34**  The **convolution** operation $\ast$ is defined as

$$(x_n) \ast (y_n) \triangleq \left( \sum_{m \in \mathbb{Z}} x_m y_{n-m} \right)_{n \in \mathbb{Z}} \qquad \forall (x_n)_{n \in \mathbb{Z}}, (y_n)_{n \in \mathbb{Z}} \in \ell_{\mathbb{R}}^2$$

**Definition 2.35**  [70] The **z-transform** $\mathbf{Z}$ of $(x_n)_{n \in \mathbb{Z}}$ is defined as

$$[\mathbf{Z}(x_n)](z) \triangleq \underbrace{\sum_{n \in \mathbb{Z}} x_n z^{-n}}_{\text{Laurent series}} \qquad \forall (x_n) \in \ell_{\mathbb{R}}^2$$

**Proposition 2.36**  [71] *Let $\ast$ be the* CONVOLUTION OPERATOR *(Definition 2.34 page 24).*

$$(x_n) \ast (y_n) = (y_n) \ast (x_n) \qquad \forall (x_n)_{n \in \mathbb{Z}}, (y_n)_{n \in \mathbb{Z}} \in \ell_{\mathbb{R}}^2 \qquad (\ast \text{ is } \text{COMMUTATIVE})$$

**Theorem 2.37**  [72] *Let $\ast$ be the convolution operator (Definition 2.34 page 24).*

$$\mathbf{Z} \underbrace{\left( (x_n) \ast (y_n) \right)}_{\text{sequence convolution}} = \underbrace{\left( \mathbf{Z}(x_n) \right) \left( \mathbf{Z}(y_n) \right)}_{\text{series multiplication}} \qquad \forall (x_n)_{n \in \mathbb{Z}}, (y_n)_{n \in \mathbb{Z}} \in \ell_{\mathbb{R}}^2$$

### 2.5.2  Discrete Time Fourier Transform

**Definition 2.38**  The **discrete-time Fourier transform** $\check{\mathbf{F}}$ of $(x_n)_{n \in \mathbb{Z}}$ is defined as

$$[\check{\mathbf{F}}(x_n)](\omega) \triangleq \sum_{n \in \mathbb{Z}} x_n e^{-i\omega n} \qquad \forall (x_n)_{n \in \mathbb{Z}} \in \ell_{\mathbb{R}}^2$$

If we compare the definition of the *Discrete Time Fourier Transform* (Definition 2.38 page 24) to the definition of the Z-transform (Definition 2.35 page 24), we see that the DTFT is just a special case of the more general Z-Transform, with $z = e^{i\omega}$. If we imagine $z \in \mathbb{C}$ as a complex plane, then $e^{i\omega}$ is a unit circle in this plane. The "frequency" $\omega$ in the DTFT is the unit circle in the much larger z-plane as illustrated in Figure 2 (page 25).

---

[70] *Laurent series*: 📖 [1], page 49
[71] 📖 [53], page 344, ⟨Proposition J.1⟩
[72] 📖 [53], pages 344–345, ⟨Theorem J.1⟩

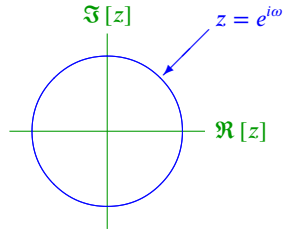

Figure 2:  Unit circle in complex-z plane

**Proposition 2.39** [73] *Let* $\check{x}(\omega) \triangleq \check{F}\big[(x_n)\big](\omega)$ *be the* DISCRETE-TIME FOURIER TRANSFORM *(Definition 2.38 page 24) of a sequence* $(x_n)_{n\in\mathbb{Z}}$ *in* $\ell^2_{\mathbb{R}}$.

$$\underbrace{\check{x}(\omega) = \check{x}(\omega + 2\pi n)}_{\text{PERIODIC with period } 2\pi} \qquad \forall n \in \mathbb{Z}$$

✎ PROOF:

$$\check{x}(\omega + 2\pi n) = \sum_{m\in\mathbb{Z}} x_m e^{-i(\omega+2\pi n)m} \qquad\qquad = \sum_{m\in\mathbb{Z}} x_m e^{-i\omega m} \underbrace{e^{-i2\pi nm}}_{1}$$

$$= \sum_{m\in\mathbb{Z}} x_m e^{-i\omega m} \qquad\qquad\qquad = \check{x}(\omega) \qquad\qquad (2\text{--}1)$$

⇨

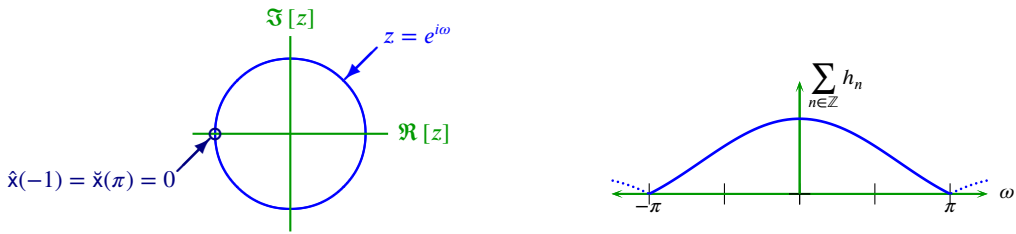

**Proposition 2.40** [74] *Let* $\hat{x}(z)$ *be the* Z-TRANSFORM *(Definition 2.35 page 24) and* $\check{x}(\omega)$ *the* DISCRETE-TIME FOURIER TRANSFORM *(Definition 2.38 page 24) of* $(x_n)$.

$$\underbrace{\left\{ \sum_{n\in\mathbb{Z}} x_n = c \right\}}_{\text{(1) time domain}} \iff \underbrace{\left\{ \hat{x}(z)\Big|_{z=1} = c \right\}}_{\text{(2) z domain}} \iff \underbrace{\left\{ \check{x}(\omega)\Big|_{\omega=0} = c \right\}}_{\text{(3) frequency domain}} \qquad \forall (x_n)_{n\in\mathbb{Z}} \in \ell^2_{\mathbb{R}}, c\in\mathbb{R}$$

---

[73] ✎ [53], pages 348–349, ⟨Proposition J.2⟩
[74] ✎ [53], pages 349–350, ⟨Proposition J.3⟩

**Proposition 2.41** [75]

$$\underbrace{\sum_{n \in \mathbb{Z}} (-1)^n x_n = c}_{\text{(1) in "time"}} \iff \underbrace{\hat{x}(z)\big|_{z=-1} = c}_{\text{(2) in "z domain"}} \iff \underbrace{\breve{x}(\omega)\big|_{\omega=\pi} = c}_{\text{(3) in "frequency"}}$$

$$\iff \underbrace{\left( \sum_{n \in \mathbb{Z}} h_{2n}, \sum_{n \in \mathbb{Z}} h_{2n+1} \right) = \left( \frac{1}{2}\left( \sum_{n \in \mathbb{Z}} h_n + c \right), \frac{1}{2}\left( \sum_{n \in \mathbb{Z}} h_n - c \right) \right)}_{\text{(4) sum of even, sum of odd}}$$

$$\forall c \in \mathbb{R}, \ (x_n)_{n \in \mathbb{Z}}, (y_n)_{n \in \mathbb{Z}} \in \ell_{\mathbb{R}}^2$$

**Lemma 2.42** *Let* $\breve{f}(\omega)$ *be the* DTFT *(Definition 2.38 page 24) of a sequence* $(x_n)_{n \in \mathbb{Z}}$.

$$\underbrace{(x_n \in \mathbb{R})_{n \in \mathbb{Z}}}_{\text{REAL-VALUED sequence}} \implies \underbrace{|\breve{x}(\omega)|^2 = |\breve{x}(-\omega)|^2}_{\text{EVEN}} \qquad \forall (x_n)_{n \in \mathbb{Z}} \in \ell_{\mathbb{R}}^2$$

✎ PROOF:

$$|\breve{x}(\omega)|^2 = |\hat{x}(z)|^2 \big|_{z=e^{i\omega}}$$

$$= \hat{x}(z)\hat{x}^*(z)\big|_{z=e^{i\omega}}$$

$$= \left[ \sum_{n \in \mathbb{Z}} x_n z^{-n} \right] \left[ \sum_{m \in \mathbb{Z}} x_m z^{-m} \right]^* \bigg|_{z=e^{i\omega}}$$

$$= \left[ \sum_{n \in \mathbb{Z}} x_n z^{-n} \right] \left[ \sum_{m \in \mathbb{Z}} x_m^*(z^*)^{-m} \right] \bigg|_{z=e^{i\omega}}$$

$$= \sum_{n \in \mathbb{Z}} \sum_{m \in \mathbb{Z}} x_n x_m^* z^{-n}(z^*)^{-m} \bigg|_{z=e^{i\omega}}$$

$$= \sum_{n \in \mathbb{Z}} \left[ |x_n|^2 + \sum_{m>n} x_n x_m^* z^{-n}(z^*)^{-m} + \sum_{m<n} x_n x_m^* z^{-n}(z^*)^{-m} \right] \bigg|_{z=e^{i\omega}}$$

$$= \sum_{n \in \mathbb{Z}} \left[ |x_n|^2 + \sum_{m>n} x_n x_m e^{i\omega(m-n)} + \sum_{m<n} x_n x_m e^{i\omega(m-n)} \right]$$

$$= \sum_{n \in \mathbb{Z}} \left[ |x_n|^2 + \sum_{m>n} x_n x_m e^{i\omega(m-n)} + \sum_{m>n} x_n x_m e^{-i\omega(m-n)} \right]$$

$$= \sum_{n \in \mathbb{Z}} \left[ |x_n|^2 + \sum_{m>n} x_n x_m \left( e^{i\omega(m-n)} + e^{-i\omega(m-n)} \right) \right]$$

$$= \sum_{n \in \mathbb{Z}} \left[ |x_n|^2 + \sum_{m>n} x_n x_m 2\cos[\omega(m-n)] \right]$$

---

[75] ✎ [27], page 123

$$= \sum_{n\in\mathbb{Z}} |x_n|^2 + 2 \sum_{n\in\mathbb{Z}} \sum_{m>n} x_n x_m \cos[\omega(m-n)]$$

Since $\cos$ is real and even, then $|\check{\mathrm{x}}(\omega)|^2$ must also be real and even.          ⬛

**Theorem 2.43** (inverse DTFT)  [76] *Let* $\check{\mathrm{x}}(\omega)$ *be the* DISCRETE-TIME FOURIER TRANSFORM *(Definition 2.38 page 24) of a sequence* $(x_n)_{n\in\mathbb{Z}} \in \ell_{\mathbb{R}}^2$. *Let* $\check{\mathrm{x}}^{-1}$ *be the inverse of* $\check{\mathrm{x}}$.

$$\underbrace{\left\{ \check{\mathrm{x}}(\omega) \triangleq \sum_{n\in\mathbb{Z}} x_n e^{-i\omega n} \right\}}_{\check{\mathrm{x}}(\omega) \triangleq \check{\mathbf{F}}(x_n)} \implies \underbrace{\left\{ x_n = \frac{1}{2\pi} \int_{\alpha-\pi}^{\alpha+\pi} \check{\mathrm{x}}(\omega) e^{i\omega n}\, d\omega \quad \forall \alpha\in\mathbb{R} \right\}}_{(x_n) = \check{\mathbf{F}}^{-1}\check{\mathbf{F}}(x_n)} \quad \forall (x_n)_{n\in\mathbb{Z}}\in\ell_{\mathbb{R}}^2$$

**Theorem 2.44** (orthonormal quadrature conditions)  [77] *Let* $\check{\mathrm{x}}(\omega)$ *be the* DISCRETE-TIME FOURIER TRANSFORM *(Definition 2.38 page 24) of a sequence* $(x_n)_{n\in\mathbb{Z}} \in \ell_{\mathbb{R}}^2$. *Let* $\bar{\delta}_n$ *be the* KRONECKER DELTA FUNCTION *at* $n$ *(Definition 3.12 page 35).*

$$\sum_{m\in\mathbb{Z}} x_m y_{m-2n}^* = 0 \iff \check{\mathrm{x}}(\omega)\check{\mathrm{y}}^*(\omega) + \check{\mathrm{x}}(\omega+\pi)\check{\mathrm{y}}^*(\omega+\pi) = 0 \quad \forall n\in\mathbb{Z}, \forall (x_n),(y_n)\in\ell_{\mathbb{R}}^2$$

$$\sum_{m\in\mathbb{Z}} x_m x_{m-2n}^* = \bar{\delta}_n \iff |\check{\mathrm{x}}(\omega)|^2 + |\check{\mathrm{x}}(\omega+\pi)|^2 = 2 \quad \forall n\in\mathbb{Z}, \forall (x_n),(y_n)\in\ell_{\mathbb{R}}^2$$

✎PROOF:  Let $z \triangleq e^{i\omega}$.

(1)  Proof that $2\sum_{n\in\mathbb{Z}}\left[\sum_{k\in\mathbb{Z}} x_k y_{k-2n}^*\right] e^{-i2\omega n} = \check{\mathrm{x}}(\omega)\check{\mathrm{y}}^*(\omega) + \check{\mathrm{x}}(\omega+\pi)\check{\mathrm{y}}^*(\omega+\pi)$:

$$2\sum_{n\in\mathbb{Z}}\left[\sum_{k\in\mathbb{Z}} x_k y_{k-2n}^*\right] e^{-i2\omega n}$$

$$= 2\sum_{k\in\mathbb{Z}} x_k \sum_{n\in\mathbb{Z}} y_{k-2n}^* z^{-2n}$$

$$= 2\sum_{k\in\mathbb{Z}} x_k \sum_{n\text{ even}} y_{k-n}^* z^{-n}$$

$$= \sum_{k\in\mathbb{Z}} x_k \sum_{n\in\mathbb{Z}} y_{k-n}^* z^{-n} \left(1 + e^{i\pi n}\right)$$

$$= \sum_{k\in\mathbb{Z}} x_k \sum_{n\in\mathbb{Z}} y_{k-n}^* z^{-n} + \sum_{k\in\mathbb{Z}} x_k \sum_{n\in\mathbb{Z}} y_{k-n}^* z^{-n} e^{i\pi n}$$

$$= \sum_{k\in\mathbb{Z}} x_k \sum_{m\in\mathbb{Z}} y_m^* z^{-(k-m)} + \sum_{k\in\mathbb{Z}} x_k \sum_{m\in\mathbb{Z}} y_m^* e^{-i(\omega+\pi)(k-m)} \quad \text{where } m \triangleq k-n$$

$$= \sum_{k\in\mathbb{Z}} x_k z^{-k} \sum_{m\in\mathbb{Z}} y_m^* z^m + \sum_{k\in\mathbb{Z}} x_k e^{-i(\omega+\pi)k} \sum_{m\in\mathbb{Z}} y_m^* e^{+i(\omega+\pi)m}$$

---

[76] 📖 [76], page 3-95, ⟨(3.6.2)⟩
[77] 📖 [30], pages 132–137, ⟨(5.1.20),(5.1.39)⟩

$$= \sum_{k\in\mathbb{Z}} x_k e^{-i\omega k} \left[ \sum_{m\in\mathbb{Z}} y_m e^{-i\omega m} \right]^* + \sum_{k\in\mathbb{Z}} x_k e^{-i(\omega+\pi)k} \left[ \sum_{m\in\mathbb{Z}} y_m e^{-i(\omega+\pi)m} \right]^*$$

$$\triangleq \breve{x}(\omega)\breve{y}^*(\omega) + \breve{x}(\omega+\pi)\breve{y}^*(\omega+\pi)$$

(2)  Proof that $\sum_{m\in\mathbb{Z}} x_m y^*_{m-2n} = 0 \implies \breve{x}(\omega)\breve{y}^*(\omega) + \breve{x}(\omega+\pi)\breve{y}^*(\omega+\pi) = 0$:

$$0 = 2 \sum_{n\in\mathbb{Z}} \left[ \sum_{k\in\mathbb{Z}} x_k y^*_{k-2n} \right] e^{-i2\omega n} \qquad \text{by left hypothesis}$$

$$= \breve{x}(\omega)\breve{y}^*(\omega) + \breve{x}(\omega+\pi)\breve{y}^*(\omega+\pi) \qquad \text{by item 1}$$

(3)  Proof that $\sum_{m\in\mathbb{Z}} x_m y^*_{m-2n} = 0 \impliedby \breve{x}(\omega)\breve{y}^*(\omega) + \breve{x}(\omega+\pi)\breve{y}^*(\omega+\pi) = 0$:

$$2 \sum_{n\in\mathbb{Z}} \left[ \sum_{k\in\mathbb{Z}} x_k y^*_{k-2n} \right] e^{-i2\omega n} = \breve{x}(\omega)\breve{y}^*(\omega) + \breve{x}(\omega+\pi)\breve{y}^*(\omega+\pi) \qquad \text{by item 1}$$

$$= 0 \qquad \text{by right hypothesis}$$

Thus by the above equation, $\sum_{n\in\mathbb{Z}} \left[ \sum_{k\in\mathbb{Z}} x_k y^*_{k-2n} \right] e^{-i2\omega n} = 0$. The only way for this to be true is if $\sum_{k\in\mathbb{Z}} x_k y^*_{k-2n} = 0$.

(4)  Proof that $\sum_{m\in\mathbb{Z}} x_m x^*_{m-2n} = \bar{\delta}_n \implies |\breve{x}(\omega)|^2 + |\breve{x}(\omega'+\pi)|^2 = 2$:
     Let $g_n \triangleq x_n$.

$$2 = 2 \sum_{n\in\mathbb{Z}} \bar{\delta}_n e^{-i2\omega n}$$

$$= 2 \sum_{n\in\mathbb{Z}} \left[ \sum_{k\in\mathbb{Z}} x_k y^*_{k-2n} \right] e^{-i2\omega n} \qquad \text{by left hypothesis}$$

$$= \breve{x}(\omega)\breve{y}^*(\omega) + \breve{x}(\omega+\pi)\breve{y}^*(\omega+\pi) \qquad \text{by item 1}$$

(5)  Proof that $\sum_{m\in\mathbb{Z}} x_m x^*_{m-2n} = \bar{\delta}_n \impliedby |\breve{x}(\omega)|^2 + |\breve{x}(\omega'+\pi)|^2 = 2$:
     Let $g_n \triangleq x_n$.

$$2 \sum_{n\in\mathbb{Z}} \left[ \sum_{k\in\mathbb{Z}} x_k y^*_{k-2n} \right] e^{-i2\omega n} = \breve{x}(\omega)\breve{y}^*(\omega) + \breve{x}(\omega+\pi)\breve{y}^*(\omega+\pi) \qquad \text{by item 1}$$

$$= 2 \qquad \text{by right hypothesis}$$

Thus by the above equation, $\sum_{n\in\mathbb{Z}} \left[ \sum_{k\in\mathbb{Z}} x_k y^*_{k-2n} \right] e^{-i2\omega n} = 1$. The only way for this to be true is if $\sum_{k\in\mathbb{Z}} x_k y^*_{k-2n} = \bar{\delta}_n$.

### 2.5.3 Frequency Response

The pole zero locations of a digital filter determine the magnitude and phase frequency response of the digital filter.[78] This can be seen by representing the poles and zeros vectors in the complex z-plane. Each of these vectors has a magnitude $M$ and a direction $\theta$. Also, each factor $(z - z_i)$ and $(z - p_i)$ can be represented as vectors as well (the difference of two vectors). Each of these factors can be represented by a magnitude/phase factor $M_i e^{i\theta_i}$. The overall magnitude and phase of $H(z)$ can then be analyzed.

Take the following filter for example:

$$
\begin{aligned}
H(z) &= \frac{b_0 + b_1 z^{-1} + b_2 z^{-2}}{1 + a_1 z^{-1} + a_2 z^{-2}} \\
&= \frac{(z - z_1)(z - z_2)}{(z - p_1)(z - p_2)} \\
&= \frac{M_1 e^{i\theta_1} \, M_2 e^{i\theta_2}}{M_3 e^{i\theta_3} \, M_4 e^{i\theta_4}} \\
&= \left( \frac{M_1 M_2}{M_3 M_4} \right) \left( \frac{e^{i\theta_1} e^{i\theta_2}}{e^{i\theta_3} e^{i\theta_4}} \right)
\end{aligned}
$$

This is illustrated in Figure 3 (page 30). The unit circle represents frequency in the Fourier domain. The frequency response of a filter is just a rotating vector on this circle. The magnitude response of the filter is just then a *vector sum*. For example, the magnitude of any $H(z)$ is as follows:

$$
|H(z)| = \frac{|(z - z_1)| \, |(z - z_2)|}{|(z - p_1)| \, |(z - p_2)|}
$$

### 2.5.4 Filter Banks

*Conjugate quadrature filters* (next definition) are used in *filter banks*. If $\hat{x}(z)$ is a *low-pass filter*, then the conjugate quadrature filter of $\hat{y}(z)$ is a *high-pass filter*.

**Definition 2.45** [79] Let $(x_n)_{n \in \mathbb{Z}}$ and $(y_n)_{n \in \mathbb{Z}}$ be *sequences* (Definition 2.32 page 23) in $\ell^2_{\mathbb{R}}$ (Definition 2.33 page 23). The sequence $(y_n)$ is a **conjugate quadrature filter** with shift $N$ with respect to

---

[78] 📚 [21], pages 90–91
[79] 📚 [127], [page 109](#), 📚 [57], [pages 256–259](#), ⟨section 4.5⟩, 📚 [130], [page 342](#), ⟨(7.2.7), (7.2.8)⟩, 📚 [123], 📚 [122], 📚 [98]

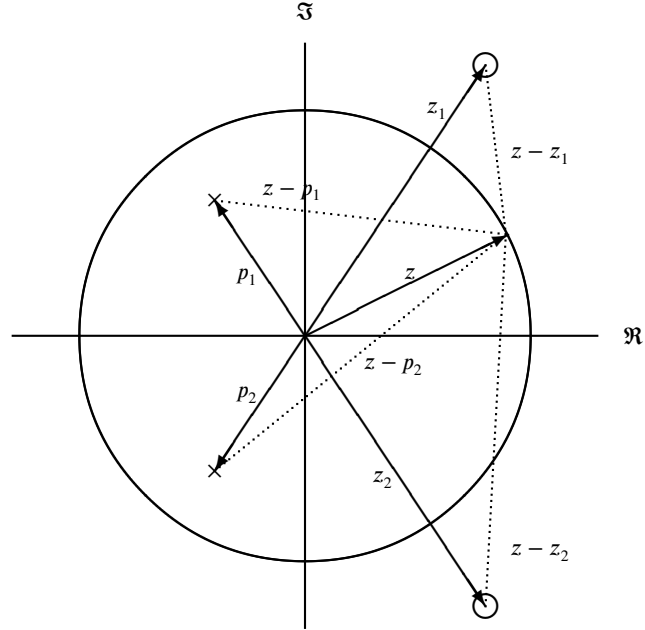

Figure 3: Vector response of digital filter

$(x_n)$ if
$$y_n = \pm(-1)^n x^*_{N-n}$$
A *conjugate quadrature filter* is also called a **CQF** or a **Smith-Barnwell filter**.
Any triple $\big((x_n), (y_n), N\big)$ in this form is said to satisfy the
> **conjugate quadrature filter condition** or the **CQF condition**.

**Theorem 2.46** [80] *Let $\breve{y}(\omega)$ and $\breve{x}(\omega)$ be the* DTFT*s (Definition 2.38 page 24) of the sequences $(y_n)_{n\in\mathbb{Z}}$ and $(x_n)_{n\in\mathbb{Z}}$, respectively, in $\ell^2_{\mathbb{R}}$ (Definition 2.33 page 23).*

$$\underbrace{y_n = \pm(-1)^n x^*_{N-n}}_{(1)\ \text{CQF in “time”}} \iff \hat{y}(z) = \pm(-1)^N z^{-N} \hat{x}^*\left(\frac{-1}{z^*}\right) \quad (2) \quad \text{CQF in “z-domain”}$$

$$\iff \breve{y}(\omega) = \pm(-1)^N e^{-i\omega N} \breve{x}^*(\omega + \pi) \quad (3) \quad \text{CQF in “frequency”}$$

$$\iff x_n = \pm(-1)^N (-1)^n y^*_{N-n} \quad (4) \quad \text{“reversed” CQF in “time”}$$

$$\iff \hat{x}(z) = \pm z^{-N} \hat{y}^*\left(\frac{-1}{z^*}\right) \quad (5) \quad \text{“reversed” CQF in “z-domain”}$$

$$\iff \breve{x}(\omega) = \pm e^{-i\omega N} \breve{y}^*(\omega + \pi) \quad (6) \quad \text{“reversed” CQF in “frequency”}$$

---

[80] ✎ [127], page 109, ✎ [95], pages 236–238, $\langle (7.58),(7.73) \rangle$, ✎ [57], pages 256–259, $\langle$ section 4.5$\rangle$, ✎ [130], page 342, $\langle (7.2.7),\ (7.2.8) \rangle$

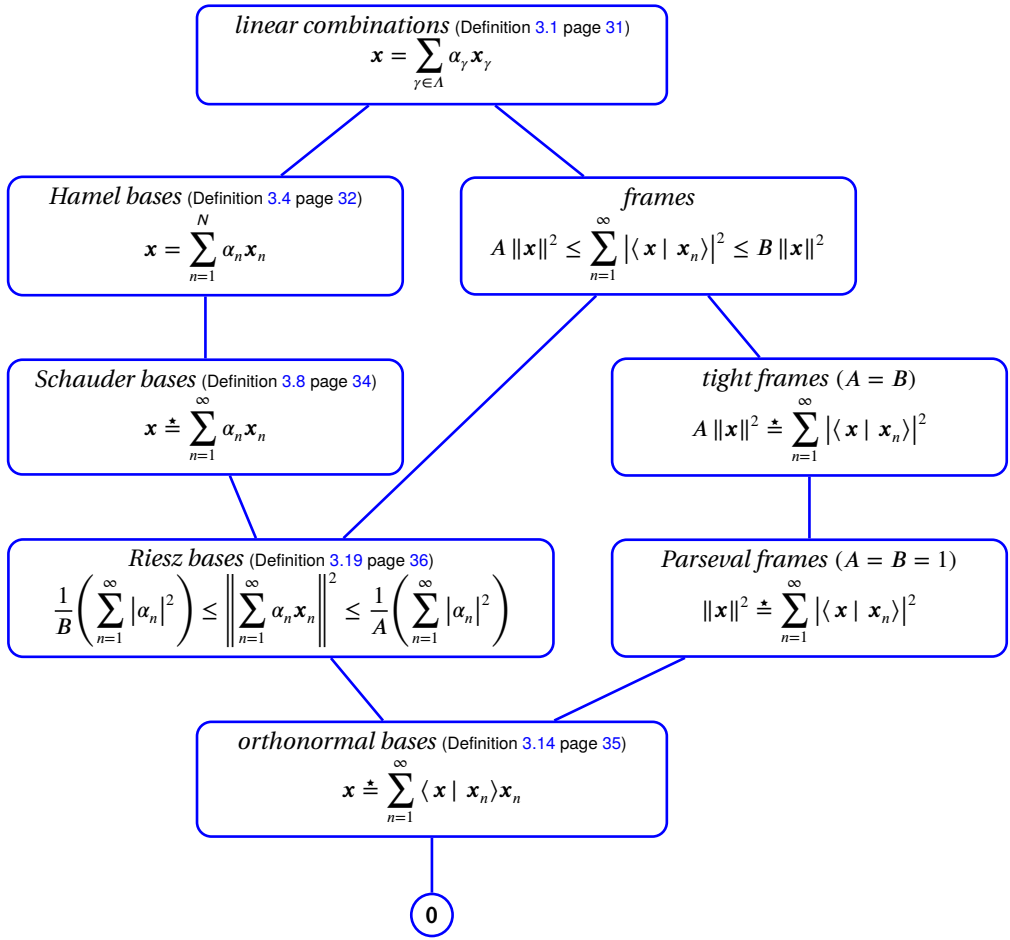

$$\forall N \in \mathbb{Z}$$

# 3 Background: basis theory

## 3.1 Linear Combinations in Linear Spaces

**Definition 3.1** [81] Let $\left\{ \boldsymbol{x}_n \in X \,\middle|\, n=1,2,\dots,N \right\}$ be a set of vectors in a linear space $\left( X, +, \cdot, (\mathbb{F}, +, \dot{\times}) \right)$.

---

[81] ✎ [16], page 11, ⟨Definition I.4.1⟩, ✎ [84], page 46

A vector $x \in X$ is a **linear combination** of the vectors in $\{x_n\}$ if there exists $\{\alpha_n \in \mathbb{F}\,|_{n=1,2,\ldots,N}\}$ such that

$$x = \sum_{n=1}^{N} \alpha_n x_n \,.$$

**Definition 3.2** [82] Let $\{x_n \in X\,|_{n=1,2,\ldots,N}\}$ be a set of vectors in a linear space $(X, +, \cdot, (\mathbb{F}, \dotplus, \dottimes))$. Let $A$ be a subset of $X$.

The **span** of $\{x_n\}$ is defined as $\qquad \mathrm{span}\{x_n\} \triangleq \left\{ \sum_{n=1}^{N} \alpha_n x_n \,\middle|\, \alpha_n \in \mathbb{F} \right\}$ .

The set $\{x_n \in X\}$ **spans** the set $A$ if $\quad A \subseteq \mathrm{span}\{x_n\}$.

**Definition 3.3** [83] Let $\{x_n \in X\,|_{n=1,2,\ldots,N}\}$ be a set of vectors in a linear space $L \triangleq (X, +, \cdot, (\mathbb{F}, \dotplus, \dottimes))$.

The set $\{x_n \in X\,|_{n=1,2,\ldots,N}\}$ is **linearly independent** in $L$ if

$$\sum_{n=1}^{N} \alpha_n x_n = 0 \qquad \Longrightarrow \qquad \alpha_1 = \alpha_2 = \cdots = \alpha_N = 0.$$

If the vectors in $\{x_n\}$ are not linearly independent, then they are **linearly dependent**. An infinite set $\{x_n \in X\,|_{n\in\mathbb{N}}\}$ is *linearly independent* if
every finite subset $\{x_{f(n)}\,|_{n=1,2,\ldots,N}\}$ is *linearly independent*.

**Definition 3.4** [84] Let $\{x_n \in X\,|_{n=1,2,\ldots,N}\}$ be a set of vectors in a linear space $L \triangleq (X, +, \cdot, (\mathbb{F}, \dotplus, \dottimes))$.

The set $\{x_n\}$ is a **Hamel basis** (also called a **linear basis**) for $L$ if
1. $\{x_n\}$ *spans* $L$ $\qquad\qquad\qquad$ and
2. $\{x_n\}$ is *linearly independent*.

If in addition $x = \sum_{n=1}^{N} \alpha_n x_n$, then $\sum_{n=1}^{N} \alpha_n x_n$ is the **expansion** of $x$ on $\{x_n\}$,

and the elements of $(\alpha_n)$ are called the **coordinates** of $x$ with respect to $\{x_n\}$.
If $\alpha_N \neq 0$, then $N$ is the **dimension** $\dim L$ of $L$.

---

## 3.2   Total sets in Topological Linear Spaces

A linear space supports the concept of the *span* of a set of vectors (Definition 3.2 page 32). In a topological linear space $\boldsymbol{\Omega} \triangleq \left( X, +, \cdot, (\mathbb{F}, \dotplus, \dot\times), \boldsymbol{T} \right)$, a set $A$ is said to be *total* in $\boldsymbol{\Omega}$ if the span of $A$ is *dense* in $\boldsymbol{\Omega}$. In this case, $A$ is said to be a *total set* or a *complete set*. However, this use of "complete" in a "*complete set*" is not equivalent to the use of "complete" in a "*complete metric space*".[85] In this text, except for these comments and Definition 3.5, "complete" refers to the metric space definition only.

If a set is both *total* and *linearly independent* (Definition 3.3 page 32) in $\boldsymbol{\Omega}$, then that set is a *Hamel basis* (Definition 3.4 page 32) for $\boldsymbol{\Omega}$.

**Definition 3.5** [86] Let $A^-$ be the *closure* of a $A$ in a *topological linear space* $\boldsymbol{\Omega} \triangleq \left( X, +, \cdot, (\mathbb{F}, \dotplus, \dot\times), \boldsymbol{T} \right)$. Let $\operatorname{span}A$ be the *span* (Definition 3.2 page 32) of a set $A$. A set of vectors $A$ is **total** (or **complete** or **fundamental**) in $\boldsymbol{\Omega}$ if
$$\left( \operatorname{span}A \right)^- = \boldsymbol{\Omega} \qquad \text{(\textit{span} of $A$ is \textit{dense} in $\boldsymbol{\Omega}$).}$$

## 3.3   Total sets in Banach spaces

Often in a linear space we have the option of appending additional structures that offer useful functionality. One of these structures is the *norm* (next definition). The norm of a vector can be described as the "*length*" or the "*magnitude*" of the vector.

**Definition 3.6** [87] Let $\left( X, +, \cdot, (\mathbb{F}, \dotplus, \dot\times) \right)$ be a linear space and $|\cdot| \in \mathbb{R}^{\mathbb{F}}$ the absolute value function.

A functional $\|\cdot\|$ in $\mathbb{R}^X$ is a **norm** if

|  |  |  |  |  |  |
|---|---|---|---|---|---|
| 1. | $\|\boldsymbol{x}\|$ | $\geq$ | $0$ | $\forall \boldsymbol{x} \in X$ | (*strictly positive*) | and |
| 2. | $\|\boldsymbol{x}\|$ | $=$ | $0 \iff \boldsymbol{x} = \mathbb{0}$ | $\forall \boldsymbol{x} \in X$ | (*nondegenerate*) | and |
| 3. | $\|\alpha \boldsymbol{x}\|$ | $=$ | $|\alpha| \, \|\boldsymbol{x}\|$ | $\forall \boldsymbol{x} \in X, \ \alpha \in \mathbb{C}$ | (*homogeneous*) | and |
| 4. | $\|\boldsymbol{x} + \boldsymbol{y}\|$ | $\leq$ | $\|\boldsymbol{x}\| + \|\boldsymbol{y}\|$ | $\forall \boldsymbol{x}, \boldsymbol{y} \in X$ | (*subadditive/ triangle inequality*). |

A **normed linear space** is the tuple $\left( X, +, \cdot, (\mathbb{F}, \dotplus, \dot\times), \|\cdot\| \right)$.

**Definition 3.7** [88] Let $\boldsymbol{B} \triangleq \left( X, +, \cdot, (\mathbb{F}, \dotplus, \dot\times), \|\cdot\| \right)$ be a *Banach space*. $\overset{\scriptscriptstyle\triangle}{=}$ represent *strong convergence* in $\boldsymbol{B}$. That is,

---

[85] ✆ [56], pages 296–297, ⟨6·Orthogonal Bases⟩, ✆ [113], page 78, ⟨Remark 3.50⟩, ✆ [63], page 21, ⟨Remark 1.26⟩

[86] ✆ [139], page 19, ⟨Definition 1.5.1⟩, ✆ [124], page 362, ⟨Definition 9.2.3⟩, ✆ [55], page 134, ⟨Definition 2.4⟩, ✆ [10], pages 149–153, ⟨Definition 9.3, Theorems 9.9 and 9.10⟩

[87] ✆ [3], pages 217–218, ✆ [13], page 53, ✆ [14], page 33, ✆ [12], page 135

[88] ✆ [10], page 138, 247, ⟨Definition 9.1⟩, ✆ [78], page 67, ⟨section 1.1⟩

$$x \stackrel{\scriptscriptstyle\bullet}{=} \sum_{n=1}^{\infty} \alpha_n x_n \quad \stackrel{\text{def}}{\Longleftrightarrow} \quad \lim_{N \to \infty} \left\| x - \sum_{n=1}^{N} \alpha_n x_n \right\| = 0 \quad \stackrel{\text{def}}{\Longleftrightarrow} \quad \left\{ \begin{array}{l} \text{for every } \varepsilon > 0, \\ \text{there exists } M \text{ such that} \\ \text{for all } N > M, \\ \left\| x - \sum_{n=1}^{N} \alpha_n x_n \right\| < \varepsilon \end{array} \right\}$$

**Definition 3.8** [89] Let $B \triangleq \big( X, +, \cdot, (\mathbb{F}, \dotplus, \dottimes), \|\cdot\| \big)$ be a *Banach space*. Let $\stackrel{\scriptscriptstyle\bullet}{=}$ represent *strong convergence* (Definition 3.7 page 33) in $B$.

The countable set $\big\{ x_n \in X \,\big|\, n \in \mathbb{N} \big\}$ is a **Schauder basis** for $B$ if for each $x \in X$

1. $\exists \, (\alpha_n \in \mathbb{F})_{n \in \mathbb{N}}$    such that    $x \stackrel{\scriptscriptstyle\bullet}{=} \sum\limits_{n=1}^{\infty} \alpha_n x_n$   (*strong convergence* in $B$) and

2. $\left\{ \sum\limits_{n=1}^{\infty} \alpha_n x_n \stackrel{\scriptscriptstyle\bullet}{=} \sum\limits_{n=1}^{\infty} \beta_n x_n \right\} \implies \big\{ (\alpha_n) = (\beta_n) \big\}$   (*coefficient functionals* are *unique*)

In this case, $\sum\limits_{n=1}^{\infty} \alpha_n x_n$ is the **expansion** of $x$ on $\big\{ x_n \,\big|\, n \in \mathbb{N} \big\}$ and

the elements of $(\alpha_n)$ are the **coefficient functionals** associated with the basis $\{ x_n \}$.
Coefficient functionals are also called **coordinate functionals**.

**Definition 3.9** [90] Let $\big\{ x_n \,\big|\, n \in \mathbb{N} \big\}$ and $\big\{ y_n \,\big|\, n \in \mathbb{N} \big\}$ be *Schauder bases* of a *Banach space* $\big( X, +, \cdot, (\mathbb{F}, \dotplus, \dottimes), \|\cdot\| \big)$. $\{ x_n \}$ is **equivalent** to $\{ y_n \}$
if there exists a *bounded invertible* operator $\mathbf{R}$ in $X^X$ such that    $\mathbf{R} x_n = y_n$      $\forall n \in \mathbb{Z}$

**Theorem 3.10** [91] *Let* $\big\{ x_n \,\big|\, n \in \mathbb{N} \big\}$ *and* $\big\{ y_n \,\big|\, n \in \mathbb{N} \big\}$ *be* Schauder bases *of a* Banach space $\big( X, +, \cdot, (\mathbb{F}, \dotplus, \dottimes), \|\cdot\| \big)$. $\big\{ \{ x_n \}$ *is* equivalent *to* $\{ y_n \} \big\}$

$$\Longleftrightarrow \qquad \left\{ \sum_{n=1}^{\infty} \alpha_n x_n \text{ is convergent} \iff \sum_{n=1}^{\infty} \alpha_n y_n \text{ is convergent} \right\}$$

## 3.4 Total sets in Hilbert spaces

### 3.4.1 Orthogonal sets in Inner product space

**Definition 3.11** [92] Let $L \triangleq \left( X, +, \cdot, (\mathbb{F}, \dotplus, \dottimes) \right)$ be a linear space.
A function $\langle \triangle \mid \triangledown \rangle \in \mathbb{F}^{X \times X}$ is an **inner product** on $L$ if

| | | | | |
|---|---|---|---|---|
| 1. | $\langle \alpha x \mid y \rangle$ | $= \alpha \langle x \mid y \rangle$ | $\forall x, y \in X, \; \forall \alpha \in \mathbb{C}$ | (*homogeneous*) and |
| 2. | $\langle x + y \mid u \rangle$ | $= \langle x \mid u \rangle + \langle y \mid u \rangle$ | $\forall x, y, u \in X$ | (*additive*) and |
| 3. | $\langle x \mid y \rangle$ | $= \langle y \mid x \rangle^*$ | $\forall x, y \in X$ | (*conjugate symmetric*). and |
| 4. | $\langle x \mid x \rangle$ | $\geq 0$ | $\forall x \in X$ | (*non-negative*) and |
| 5. | $\langle x \mid x \rangle$ | $= 0 \iff x = 0$ | $\forall x \in X$ | (*non-isotropic*) |

An inner product is also called a **scalar product**.
The tuple $\left( X, +, \cdot, (\mathbb{F}, \dotplus, \dottimes), \langle \triangle \mid \triangledown \rangle \right)$ is an **inner product space**.

**Definition 3.12**
The **Kronecker delta function** $\bar{\delta}_n$ is defined as
$$\bar{\delta}_n \triangleq \begin{cases} 1 & \text{for } n = 0 \quad \text{and} \\ 0 & \text{for } n \neq 0. \end{cases} \qquad \forall n \in \mathbb{Z}$$

**Definition 3.13** Let $\left( X, +, \cdot, (\mathbb{F}, \dotplus, \dottimes), \langle \triangle \mid \triangledown \rangle \right)$ be an *inner product space* (Definition 3.11 page 35).
Two vectors $x$ and $y$ in $X$ are **orthogonal** if
$$\langle x \mid y \rangle = \begin{cases} 0 & \text{for } x \neq y \\ c \in \mathbb{F} \setminus 0 & \text{for } x = y \end{cases}$$
The notation $x \perp y$ implies $x$ and $y$ are orthogonal. A set $Y \in 2^X$ is **orthogonal** if $x \perp y \;\; \forall x, y \in Y$. A set $Y$ is **orthonomal** if it is *orthogonal* and $\langle y \mid y \rangle = 1 \;\; \forall y \in Y$. A sequence $(x_n \in X)_{n \in \mathbb{Z}}$ is **orthogonal** if $\langle x_n \mid x_m \rangle = c \bar{\delta}_{nm}$ for some $c \in \mathbb{R} \setminus 0$. A sequence $(x_n \in X)_{n \in \mathbb{Z}}$ is **orthonormal** if $\langle x_n \mid x_m \rangle = \bar{\delta}_{nm}$.

### 3.4.2 Orthonormal bases in Hilbert spaces

**Definition 3.14** Let $\left\{ x_n \in X \mid_{n=1,2,\ldots,N} \right\}$ be a set of vectors in an *inner product space* (Definition 3.11 page 35) $\boldsymbol{\Omega} \triangleq \left( X, +, \cdot, (\mathbb{F}, \dotplus, \dottimes), \langle \triangle \mid \triangledown \rangle \right)$.
The set $\{ x_n \}$ is an **orthogonal basis** for $\boldsymbol{\Omega}$ if $\{ x_n \}$ is *orthogonal* and is
a *Schauder basis* for $\boldsymbol{\Omega}$.
The set $\{ x_n \}$ is an **orthonormal basis** for $\boldsymbol{\Omega}$ if $\{ x_n \}$ is *orthonormal* and is
a *Schauder basis* for $\boldsymbol{\Omega}$.

---

**Definition 3.15** [93] Let $H \triangleq \big(X, +, \cdot, (\mathbb{F}, \dotplus, \dottimes), \langle \triangle \mid \triangledown \rangle\big)$ be a Hilbert space.

Suppose there exists a set $\big\{ x_n \in X \big|_{n \in \mathbb{N}} \big\}$ such that $x \triangleq \sum_{n=1}^{\infty} \langle x \mid x_n \rangle \, x_n$.

Then the quantities $\langle x \mid x_n \rangle$ are called the **Fourier coefficients** of $x$ and the sum $\sum_{n=1}^{\infty} \langle x \mid x_n \rangle \, x_n$ is called the **Fourier expansion** of $x$ or the **Fourier series** for $x$.

**Theorem 3.16** (The Fourier Series Theorem) [94] *Let* $\big\{ x_n \in X \big|_{n \in \mathbb{N}} \big\}$ *be a set of vectors in a* HILBERT SPACE $H \triangleq \big(X, +, \cdot, (\mathbb{F}, \dotplus, \dottimes), \langle \triangle \mid \triangledown \rangle\big)$ *and let* $\|x\| \triangleq \sqrt{\langle x \mid x \rangle}$. (A) $\{x_n\}$ *is* ORTHONORMAL *in* $H$ $\implies$

$$
\left\{
\begin{array}{llll}
& (1). & (\mathrm{span}\{x_n\})^{-} = H & (\langle x_n \rangle \text{ is TOTAL in } H) \\
\iff & (2). & \langle x \mid y \rangle \triangleq \sum_{n=1}^{\infty} \langle x \mid x_n \rangle \langle y \mid x_n \rangle^{*} \quad \forall x, y \in X & (\text{GENERALIZED PARSEVAL'S IDENTITY}) \\
\iff & (3). & \|x\|^2 \triangleq \sum_{n=1}^{\infty} \big| \langle x \mid x_n \rangle \big|^2 \quad \forall x \in X & (\text{PARSEVAL'S IDENTITY}) \\
\iff & (4). & x \triangleq \sum_{n=1}^{\infty} \langle x \mid x_n \rangle \, x_n \quad \forall x \in X & (\text{FOURIER SERIES EXPANSION})
\end{array}
\right\}
$$

**Theorem 3.17** [95] *Let* $H$ *be a* HILBERT SPACE.

     $H$ *has a* SCHAUDER BASIS    $\iff$    $H$ *is* SEPARABLE

**Theorem 3.18** [96] *Let* $H$ *be a* HILBERT SPACE.

     $H$ *has an* ORTHONORMAL BASIS    $\iff$    $H$ *is* SEPARABLE

## 3.5 Riesz bases in Hilbert space

**Definition 3.19** [97] Let $\big\{ x_n \in X \big|_{n \in \mathbb{N}} \big\}$ be a set of vectors in a *separable Hilbert space* $H \triangleq \big(X, +, \cdot, (\mathbb{F}, \dotplus, \dottimes), \langle \triangle \mid \triangledown \rangle\big)$. $\{x_n\}$ is a **Riesz basis** for $H$ if $\{x_n\}$ is *equivalent* (Definition 3.9 page 34) to some *orthonormal basis* (Definition 3.14 page 35) in $H$.

**Definition 3.20** [98] Let $(x_n \in X)_{n \in \mathbb{N}}$ be a sequence of vectors in a *separable Hilbert space* $H \triangleq (X, +, \cdot, (\mathbb{F}, \dotplus, \dot\times), \langle \triangle \mid \triangledown \rangle)$. The sequence $(x_n)$ is a **Riesz sequence** for $H$ if

$$\exists A, B \in \mathbb{R}^+ \quad \text{such that} \quad A\left(\sum_{n=1}^{\infty} |\alpha_n|^2\right) \leq \left\|\sum_{n=1}^{\infty} \alpha_n x_n\right\|^2 \leq B\left(\sum_{n=1}^{\infty} |\alpha_n|^2\right) \qquad \forall (\alpha_n) \in \ell_{\mathbb{F}}^2.$$

**Lemma 3.21** [99] *Let $\{x_n|_{n \in \mathbb{N}}\}$ be a sequence in a* Hilbert space *$X \triangleq (X, +, \cdot, (\mathbb{F}, \dotplus, \dot\times), \langle \triangle \mid \triangledown \rangle)$.*
*Let $\{y_n|_{n \in \mathbb{N}}\}$ be a sequence in a* Hilbert space *$Y \triangleq (Y, +, \cdot, (\mathbb{F}, \dotplus, \dot\times), \langle \triangle \mid \triangledown \rangle)$. Let*

$$\left\{ \begin{array}{ll} (i). & \{x_n\} \text{ is } \text{total } in \; X \hfill and \\[2mm] (ii). & \text{There exists } A > 0 \text{ such that } A \sum_{n \in C} |a_n|^2 \leq \left\|\sum_{n \in C} a_n x_n\right\|^2 \quad \text{for finite } C \quad and \\[2mm] (iii). & \text{There exists } B > 0 \text{ such that } \left\|\sum_{n=1}^{\infty} b_n y_n\right\|^2 \leq B \sum_{n=1}^{\infty} |b_n|^2 \quad \forall (b_n)_{n \in \mathbb{N}} \in \ell_{\mathbb{F}}^2 \end{array} \right\} \implies$$

$$\left\{ \begin{array}{ll} (1). & \mathbf{R}^\circ \text{ is a linear bounded operator that maps from } \operatorname{span}\{x_n\} \text{ to } \operatorname{span}\{y_n\} \\ & \text{where } \mathbf{R}^\circ \sum_{n \in C} c_n x_n \triangleq \sum_{n \in C} c_n y_n, \text{ for some sequence } (c_n) \text{ and finite set } C \quad and \\[2mm] (2). & \mathbf{R} \text{ has a unique extension to a bounded operator } \mathbf{R} \text{ that maps from } X \text{ to } Y \quad and \\[2mm] (3). & \|\|\mathbf{R}^\circ\|\| \leq \frac{B}{A} \hfill and \\[2mm] (4). & \|\|\mathbf{R}\|\| \leq \frac{B}{A} \end{array} \right\}$$

**Theorem 3.22** [100] *Let $\{x_n \in X|_{n \in \mathbb{N}}\}$ be a set of vectors in a* separable Hilbert space *$H \triangleq (X, +, \cdot, (\mathbb{F}, \dotplus, \dot\times), \langle \triangle \mid \triangledown \rangle)$.*

$$\left\{ \begin{array}{l} \{x_n\} \text{ is a } \text{Riesz basis} \\ \text{for } H \end{array} \right\} \iff \left\{ \begin{array}{ll} (1). & \{x_n\} \text{ is } \text{total } in \; H \hfill and \\ (2). & \exists_{A,B \in \mathbb{R}^+} \text{ such that } \quad \forall (\alpha_n) \in \ell_{\mathbb{F}}^2, \\ & A \sum_{n=1}^{\infty} |\alpha_n|^2 \leq \left\|\sum_{n=1}^{\infty} \alpha_n x_n\right\|^2 \leq B \sum_{n=1}^{\infty} |\alpha_n|^2 \end{array} \right\}$$

**Theorem 3.23** [101] *Let $H \triangleq (X, +, \cdot, (\mathbb{F}, \dotplus, \dot\times), \langle \triangle \mid \triangledown \rangle)$ be a* separable Hilbert space.

---

$$\left\{ \begin{array}{l} (\!(x_n \in H)\!)_{n \in \mathbb{Z}} \text{ is a} \\ \text{Riesz basis } \textit{for } H \end{array} \right\} \implies \left\{ \begin{array}{l} \textit{There exists } (\!(y_n \in H)\!)_{n \in \mathbb{Z}} \textit{ such that} \\ \quad (1). \quad (\!(x_n)\!) \textit{ and } (\!(y_n)\!) \textit{ are } \text{biorthogonal} \quad \textit{and} \\ \quad (2). \quad (\!(y_n)\!) \textit{ is also a } \text{Riesz basis } \textit{for } H \quad \textit{and} \\ \quad (3). \quad \exists B \geq A > 0 \quad \textit{such that} \\ \qquad A \sum_{n=1}^{\infty} |a_n|^2 \leq \left\| \sum_{n=1}^{\infty} a_n x_n \right\|^2 = \left\| \sum_{n=1}^{\infty} a_n y_n \right\|^2 \leq B \sum_{n=1}^{\infty} |a_n|^2 \\ \qquad \scriptstyle \forall (a_n)_{n \in \mathbb{N}} \in \ell_\mathbb{F}^2 \end{array} \right\}$$

**Proposition 3.24** [102] *Let* $\left\{ x_n \big|_{n \in \mathbb{N}} \right\}$ *be a set of vectors in a* Hilbert space
$H \triangleq \left( X, +, \cdot, (\mathbb{F}, \dotplus, \dot{\times}), \langle \triangle \,|\, \nabla \rangle \right)$.

$$\left\{ \begin{array}{l} \{x_n\} \textit{ is a } \text{Riesz basis } \textit{for } H \textit{ with} \\ A \sum_{n=1}^{\infty} |a_n|^2 \leq \left\| \sum_{n=1}^{\infty} a_n x_n \right\|^2 \leq B \sum_{n=1}^{\infty} |a_n|^2 \\ \scriptstyle \forall \{a_n\} \in \ell_\mathbb{F}^2 \end{array} \right\} \implies \left\{ \begin{array}{l} \{x_n\} \textit{ is a } \text{frame } \textit{for } H \textit{ with} \\ \underbrace{\frac{1}{B} \|x\|^2 \leq \sum_{n=1}^{\infty} |\langle x \,|\, x_n \rangle|^2 \leq \frac{1}{A} \|x\|^2}_{\text{stability condition}} \\ \scriptstyle \forall x \in H \end{array} \right\}$$

# 4 Main results: transversal operators

## 4.1 Definitions

The **transversal operators T** and **D** (next definition) are the main focus of this paper.

**Definition 4.1** [103]

1. **T** is the **translation operator** on $\mathbb{C}^\mathbb{C}$ defined as

   $$\mathbf{T}_\tau f(x) \triangleq f(x - \tau) \quad \text{and} \quad \mathbf{T} \triangleq \mathbf{T}_1 \qquad \forall f \in \mathbb{C}^\mathbb{C}$$

2. **D** is the **dilation operator** on $\mathbb{C}^\mathbb{C}$ defined as

   $$\mathbf{D}_\alpha f(x) \triangleq f(\alpha x) \quad \text{and} \quad \mathbf{D} \triangleq \sqrt{2} \mathbf{D}_2 \qquad \forall f \in \mathbb{C}^\mathbb{C}$$

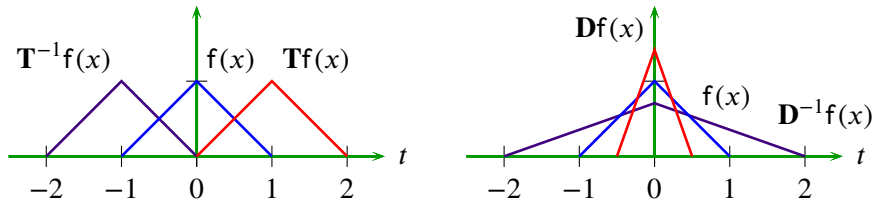

---

## 4.2 Properties

### 4.2.1 Algebraic properties

**Proposition 4.2** [104] *Let* $\mathbf{T}$ *be the* TRANSLATION OPERATOR *(Definition 4.1 page 38).*

$$\sum_{n\in\mathbb{Z}} \mathbf{T}^n f(x) = \sum_{n\in\mathbb{Z}} \mathbf{T}^n f(x+1) \qquad \forall f \in \mathbb{R}^{\mathbb{R}} \qquad \left( \sum_{n\in\mathbb{Z}} \mathbf{T}^n f(x) \text{ is PERIODIC with period } 1 \right)$$

✎PROOF:

$$\sum_{n\in\mathbb{Z}} \mathbf{T}^n f(x+1) = \sum_{n\in\mathbb{Z}} f(x-n+1) \qquad \text{by definition of } \mathbf{T} \text{ Definition 4.1 page 38}$$

$$= \sum_{m\in\mathbb{Z}} f(x-m) \qquad \text{where } m \triangleq n-1$$

$$= \sum_{m\in\mathbb{Z}} \mathbf{T}^m f(x) \qquad \text{by definition of } \mathbf{T} \text{ Definition 4.1 page 38}$$

          ✎

In a linear space, every operator has an *inverse*. Although the inverse always exists as a relation, it may not exist as a function or as an operator. But in some cases the inverse of an operator is itself an operator. The inverses of the operators $\mathbf{T}$ and $\mathbf{D}$ both exists as operators, as demonstrated by Proposition 4.3 (next).

**Proposition 4.3** [105] *Let* $\mathbf{T}$ *and* $\mathbf{D}$ *be as defined in Definition 4.1 page 38.*
    $\mathbf{T}$ *has an inverse* $\mathbf{T}^{-1}$ *in* $\mathbb{C}^{\mathbb{C}}$ *expressed by the relation*

    $\mathbf{T}^{-1} f(x) = f(x+1) \qquad \forall f \in \mathbb{C}^{\mathbb{C}} \qquad \textbf{(translation operator inverse).}$

    $\mathbf{D}$ *has an inverse* $\mathbf{D}^{-1}$ *in* $\mathbb{C}^{\mathbb{C}}$ *expressed by the relation*

    $\mathbf{D}^{-1} f(x) = \frac{\sqrt{2}}{2} f\left(\frac{1}{2}x\right) \qquad \forall f \in \mathbb{C}^{\mathbb{C}} \qquad \textbf{(dilation operator inverse).}$

✎PROOF:

1. Proof that $\mathbf{T}^{-1}$ is the inverse of $\mathbf{T}$:

$$\mathbf{T}^{-1}\mathbf{T}f(x) = \mathbf{T}^{-1}f(x-1) \qquad \text{by Proposition 4.3 page 39}$$

$$= f([x+1]-1)$$

$$= f(x)$$

$$= f([x-1]+1)$$

$$= \mathbf{T}f(x+1) \qquad \text{by defintion of } \mathbf{T} \text{ (Definition 4.1 page 38)}$$

$$= \mathbf{T}\mathbf{T}^{-1}f(x)$$

$$\implies \mathbf{T}^{-1}\mathbf{T} = \mathbf{I} = \mathbf{T}\mathbf{T}^{-1}$$

---

[104] ✎ [53], page 3
[105] ✎ [53], page 3

2. Proof that $\mathbf{D}^{-1}$ is the inverse of $\mathbf{D}$:

$$
\begin{aligned}
\mathbf{D}^{-1}\mathbf{D}f(x) &= \mathbf{D}^{-1}\sqrt{2}f(2x) && \text{by defintion of } \mathbf{D} \text{ (Definition 4.1 page 38)}\\
&= \left(\tfrac{\sqrt{2}}{2}\right)\sqrt{2}f\left(2\left[\tfrac{1}{2}x\right]\right)\\
&= f(x)\\
&= \sqrt{2}\left[\tfrac{\sqrt{2}}{2}f\left(\tfrac{1}{2}[2x]\right)\right]\\
&= \mathbf{D}\left[\tfrac{\sqrt{2}}{2}f\left(\tfrac{1}{2}x\right)\right] && \text{by defintion of } \mathbf{D} \text{ (Definition 4.1 page 38)}\\
&= \mathbf{D}\mathbf{D}^{-1}f(x)\\
\implies \mathbf{D}^{-1}\mathbf{D} &= \mathbf{I} = \mathbf{D}\mathbf{D}^{-1}
\end{aligned}
$$

**Proposition 4.4** [106] *Let* $\mathbf{T}$ *and* $\mathbf{D}$ *be as defined in Definition 4.1 page 38. Let* $\mathbf{D}^0 = \mathbf{T}^0 \triangleq \mathbf{I}$ *be the* IDENTITY OPERATOR.

$$\mathbf{D}^j\mathbf{T}^nf(x) = 2^{j/2}f\left(2^jx - n\right) \qquad \forall j,n\in\mathbb{Z},\, f\in\mathbb{C}^{\mathbb{C}}$$

### 4.2.2 Linear space properties

**Definition 4.5** [107] Let $+$ be an addition operator on a tuple $\langle\!\langle x_n\rangle\!\rangle_m^N$.

The **summation** of $\langle\!\langle x_n\rangle\!\rangle$ from index $m$ to index $N$ with respect to $+$ is

$$
\sum_{n=m}^N x_n \triangleq
\begin{cases}
0 & \text{for } N < m\\
\left(\displaystyle\sum_{n=m}^{N-1} x_n\right) + x_N & \text{for } N \geq m
\end{cases}
$$

An infinite summation $\sum_{n=1}^\infty \phi_n$ is meaningless outside some topological space (e.g. metric space, normed space, etc.). The sum $\sum_{n=1}^\infty \phi_n$ is an abbreviation for $\lim_{N\to\infty}\sum_{n=1}^N \phi_n$ (the limit of partial sums). And the concept of limit is also itself meaningless outside of a topological space.

**Definition 4.6** [108] Let $(X,\mathcal{T})$ be a topological space and $\lim$ be the limit induced by the topology $\mathcal{T}$.

$$\sum_{n=1}^{\infty} \boldsymbol{x}_n \triangleq \sum_{n \in \mathbb{N}} \boldsymbol{x}_n \triangleq \lim_{N \to \infty} \sum_{n=1}^{N} \boldsymbol{x}_n$$

$$\sum_{n=-\infty}^{\infty} \boldsymbol{x}_n \triangleq \sum_{n \in \mathbb{Z}} \boldsymbol{x}_n \triangleq \lim_{N \to \infty} \left( \sum_{n=0}^{N} \boldsymbol{x}_n \right) + \left( \lim_{N \to -\infty} \sum_{n=-1}^{N} \boldsymbol{x}_n \right)$$

In general the operators $\mathbf{T}$ and $\mathbf{D}$ are *noncommutative* ($\mathbf{TD} \neq \mathbf{DT}$), as demonstrated by Proposition 4.8 and by the following illustration.

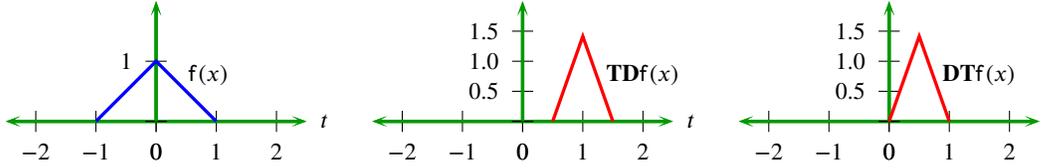

**Proposition 4.7**  *Let $\mathbf{T}$ and $\mathbf{D}$ be as in Definition 4.1 page 38.*

$$\mathbf{D}^j \mathbf{T}^n \big[ \mathbf{fg} \big] = 2^{-j/2} \big[ \mathbf{D}^j \mathbf{T}^n \mathbf{f} \big] \big[ \mathbf{D}^j \mathbf{T}^n \mathbf{g} \big] \qquad \forall j, n \in \mathbb{Z}, \, \mathbf{f} \in \mathbb{C}^{\mathbb{C}}$$

✎ PROOF:

$$\begin{aligned}
\mathbf{D}^j \mathbf{T}^n \big[ \mathbf{f}(x)\mathbf{g}(x) \big] &= 2^{j/2} \mathbf{f}(2^j x - n) \mathbf{g}(2^j x - n) && \text{by Proposition 4.4 page 40} \\
&= 2^{-j/2} \big[ 2^{j/2} \mathbf{f}(2^j x - n) \big] \big[ 2^{j/2} \mathbf{g}(2^j x - n) \big] \\
&= 2^{-j/2} \big[ \mathbf{D}^j \mathbf{T}^n \mathbf{f}(x) \big] \big[ \mathbf{D}^j \mathbf{T}^n \mathbf{g}(x) \big] && \text{by Proposition 4.4 page 40}
\end{aligned}$$

✎

**Proposition 4.8**  (commutator relation)  [109] *Let $\mathbf{T}$ and $\mathbf{D}$ be as in Definition 4.1 page 38.*

$$\begin{aligned}
\mathbf{D}^j \mathbf{T}^n &= \mathbf{T}^{2^{-j/2}n} \mathbf{D}^j && \forall j, n \in \mathbb{Z} \\
\mathbf{T}^n \mathbf{D}^j &= \mathbf{D}^j \mathbf{T}^{2^j n} && \forall n, j \in \mathbb{Z}
\end{aligned}$$

✎ PROOF:

$$\begin{aligned}
\mathbf{D}^j \mathbf{T}^{2^j n} \mathbf{f}(x) &= 2^{j/2} \mathbf{f}(2^j x - 2^j n) && \text{by Proposition 4.7 page 41} \\
&= 2^{j/2} \mathbf{f}\big( 2^j [x - n] \big) && \text{by } \textit{distributivity} \text{ property of the field } (\mathbb{R}, +, \times) \\
&= \mathbf{T}^n 2^{j/2} \mathbf{f}(2^j x) && \text{by definition of } \mathbf{T} \text{ (Definition 4.1 page 38)} \\
&= \mathbf{T}^n \mathbf{D}^j \mathbf{f}(x) && \text{by definition of } \mathbf{D} \text{ (Definition 4.1 page 38)}
\end{aligned}$$

$$\begin{aligned}
\mathbf{D}^j \mathbf{T}^n \mathbf{f}(x) &= 2^{j/2} \mathbf{f}(2^j x - n) && \text{by Proposition 4.7 page 41} \\
&= 2^{j/2} \mathbf{f}\big( 2^j [x - 2^{-j/2}n] \big) && \text{by } \textit{distributivity} \text{ property of the field } (\mathbb{R}, +, \times) \\
&= \mathbf{T}^{2^{-j/2}n} 2^{j/2} \mathbf{f}(2^j x) && \text{by definition of } \mathbf{T} \text{ (Definition 4.1 page 38)} \\
&= \mathbf{T}^{2^{-j/2}n} \mathbf{D}^j \mathbf{f}(x) && \text{by definition of } \mathbf{D} \text{ (Definition 4.1 page 38)}
\end{aligned}$$

---

[109] ✎ [25], page 42, ⟨equation (2.9)⟩, ✎ [28], page 21, ✎ [50], page 641, ✎ [51], page 110

### 4.2.3  Inner-product space properties

In an inner product space, every operator has an *adjoint* (Proposition 1.33 page 10) and this adjoint is always itself an operator. In the case where the adjoint and inverse of an operator **U** coincide, then **U** is said to be *unitary* (Definition 1.48 page 13). And in this case, **U** has several nice properties (see Proposition 4.14 and Theorem 4.15 page 45). Proposition 4.9 (next) gives the adjoints of **D** and **T**, and Proposition 4.10 (page 43) demonstrates that both **D** and **T** are unitary. Other examples of unitary operators include the *Fourier Transform operator* $\tilde{\mathbf{F}}$ and the *rotation matrix operator*.

**Proposition 4.9**  *Let **T** be the translation operator (Definition 4.1 page 38) with adjoint $\mathbf{T}^*$ and **D** the dilation operator with adjoint $\mathbf{D}^*$.*

$$
\begin{aligned}
\mathbf{T}^*\mathsf{f}(x) &= \mathsf{f}(x+1) &&\forall \mathsf{f} \in \mathcal{L}_{\mathbb{R}}^2 &&\textit{(translation operator adjoint)}\\
\mathbf{D}^*\mathsf{f}(x) &= \tfrac{\sqrt{2}}{2}\mathsf{f}\!\left(\tfrac{1}{2}x\right) &&\forall \mathsf{f} \in \mathcal{L}_{\mathbb{R}}^2 &&\textit{(dilation operator adjoint)}
\end{aligned}
$$

✎ PROOF:

1. Proof that $\mathbf{T}^*\mathsf{f}(x) = \mathsf{f}(x+1)$:

$$
\begin{aligned}
\langle \mathsf{g}(x) \mid \mathbf{T}^*\mathsf{f}(x)\rangle &= \langle \mathsf{g}(u) \mid \mathbf{T}^*\mathsf{f}(u)\rangle &&\text{by change of dummy variable } x \to u\\
&= \langle \mathbf{T}\mathsf{g}(u) \mid \mathsf{f}(u)\rangle &&\text{by definition of adjoint } \mathbf{T}^*\\
&= \langle \mathsf{g}(u-1) \mid \mathsf{f}(u)\rangle &&\text{by definition of } \mathbf{T}\text{ (Definition 4.1 page 38)}\\
&= \langle \mathsf{g}(x) \mid \mathsf{f}(x+1)\rangle &&\text{where } x \triangleq u-1 \implies u = x+1\\
\implies \mathbf{T}^*\mathsf{f}(x) &= \mathsf{f}(x+1)
\end{aligned}
$$

2. Proof that $\mathbf{D}^*\mathsf{f}(x) = \tfrac{\sqrt{2}}{2}\mathsf{f}\!\left(\tfrac{1}{2}x\right)$:

$$
\begin{aligned}
\langle \mathsf{g}(x) \mid \mathbf{D}^*\mathsf{f}(x)\rangle &= \langle \mathsf{g}(u) \mid \mathbf{D}^*\mathsf{f}(u)\rangle &&\text{by change of dummy variable } t \to u\\
&= \langle \mathbf{D}\mathsf{g}(u) \mid \mathsf{f}(u)\rangle &&\text{by definition of } \mathbf{D}^*\\
&= \left\langle \sqrt{2}\mathsf{g}(2u) \,\middle|\, \mathsf{f}(u)\right\rangle &&\text{by definition of } \mathbf{D}\text{ (Definition 4.1 page 38)}\\
&= \int_{u \in \mathbb{R}} \sqrt{2}\mathsf{g}(2u)\mathsf{f}^*(u)\,\mathrm{d}u &&\text{by definition of } \langle \triangle \mid \triangledown\rangle\\
&= \int_{x \in \mathbb{R}} \mathsf{g}(x)\left[\sqrt{2}\mathsf{f}\!\left(\tfrac{x}{2}\right)\tfrac{1}{2}\right]^*\mathrm{d}x &&\text{where } x = 2u\\
&= \left\langle \mathsf{g}(x) \mid \tfrac{\sqrt{2}}{2}\mathsf{f}\!\left(\tfrac{x}{2}\right)\right\rangle &&\text{by definition of } \langle \triangle \mid \triangledown\rangle\\
\implies \mathbf{D}^*\mathsf{f}(x) &= \tfrac{\sqrt{2}}{2}\mathsf{f}\!\left(\tfrac{x}{2}\right)
\end{aligned}
$$

**Proposition 4.10** [110] *Let* $\mathbf{T}$ *and* $\mathbf{D}$ *be as in Definition 4.1 page 38. Let* $\mathbf{T}^{-1}$ *and* $\mathbf{D}^{-1}$ *be as in Proposition 4.3 page 39.*

$\quad\quad \mathbf{T}$ *is* UNITARY *in* $L_{\mathbb{R}}^2$  $(\mathbf{T}^{-1} = \mathbf{T}^*$ *in* $L_{\mathbb{R}}^2).$
$\quad\quad \mathbf{D}$ *is* UNITARY *in* $L_{\mathbb{R}}^2$  $(\mathbf{D}^{-1} = \mathbf{D}^*$ *in* $L_{\mathbb{R}}^2).$

✎PROOF:

$\quad \mathbf{T}^{-1} = \mathbf{T}^*$            by Proposition 4.3 page 39 and Proposition 4.9 page 42
$\quad\quad \implies \quad \mathbf{T}$ *is unitary*      by the definition of *unitary* operators (Definition 1.48 page 13)

$\quad \mathbf{D}^{-1} = \mathbf{D}^*$            by Proposition 4.3 page 39 and Proposition 4.9 page 42
$\quad\quad \implies \quad \mathbf{D}$ *is unitary*      by the definition of *unitary* operators (Definition 1.48 page 13)

$\quad\quad\quad\quad\quad\quad\quad\quad\quad\quad\quad\quad\quad\quad\quad\quad\quad\quad\quad\quad\quad\quad\quad\quad\quad\quad$ ✍

### 4.2.4  Normed linear space properties

**Proposition 4.11**  *Let* $\mathbf{D}$ *be the* DILATION OPERATOR *(Definition 4.1 page 38).*

$$\left\{ \begin{array}{ll} (1). & \mathbf{D}f(x) = \sqrt{2}f(x) \quad \text{and} \\ (2). & f(x) \text{ is CONTINUOUS} \end{array} \right\} \quad \implies \quad \{f(x) \text{ is a CONSTANT}\} \quad \forall f \in L_{\mathbb{R}}^2$$

✎PROOF:

(1)  Proof that (1) $\impliedby$ *constant* property:

$\quad\quad \mathbf{D}f(x) \triangleq \sqrt{2}f(2x)$          by definition of $\mathbf{D}$ (Definition 4.1 page 38)

$\quad\quad\quad\quad = \sqrt{2}f(x)$            by *constant* hypothesis

(2)  Proof that (2) $\impliedby$ *constant* property:

$\quad\quad \|f(x) - f(x+h)\| = \|f(x) - f(x)\|$        by *constant* hypothesis

$\quad\quad\quad\quad = \|0\|$

$\quad\quad\quad\quad = 0$      by *nondegenerate* property of $\|\cdot\|$ (Definition 3.6 page 33)

$\quad\quad\quad\quad \le \varepsilon$

$\quad\quad\quad\quad \implies \forall h > 0, \ \exists \varepsilon \quad \text{such that} \quad \|f(x) - f(x+h)\| < \varepsilon$

$\quad\quad\quad\quad \overset{\text{def}}{\iff} f(x) \text{ is } continuous$

(3)  Proof that (1,2) $\implies$ *constant* property:

---

[110] ✍ [25], page 41, ⟨Lemma 2.5.1⟩, ✍ [137], page 18, ⟨Lemma 2.5⟩

(a)   Suppose there exists $x, y \in \mathbb{R}$ such that $\mathsf{f}(x) \neq \mathsf{f}(y)$.

(b)   Let $(x_n)_{n\in\mathbb{N}}$ be a sequence with limit $x$ and $(y_n)_{n\in\mathbb{N}}$ a sequence with limit $y$

(c)   Then

$$
\begin{aligned}
0 &< \|\mathsf{f}(x) - \mathsf{f}(y)\| && \text{by assumption in item 3a page 44} \\
&= \lim_{n\to\infty} \|\mathsf{f}(x_n) - \mathsf{f}(y_n)\| && \text{by (2) and definition of } (x_n) \text{ and } (y_n) \text{ (item 3b page 44)} \\
&= \lim_{n\to\infty} \|\mathsf{f}(2^m x_n) - \mathsf{f}(2^\ell y_n)\| && \forall m, \ell \in \mathbb{Z}; \text{ by (1)} \\
&= 0
\end{aligned}
$$

(d)   But this is a *contradiction*, so $\mathsf{f}(x) = \mathsf{f}(y)$ for all $x, y \in \mathbb{R}$, and $\mathsf{f}(x)$ is *constant*.

Note that in Proposition 4.11, it is not possible to remove the *continuous* constraint outright (next two counterexamples).

**Counterexample 4.12**   Let $\mathsf{f}(x)$ be a function in $\mathbb{R}^{\mathbb{R}}$.

Let $\mathsf{f}(x) \triangleq \begin{cases} 0 & \text{for } x = 0 \\ 1 & \text{otherwise.} \end{cases}$

Then $\mathbf{D}\mathsf{f}(x) \triangleq \sqrt{2}\mathsf{f}(2x) = \sqrt{2}\mathsf{f}(x)$, but $\mathsf{f}(x)$ is *not constant*.

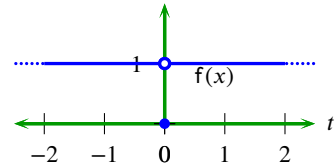

**Counterexample 4.13**   Let $\mathsf{f}(x)$ be a function in $\mathbb{R}^{\mathbb{R}}$. Let $\mathbb{Q}$ be the set of *rational numbers* and $\mathbb{R}\setminus\mathbb{Q}$ the set of *irrational numbers*.

Let $\mathsf{f}(x) \triangleq \begin{cases} 1 & \text{for } x \in \mathbb{Q} \\ -1 & \text{for } x \in \mathbb{R}\setminus\mathbb{Q}. \end{cases}$

Then $\mathbf{D}\mathsf{f}(x) \triangleq \sqrt{2}\mathsf{f}(2x) = \sqrt{2}\mathsf{f}(x)$, but $\mathsf{f}(x)$ is *not constant*.

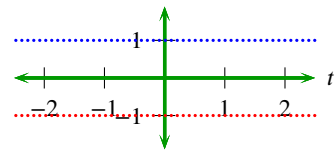

**Proposition 4.14**   (Operator norm)   *Let* $\mathbf{T}$ *and* $\mathbf{D}$ *be as in Definition 4.1 page 38. Let* $\mathbf{T}^{-1}$ *and* $\mathbf{D}^{-1}$ *be as in Proposition 4.3 page 39. Let* $\mathbf{T}^*$ *and* $\mathbf{D}^*$ *be as in Proposition 4.9 page 42. Let* $\|\cdot\|$ *and* $\langle \triangle \mid \triangledown \rangle$ *be as in Definition 2.2 page 15. Let* $\|\|\cdot\|\|$ *be the operator norm (Definition 1.22 page 8) induced by* $\|\cdot\|$.

$$\|\|\mathbf{T}\|\| = \|\|\mathbf{D}\|\| = \|\|\mathbf{T}^*\|\| = \|\|\mathbf{D}^*\|\| = \|\|\mathbf{T}^{-1}\|\| = \|\|\mathbf{D}^{-1}\|\| = 1$$

✎ PROOF:   These results follow directly from the fact that $\mathbf{T}$ and $\mathbf{D}$ are *unitary* (Proposition 4.10 page 43) and from Theorem 1.50 page 14 and Theorem 1.51 page 14.

**Theorem 4.15**  *Let* $\mathbf{T}$ *and* $\mathbf{D}$ *be as in Definition 4.1 page 38. Let* $\mathbf{T}^{-1}$ *and* $\mathbf{D}^{-1}$ *be as in Proposition 4.3 page 39. Let* $\|\cdot\|$ *and* $\langle\,\triangle\mid\triangledown\,\rangle$ *be as in Definition 2.2 page 15.*

| | | | | | | |
|---|---|---|---|---|---|---|
| 1. | $\|\mathbf{T}f\|$ | $=$ | $\|\mathbf{D}f\|$ | $=$ | $\|f\|$ | $\forall f \in L^2_{\mathbb{R}}$   *(ISOMETRIC IN LENGTH)* |
| 2. | $\|\mathbf{T}f - \mathbf{T}g\|$ | $=$ | $\|\mathbf{D}f - \mathbf{D}g\|$ | $=$ | $\|f - g\|$ | $\forall f,g \in L^2_{\mathbb{R}}$   *(ISOMETRIC IN DISTANCE)* |
| 3. | $\|\mathbf{T}^{-1}f - \mathbf{T}^{-1}g\|$ | $=$ | $\|\mathbf{D}^{-1}f - \mathbf{D}^{-1}g\|$ | $=$ | $\|f - g\|$ | $\forall f,g \in L^2_{\mathbb{R}}$   *(ISOMETRIC IN DISTANCE)* |
| 4. | $\langle \mathbf{T}f \mid \mathbf{T}g \rangle$ | $=$ | $\langle \mathbf{D}f \mid \mathbf{D}g \rangle$ | $=$ | $\langle f \mid g \rangle$ | $\forall f,g \in L^2_{\mathbb{R}}$   *(SURJECTIVE)* |
| 5. | $\langle \mathbf{T}^{-1}f \mid \mathbf{T}^{-1}g \rangle$ | $=$ | $\langle \mathbf{D}^{-1}f \mid \mathbf{D}^{-1}g \rangle$ | $=$ | $\langle f \mid g \rangle$ | $\forall f,g \in L^2_{\mathbb{R}}$   *(SURJECTIVE)* |

✎**PROOF:**   These results follow directly from the fact that $\mathbf{T}$ and $\mathbf{D}$ are *unitary* (Proposition 4.10 page 43) and from Theorem 1.50 page 14 and Theorem 1.51 page 14.   ✐

**Proposition 4.16**  *Let* $\mathbf{T}$ *be as in Definition 4.1 page 38. Let* $\mathbf{A}^*$ *be the adjoint* (Proposition 1.33 page 10) *of an operator* $\mathbf{A}$.

$$\left( \sum_{n\in\mathbb{Z}} \mathbf{T}^n \right) = \left( \sum_{n\in\mathbb{Z}} \mathbf{T}^n \right)^* \qquad \left( \textit{The operator } \left[ \sum_{n\in\mathbb{Z}} \mathbf{T}^n \right] \textit{ is SELF-ADJOINT} \right)$$

✎**PROOF:**

$$\left\langle \left( \sum_{n\in\mathbb{Z}} \mathbf{T}^n \right) f(x) \,\middle|\, g(x) \right\rangle = \left\langle \sum_{n\in\mathbb{Z}} f(x-n) \,\middle|\, g(x) \right\rangle \qquad \text{by definition of } \mathbf{T} \text{ (Definition 4.1 page 38)}$$

$$= \left\langle \sum_{n\in\mathbb{Z}} f(x+n) \,\middle|\, g(x) \right\rangle \qquad \text{by } \textit{commutative} \text{ property of addition}$$

$$= \sum_{n\in\mathbb{Z}} \langle f(x+n) \mid g(x) \rangle \qquad \text{by } \textit{additive} \text{ property of } \langle\,\triangle\mid\triangledown\,\rangle$$

$$= \sum_{n\in\mathbb{Z}} \langle f(u) \mid g(u-n) \rangle \qquad \text{where } u = x + n$$

$$= \left\langle f(u) \,\middle|\, \sum_{n\in\mathbb{Z}} g(u-n) \right\rangle \qquad \text{by } \textit{additive} \text{ property of } \langle\,\triangle\mid\triangledown\,\rangle$$

$$= \left\langle f(x) \,\middle|\, \sum_{n\in\mathbb{Z}} g(x-n) \right\rangle \qquad \text{by change of dummy variable: } t \to u$$

$$= \left\langle f(x) \,\middle|\, \sum_{n\in\mathbb{Z}} \mathbf{T}^n g(x) \right\rangle \qquad \text{by definition of } \mathbf{T} \text{ (Definition 4.1 page 38)}$$

$$\iff \left( \sum_{n\in\mathbb{Z}} \mathbf{T}^n \right) = \left( \sum_{n\in\mathbb{Z}} \mathbf{T}^n \right)^* \qquad \text{by definition of operator adjoint (page 10)}$$

$$\iff \left( \sum_{n\in\mathbb{Z}} \mathbf{T}^n \right) \text{ is } \textit{self-adjoint} \qquad \text{by definition of } \textit{self-adjoint} \text{ (Definition 1.37 page 11)}$$

✐

### 4.2.5  Fourier transform properties

**Proposition 4.17** *Let* **T** *and* **D** *be as in Definition 4.1 page 38. Let* **B** *be the* TWO-SIDED LAPLACE TRANSFORM *defined as*

$$[\mathbf{B}f](s) \triangleq \frac{1}{\sqrt{2\pi}} \int_{\mathbb{R}} f(x) e^{-sx} \, dx \ .$$

*1.*     $\mathbf{B}\mathbf{T}^n \ = \ e^{-sn}\mathbf{B}$          $\forall n \in \mathbb{Z}$

*2.*     $\mathbf{B}\mathbf{D}^j \ = \ \mathbf{D}^{-j}\mathbf{B}$          $\forall j \in \mathbb{Z}$

*3.*     $\mathbf{D}\mathbf{B} \ = \ \mathbf{B}\mathbf{D}^{-1}$          $\forall n \in \mathbb{Z}$

*4.*   $\mathbf{B}\mathbf{D}^{-1}\mathbf{B}^{-1} \ = \ \mathbf{B}^{-1}\mathbf{D}^{-1}\mathbf{B} \ = \ \mathbf{D}$   $\forall n \in \mathbb{Z}$   ($\mathbf{D}^{-1}$ *is* SIMILAR *to* $\mathbf{D}$)

*5.*     $\mathbf{D}\mathbf{B}\mathbf{D} \ = \ \mathbf{D}^{-1}\mathbf{B}\mathbf{D}^{-1} \ = \ \mathbf{B}$   $\forall n \in \mathbb{Z}$

PROOF:

$\mathbf{B}\mathbf{T}^n f(x) = \mathbf{B}f(x-n)$          by definition of **T** (Definition 4.1 page 38)

$\qquad = \frac{1}{\sqrt{2\pi}} \int_{\mathbb{R}} f(x-n) e^{-sx} \, dx$          by definition of **B**

$\qquad = \frac{1}{\sqrt{2\pi}} \int_{\mathbb{R}} f(u) e^{-s(u+n)} \, du$          where $u \triangleq x - n$

$\qquad = e^{-sn} \left[ \frac{1}{\sqrt{2\pi}} \int_{\mathbb{R}} f(u) e^{-su} \, du \right]$

$\qquad = e^{-sn} \, \mathbf{B}f(x)$          by definition of **B**

$\mathbf{B}\mathbf{D}^j f(x) = \mathbf{B}\left[ 2^{j/2} f\left( 2^j x \right) \right]$          by definition of **D** (Definition 4.1 page 38)

$\qquad = \frac{1}{\sqrt{2\pi}} \int_{\mathbb{R}} \left[ 2^{j/2} f\left( 2^j x \right) \right] e^{-sx} \, dx$          by definition of **B**

$\qquad = \frac{1}{\sqrt{2\pi}} \int_{\mathbb{R}} \left[ 2^{j/2} f(u) \right] e^{-s2^{-j}} 2^{-j} \, du$          let $u \triangleq 2^j x \implies x = 2^{-j} u$

$\qquad = \frac{\sqrt{2}}{2} \frac{1}{\sqrt{2\pi}} \int_{\mathbb{R}} f(u) e^{-s2^{-j}u} \, du$

$\qquad = \mathbf{D}^{-1} \left[ \frac{1}{\sqrt{2\pi}} \int_{\mathbb{R}} f(u) e^{-su} \, du \right]$          by Proposition 4.9 page 42 and Proposition 4.10 page 43

$\qquad = \mathbf{D}^{-j} \, \mathbf{B}f(x)$          by definition of **B**

$\mathbf{D}\mathbf{B} f(x) = \mathbf{D} \left[ \frac{1}{\sqrt{2\pi}} \int_{\mathbb{R}} f(x) e^{-sx} \, dx \right]$          by definition of **B**

$\qquad = \frac{\sqrt{2}}{\sqrt{2\pi}} \int_{\mathbb{R}} f(x) e^{-2sx} \, dx$          by definition of **D** (Definition 4.1 page 38)

$\qquad = \frac{\sqrt{2}}{\sqrt{2\pi}} \int_{\mathbb{R}} f\left( \frac{u}{2} \right) e^{-su} \, \frac{1}{2} \, du$          let $u \triangleq 2x \implies x = \frac{1}{2}u$

$$= \frac{1}{\sqrt{2\pi}} \int_{\mathbb{R}} \left[ \frac{\sqrt{2}}{2} f\left(\frac{u}{2}\right) \right] e^{-su} \, du$$

$$= \frac{1}{\sqrt{2\pi}} \int_{\mathbb{R}} \left[ \mathbf{D}^{-1} f \right](u) \, e^{-su} \, du \qquad \text{by Proposition 4.9 page 42 and Proposition 4.10 page 43}$$

$$= \mathbf{B}\mathbf{D}^{-1} f(x) \qquad \text{by definition of } \mathbf{B}$$

$$\mathbf{B}^{-1}\mathbf{D}^{-1}\mathbf{B} = \mathbf{B}^{-1}\mathbf{B}\mathbf{D} \qquad \text{by previous result}$$

$$= \mathbf{D} \qquad \text{by definition of operator inverse (Definition 1.15 page 6)}$$

$$\mathbf{B}\mathbf{D}^{-1}\mathbf{B}^{-1} = \mathbf{D}\mathbf{B}\mathbf{B}^{-1} \qquad \text{by previous result}$$

$$= \mathbf{D} \qquad \text{by definition of operator inverse (Definition 1.15 page 6)}$$

$$\mathbf{D}\mathbf{B}\mathbf{D} = \mathbf{D}\mathbf{D}^{-1}\mathbf{B} \qquad \text{by previous result}$$

$$= \mathbf{B} \qquad \text{by definition of operator inverse (Definition 1.15 page 6)}$$

$$\mathbf{D}^{-1}\mathbf{B}\mathbf{D}^{-1} = \mathbf{D}^{-1}\mathbf{D}\mathbf{B} \qquad \text{by previous result}$$

$$= \mathbf{B} \qquad \text{by definition of operator inverse (Definition 1.15 page 6)}$$

✑

**Corollary 4.18**  *Let* $\mathbf{T}$ *and* $\mathbf{D}$ *be as in Definition 4.1 page 38. Let* $\tilde{f}(\omega) \triangleq \tilde{\mathbf{F}}f(x)$ *be the* FOURIER TRANSFORM *(Definition 2.22 page 20) of some function* $f \in L^2_{\mathbb{R}}$ *(Definition 2.2 page 15).*

1.  $\tilde{\mathbf{F}}\mathbf{T}^n = e^{-i\omega n}\tilde{\mathbf{F}}$
2.  $\tilde{\mathbf{F}}\mathbf{D}^j = \mathbf{D}^{-j}\tilde{\mathbf{F}}$
3.  $\mathbf{D}\tilde{\mathbf{F}} = \tilde{\mathbf{F}}\mathbf{D}^{-1}$
4.  $\mathbf{D} = \tilde{\mathbf{F}}\mathbf{D}^{-1}\tilde{\mathbf{F}}^{-1} = \tilde{\mathbf{F}}^{-1}\mathbf{D}^{-1}\tilde{\mathbf{F}}$
5.  $\tilde{\mathbf{F}} = \mathbf{D}\tilde{\mathbf{F}}\mathbf{D} = \mathbf{D}^{-1}\tilde{\mathbf{F}}\mathbf{D}^{-1}$

✎PROOF:   These results follow directly from Proposition 4.17 page 46.

$$\tilde{\mathbf{F}} = \mathbf{B}|_{s=i\omega}$$

✑

**Proposition 4.19**  *Let* $\mathbf{T}$ *and* $\mathbf{D}$ *be as in Definition 4.1 page 38. Let* $\tilde{f}(\omega) \triangleq \tilde{\mathbf{F}}f(x)$ *be the* FOURIER TRANSFORM *(Definition 2.22 page 20) of some function* $f \in L^2_{\mathbb{R}}$ *(Definition 2.2 page 15).*

$$\tilde{\mathbf{F}}\mathbf{D}^j\mathbf{T}^nf(x) = \frac{1}{2^{j/2}}e^{-i\frac{\omega}{2^j}n}\tilde{f}\left(\frac{\omega}{2^j}\right)$$

✎PROOF:

$$\tilde{\mathbf{F}}\mathbf{D}^j\mathbf{T}^nf(x) = \mathbf{D}^{-j}\tilde{\mathbf{F}}\mathbf{T}^nf(x) \qquad \text{by Corollary 4.18 page 47 (3)}$$

$$= \mathbf{D}^{-j}e^{-i\omega n}\tilde{\mathbf{F}}f(x) \qquad \text{by Corollary 4.18 page 47 (3)}$$

$$= \mathbf{D}^{-j}e^{-i\omega n}\tilde{f}(\omega)$$

$$= 2^{-j/2}e^{-i2^{-j}\omega n}\tilde{f}\left(2^{-j}\omega\right) \qquad \text{by Proposition 4.3 page 39}$$

✑

**Proposition 4.20**  *Let* **T** *be the translation operator* (Definition *4.1 page 38*). *Let* $\tilde{\mathsf{f}}(\omega) \triangleq \hat{\mathbf{F}}\mathsf{f}(x)$ *be the* FOURIER TRANSFORM (Definition *2.22 page 20*) *of a function* $\mathsf{f} \in \boldsymbol{L}_{\mathbb{R}}^2$. *Let* $\breve{\mathsf{a}}(\omega)$ *be the* DTFT (Definition *2.38 page 24*) *of a sequence* $(a_n)_{n \in \mathbb{Z}} \in \boldsymbol{\ell}_{\mathbb{R}}^2$ (Definition *2.33 page 23*).

$$\tilde{\mathbf{F}} \sum_{n \in \mathbb{Z}} a_n \mathbf{T}^n \phi(x) = \breve{\mathsf{a}}(\omega)\tilde{\phi}(\omega) \qquad \forall (a_n) \in \boldsymbol{\ell}_{\mathbb{R}}^2, \phi(x) \in \boldsymbol{L}_{\mathbb{R}}^2$$

✎PROOF:

$$\tilde{\mathbf{F}} \sum_{n \in \mathbb{Z}} a_n \mathbf{T}^n \phi(x) = \sum_{n \in \mathbb{Z}} a_n \tilde{\mathbf{F}}\mathbf{T}^n \phi(x)$$

$$= \sum_{n \in \mathbb{Z}} a_n e^{-i\omega n} \tilde{\mathbf{F}} \phi(x) \qquad \text{by Corollary 4.18 page 47}$$

$$= \left[ \sum_{n \in \mathbb{Z}} a_n e^{-i\omega n} \right] \tilde{\phi}(\omega) \qquad \text{by definition of } \tilde{\phi}(\omega)$$

$$= \breve{\mathsf{a}}(\omega)\tilde{\phi}(\omega) \qquad \text{by definition of } DTFT \text{ (Definition 2.38 page 24)}$$

✎

**Theorem 4.21**  (Poisson Summation Formula—PSF) [111] *Let* $\tilde{\mathsf{f}}(\omega)$ *be the* FOURIER TRANSFORM (Definition *2.22 page 20*) *of a function* $\mathsf{f}(x) \in \boldsymbol{L}_{\mathbb{R}}^2$.

$$\underbrace{\sum_{n \in \mathbb{Z}} \mathbf{T}_\tau^n \mathsf{f}(x) = \sum_{n \in \mathbb{Z}} \mathsf{f}(x + n\tau)}_{\text{summation in "time"}} = \underbrace{\sqrt{\frac{2\pi}{\tau}} \, \hat{\mathbf{F}}^{-1} \mathbf{S} \tilde{\mathbf{F}}[\mathsf{f}(x)]}_{\text{operator notation}} = \underbrace{\frac{\sqrt{2\pi}}{\tau} \sum_{n \in \mathbb{Z}} \tilde{\mathsf{f}}\left(\frac{2\pi}{\tau}n\right) e^{i\frac{2\pi}{\tau}nx}}_{\text{summation in "frequency"}}$$

*where* $\mathbf{S} \in \boldsymbol{\ell}_{\mathbb{R}}^{2\,\boldsymbol{L}_{\mathbb{R}}^2}$ *is the* SAMPLING OPERATOR *defined as*

$$[\mathbf{S}\mathsf{f}(x)](n) \triangleq \mathsf{f}\left(\frac{2\pi}{\tau}n\right) \qquad \forall \mathsf{f} \in \boldsymbol{L}_{(\mathbb{R},\mathscr{B},\mu)}^2, \, \tau \in \mathbb{R}^+$$

✎PROOF:

(1)  lemma: If $\mathsf{h}(x) \triangleq \sum_{n \in \mathbb{Z}} \mathsf{f}(x + n\tau)$ then $\mathsf{h} \equiv \hat{\mathbf{F}}^{-1}\hat{\mathbf{F}}\,\mathsf{h}$. Proof:

Note that $\mathsf{h}(x)$ is periodic with period $\tau$. Because $\mathsf{h}$ is periodic, it is in the domain of $\hat{\mathbf{F}}$ and thus $\mathsf{h} \equiv \hat{\mathbf{F}}^{-1}\hat{\mathbf{F}}\,\mathsf{h}$.

---

[111] ✎ [5], page 624, ✎ [83], page 389, ✎ [89], page 254, ✎ [111], pages 194–195, ✎ [37], page 337

(2) Proof of PSF (this theorem—Theorem 4.21):

$$\sum_{n\in\mathbb{Z}} f(x+n\tau) = \hat{\mathbf{F}}^{-1}\hat{\mathbf{F}}\sum_{n\in\mathbb{Z}} f(x+n\tau) \qquad \text{by item 1 page 48}$$

$$= \hat{\mathbf{F}}^{-1}\underbrace{\left[\frac{1}{\sqrt{\tau}}\int_0^\tau \left(\sum_{n\in\mathbb{Z}} f(x+n\tau)\right)e^{-i\frac{2\pi}{\tau}kx}\,\mathrm{d}x\right]}_{\hat{\mathbf{F}}\left[\sum_{n\in\mathbb{Z}} f(x+n\tau)\right]} \qquad \text{by def. of } \hat{\mathbf{F}} \text{ (Definition 2.16 page 18)}$$

$$= \hat{\mathbf{F}}^{-1}\left[\frac{1}{\sqrt{\tau}}\sum_{n\in\mathbb{Z}}\int_0^\tau f(x+n\tau)e^{-i\frac{2\pi}{\tau}kx}\,\mathrm{d}x\right]$$

$$= \hat{\mathbf{F}}^{-1}\left[\frac{1}{\sqrt{\tau}}\sum_{n\in\mathbb{Z}}\int_{u=n\tau}^{u=(n+1)\tau} f(u)e^{-i\frac{2\pi}{\tau}k(u-n\tau)}\,\mathrm{d}u\right] \qquad \text{where } u \triangleq x+n\tau \implies x=u-n\tau$$

$$= \hat{\mathbf{F}}^{-1}\left[\frac{1}{\sqrt{\tau}}\sum_{n\in\mathbb{Z}} e^{i2\pi kn}\int_{u=n\tau}^{u=(n+1)\tau} f(u)e^{-i\frac{2\pi}{\tau}ku}\,\mathrm{d}u\right]$$

$$= \sqrt{\frac{2\pi}{\tau}}\hat{\mathbf{F}}^{-1}\underbrace{\left[\frac{1}{\sqrt{2\pi}}\int_{u\in\mathbb{R}} f(u)e^{-i\left(\frac{2\pi}{\tau}k\right)u}\,\mathrm{d}u\right]}_{[\bar{\mathbf{F}}f]\left(\frac{2\pi}{\tau}k\right)} \qquad \text{by Theorem 2.17 page 18}$$

$$= \sqrt{\frac{2\pi}{\tau}}\hat{\mathbf{F}}^{-1}\left[\left[\bar{\mathbf{F}}f(x)\right]\left(\frac{2\pi}{\tau}k\right)\right] \qquad \text{by definition of } \tilde{\mathbf{F}} \text{ (page 20)}$$

$$= \sqrt{\frac{2\pi}{\tau}}\hat{\mathbf{F}}^{-1}\mathbf{S}\bar{\mathbf{F}}f \qquad \text{by definition of } \mathbf{S}$$

$$= \frac{\sqrt{2\pi}}{\tau}\sum_{n\in\mathbb{Z}}\tilde{f}\left(\frac{2\pi}{\tau}n\right)e^{i\frac{2\pi}{\tau}nx} \qquad \text{by Theorem 2.17 page 18}$$

**Theorem 4.22** (Inverse Poisson Summation Formula—IPSF) [112] *Let $\tilde{f}(\omega)$ be the* FOURIER TRANSFORM *(Definition 2.22 page 20) of a function* $f(x) \in \boldsymbol{L}_\mathbb{R}^2$.

$$\underbrace{\sum_{n\in\mathbb{Z}}\mathbf{T}_{2\pi/\tau}^n\tilde{f}(\omega) \triangleq \sum_{n\in\mathbb{Z}}\tilde{f}\left(\omega-\frac{2\pi}{\tau}n\right)}_{\textit{summation in "frequency"}} \quad = \quad \underbrace{\frac{\tau}{\sqrt{2\pi}}\sum_{n\in\mathbb{Z}}f(n\tau)e^{-i\omega n\tau}}_{\textit{summation in "time"}}$$

✎ PROOF:

---

[112] ☞ [45], page 88

(1) **lemma**: If $h(\omega) \triangleq \sum_{n \in \mathbb{Z}} \tilde{f}\left(\omega + \frac{2\pi}{\tau}n\right)$, then $h \equiv \hat{\mathbf{F}}^{-1}\hat{\mathbf{F}}\,h$. Proof:
   Note that $h(\omega)$ is periodic with period $2\pi/T$:

$$h\left(\omega + \frac{2\pi}{\tau}\right) \triangleq \sum_{n \in \mathbb{Z}} \tilde{f}\left(\omega + \frac{2\pi}{\tau} + \frac{2\pi}{\tau}n\right) = \sum_{n \in \mathbb{Z}} \tilde{f}\left(\omega + (n+1)\frac{2\pi}{\tau}\right) = \sum_{n \in \mathbb{Z}} \tilde{f}\left(\omega + \frac{2\pi}{\tau}n\right) \triangleq h(\omega)$$

   Because $h$ is periodic, it is in the domain of $\hat{\mathbf{F}}$ and is equivalent to $\hat{\mathbf{F}}^{-1}\hat{\mathbf{F}}\,h$.

(2) Proof of IPSF (this theorem—Theorem 4.22):

$$\sum_{n \in \mathbb{Z}} \tilde{f}\left(\omega + \frac{2\pi}{\tau}n\right)$$

$$= \hat{\mathbf{F}}^{-1}\hat{\mathbf{F}} \sum_{n \in \mathbb{Z}} \tilde{f}\left(\omega + \frac{2\pi}{\tau}n\right) \qquad \text{by item 1 page 50}$$

$$= \hat{\mathbf{F}}^{-1}\underbrace{\left[\sqrt{\frac{\tau}{2\pi}} \int_0^{\frac{2\pi}{\tau}} \sum_{n \in \mathbb{Z}} \tilde{f}\left(\omega + \frac{2\pi}{\tau}n\right) e^{-i\omega \frac{2\pi}{2\pi/\tau}k}\, d\omega\right]}_{\hat{\mathbf{F}}\left[\sum_{n \in \mathbb{Z}} \tilde{f}\left(\omega + \frac{2\pi}{\tau}n\right)\right]} \qquad \text{by definition of } \hat{\mathbf{F}} \text{ page 18}$$

$$= \hat{\mathbf{F}}^{-1}\left[\sqrt{\frac{\tau}{2\pi}} \sum_{n \in \mathbb{Z}} \int_0^{\frac{2\pi}{\tau}} \tilde{f}\left(\omega + \frac{2\pi}{\tau}n\right) e^{-i\omega Tk}\, d\omega\right]$$

$$= \hat{\mathbf{F}}^{-1}\left[\sqrt{\frac{\tau}{2\pi}} \sum_{n \in \mathbb{Z}} \int_{u=\frac{2\pi}{\tau}n}^{u=\frac{2\pi}{\tau}(n+1)} \tilde{f}(u) e^{-i\left(u-\frac{2\pi}{\tau}n\right)Tk}\, du\right] \qquad \text{where } u \triangleq \omega + \frac{2\pi}{\tau}n \implies \omega = u - \frac{2\pi}{\tau}n$$

$$= \hat{\mathbf{F}}^{-1}\left[\sqrt{\frac{\tau}{2\pi}} \sum_{n \in \mathbb{Z}} e^{i2\pi nk} \int_{\frac{2\pi}{\tau}n}^{\frac{2\pi}{\tau}(n+1)} \tilde{f}(u) e^{-iu\tau k}\, du\right]$$

$$= \hat{\mathbf{F}}^{-1}\left[\sqrt{\frac{\tau}{2\pi}} \int_{\mathbb{R}} \tilde{f}(u) e^{-iu\tau k}\, du\right]$$

$$= \sqrt{\tau}\, \hat{\mathbf{F}}^{-1}\underbrace{\left[\frac{1}{\sqrt{2\pi}} \int_{\mathbb{R}} \tilde{f}(u) e^{iu(-\tau k)}\, du\right]}_{[\hat{\mathbf{F}}^{-1}\tilde{f}](-k\tau)}$$

$$= \sqrt{\tau}\, \hat{\mathbf{F}}^{-1}\left[[\hat{\mathbf{F}}^{-1}\tilde{f}](-k\tau)\right] \qquad \text{by value of } \tilde{\mathbf{F}}^{-1} \text{ (Theorem 2.24 page 20)}$$

$$= \sqrt{\tau}\, \hat{\mathbf{F}}^{-1}\mathbf{S}\hat{\mathbf{F}}^{-1}\,\tilde{f} \qquad \text{by definition of } \mathbf{S}$$

$$= \sqrt{\tau}\, \hat{\mathbf{F}}^{-1}\mathbf{S}f(x) \qquad \text{by definition of } \tilde{\mathbf{F}} \text{ (Definition 2.22 page 20)}$$

$$= \sqrt{\tau}\, \hat{\mathbf{F}}^{-1}f(-k\tau) \qquad \text{by definition of } \mathbf{S}$$

$$= \sqrt{\tau}\, \frac{1}{\sqrt{\frac{2\pi}{\tau}}} \sum_{k \in \mathbb{Z}} f(-k\tau) e^{i2\pi \frac{1}{2\pi} k\omega} \qquad \text{by definition of } \hat{\mathbf{F}}^{-1} \text{ (Theorem 2.17 page 18)}$$

$$= \frac{\tau}{\sqrt{\frac{2\pi}{\tau}}} \sum_{k \in \mathbb{Z}} f(-k\tau) e^{ik\tau\omega} \qquad \text{by definition of } \hat{\mathbf{F}}^{-1} \text{ (Theorem 2.17 page 18)}$$

$$= \frac{\tau}{\sqrt{2\pi}} \sum_{m \in \mathbb{Z}} \mathsf{f}(m\tau) e^{-i\omega m\tau} \qquad\qquad \text{let } m \triangleq -k$$

**Remark 4.23** The left hand side of the *Poisson Summation Formula* (Theorem 4.21 page 48) is very similar to the *Zak Transform* **Z**: [113]

$$(\mathbf{Z}\mathsf{f})(t, \omega) \triangleq \sum_{n \in \mathbb{Z}} \mathsf{f}(x + n\tau) e^{i2\pi n\omega}$$

**Remark 4.24** A generalization of the *Poisson Summation Formula* (Theorem 4.21 page 48) is the **Selberg Trace Formula**.[114]

### 4.2.6 Basis theory properties

Definition 4.25 and Definition 4.26 define four quantities. In this text, the quantities' notation and terminology are similar to that used in the study of *random processes*.

**Definition 4.25** [115] Let $\langle \triangle \mid \triangledown \rangle$ be the *standard inner product* in $\mathbf{L}^2_{\mathbb{R}}$ (Definition 2.2 page 15).

$$
\begin{aligned}
R_{\mathsf{fg}}(n) &\triangleq \langle \mathsf{f}(x) \mid \mathbf{T}^n \mathsf{g}(x) \rangle, & \mathsf{f}, \mathsf{g} \in L^2_{\mathbb{F}}, & \quad\text{is the } \textbf{cross-correlation function} \text{ of } \mathsf{f} \text{ and } \mathsf{g}. \\
R_{\mathsf{ff}}(n) &\triangleq \langle \mathsf{f}(x) \mid \mathbf{T}^n \mathsf{f}(x) \rangle, & \mathsf{f} \in L^2_{\mathbb{F}}, & \quad\text{is the } \textbf{autocorrelation function} \text{ of } \mathsf{f}.
\end{aligned}
$$

**Definition 4.26** [116] Let $R_{\mathsf{fg}}(n)$ and $R_{\mathsf{ff}}(n)$ be the sequences defined in Definition 4.25 page 51. Let $\mathbf{Z}\langle\!\langle x_n \rangle\!\rangle$ be the *z-transform* (Definition 2.35 page 24) of a sequence $\langle\!\langle x_n \rangle\!\rangle_{n \in \mathbb{Z}}$.

$$
\begin{aligned}
S_{\mathsf{fg}}(z) &\triangleq \mathbf{Z}\Big[ R_{\mathsf{fg}}(n) \Big], & {}_{\mathsf{f}, \mathsf{g} \in L^2_{\mathbb{F}}}, & \quad\text{is the } \textbf{complex cross-power spectrum} \text{ of } \mathsf{f} \text{ and } \mathsf{g}. \\
S_{\mathsf{ff}}(z) &\triangleq \mathbf{Z}\Big[ R_{\mathsf{ff}}(n) \Big], & {}_{\mathsf{f}, \mathsf{g} \in L^2_{\mathbb{F}}}, & \quad\text{is the } \textbf{complex auto-power spectrum} \text{ of } \mathsf{f}.
\end{aligned}
$$

**Definition 4.27** [117] Let $S_{\mathsf{fg}}(z)$ and $S_{\mathsf{ff}}(z)$ be the functions defined in Definition 4.26 page 51.

$$
\begin{aligned}
\tilde{S}_{\mathsf{fg}}(\omega) &\triangleq S_{\mathsf{fg}}\big( e^{i\omega} \big), & {}_{\forall \mathsf{f}, \mathsf{g} \in L^2_{\mathbb{F}}}, & \quad\text{is the } \textbf{cross-power spectrum} \text{ of } \mathsf{f} \text{ and } \mathsf{g}. \\
\tilde{S}_{\mathsf{ff}}(\omega) &\triangleq S_{\mathsf{ff}}\big( e^{i\omega} \big), & {}_{\forall \mathsf{f} \in L^2_{\mathbb{F}}}, & \quad\text{is the } \textbf{auto-power spectrum} \text{ of } \mathsf{f}.
\end{aligned}
$$

**Theorem 4.28** [118] *Let $\tilde{\mathsf{f}}(\omega)$ be the* FOURIER TRANSFORM *(Definition 2.22 page 20) of a function* $\mathsf{f}(x) \in L^2_{\mathbb{F}}$.

$$
\begin{aligned}
\tilde{S}_{\mathsf{fg}}(\omega) &= 2\pi \sum_{n \in \mathbb{Z}} \tilde{\mathsf{f}}(\omega + 2\pi n)\tilde{\mathsf{g}}^*(\omega + 2\pi n) & {}_{\forall \mathsf{f}, \mathsf{g} \in L^2_{\mathbb{F}}} \\
\tilde{S}_{\mathsf{ff}}(\omega) &= 2\pi \sum_{n \in \mathbb{Z}} \big| \tilde{\mathsf{f}}(\omega + 2\pi n) \big|^2 & {}_{\forall \mathsf{f} \in L^2_{\mathbb{F}}}
\end{aligned}
$$

---

✎ PROOF:  Let $z \triangleq e^{i\omega}$.

$$\tilde{S}_{fg}(\omega) \triangleq S_{fg}(z) \qquad\qquad \text{by definition of } \tilde{S}_{fg} \text{ Definition 4.27 page 51}$$

$$= \sum_{n \in \mathbb{Z}} R_{fg}(n) z^{-n} \qquad\qquad \text{by definition of } S_{fg} \text{ Definition 4.26 page 51}$$

$$= \sum_{n \in \mathbb{Z}} \langle f(x) \mid g(x-n) \rangle \, z^{-n} \qquad\qquad \text{by Definition 4.27 page 51}$$

$$= \sum_{n \in \mathbb{Z}} \langle \tilde{F}[f(x)] \mid \tilde{F}[g(x-n)] \rangle \, z^{-n} \qquad\qquad \text{by } \textit{unitary} \text{ property of } \tilde{F} \text{ (Theorem 2.27 page 21)}$$

$$= \sum_{n \in \mathbb{Z}} \langle \tilde{f}(v) \mid e^{-ivn} \tilde{g}(v) \rangle \, z^{-n} \qquad\qquad \text{by Theorem 2.28 page 22}$$

$$= \sum_{n \in \mathbb{Z}} \sqrt{2\pi} \left[ \frac{1}{\sqrt{2\pi}} \int_{\mathbb{R}} \tilde{f}(v) \tilde{g}^*(v) e^{ivu} \, dv \right]_{u=n} z^{-n} \qquad\qquad \text{by Definition 2.2 page 15}$$

$$= \sqrt{2\pi} \sum_{n \in \mathbb{Z}} \left[ \tilde{F}^{-1}\left( \sqrt{2\pi} \tilde{f}(v) \tilde{g}^*(v) \right) \right]_{u=n} e^{-i\omega n} \qquad\qquad \text{by Theorem 2.24 page 20}$$

$$= 2\pi \sum_{n \in \mathbb{Z}} \tilde{f}(\omega+2\pi n) \tilde{g}^*(\omega+2\pi n) \qquad\qquad \text{by IPSF (Theorem 4.22 page 49)}, \ \tau = 1$$

$$\tilde{S}_{ff}(\omega) = \tilde{S}_{fg}(\omega)\big|_{g=f} \qquad\qquad \text{by Definition 4.27 page 51}$$

$$= 2\pi \sum_{n \in \mathbb{Z}} \tilde{f}(\omega+2\pi n) \tilde{g}^*(\omega+2\pi n)\bigg|_{g=f} \qquad\qquad \text{by previous result}$$

$$= 2\pi \sum_{n \in \mathbb{Z}} \tilde{f}(\omega+2\pi n) \tilde{f}^*(\omega+2\pi n)$$

$$= 2\pi \sum_{n \in \mathbb{Z}} \left| \tilde{f}(\omega+2\pi n) \right|^2$$

✌

**Theorem 4.29** [119] *Let* $\tilde{S}_{\phi\phi}$ *be the* SPECTRAL DENSITY FUNCTION *(Definition 4.27 page 51) of a function* $\phi(x) \in L^2_{\mathbb{R}}$. *Let* $0 < A < B$. *Let* $\|\cdot\|$ *be defined as in Definition 2.2 page 15.*

$$\underbrace{\left\{ A \sum_{n \in \mathbb{N}} |a_n|^2 \le \|a_n \mathbf{T}^n \phi(x)\|^2 \le B \sum_{n \in \mathbb{N}} |a_n|^2 \quad \forall_{(a_n) \in \ell^2_F} \right\}}_{(\mathbf{T}^n \phi) \text{ is a RIESZ BASIS (Theorem 3.22 page 37) for } L^2_{\mathbb{R}}} \iff \left\{ A \le \tilde{S}_{\phi\phi}(\omega) \le B \right\}$$

✎ PROOF:

(1) lemma:

$$\left\lVert \sum_{n\in\mathbb{Z}} a_n \mathbf{T}^n \phi(x) \right\rVert^2 = \left\lVert \tilde{\mathbf{F}} \sum_{n\in\mathbb{Z}} a_n \mathbf{T}^n \phi(x) \right\rVert^2 \qquad \text{because } \tilde{\mathbf{F}} \text{ is } \textit{unitary} \text{ (Theorem 2.25 page 20)}$$

$$= \left\lVert \breve{\mathtt{a}}(\omega)\tilde{\phi}(\omega) \right\rVert^2 \qquad \text{by Proposition 4.20 page 48}$$

$$= \int_{\mathbb{R}} \left\lvert \breve{\mathtt{a}}(\omega)\tilde{\phi}(\omega) \right\rvert^2 \, d\omega \qquad \text{by definition of } \lVert\cdot\rVert$$

$$= \sum_{n\in\mathbb{Z}} \int_0^{2\pi} \left\lvert \breve{\mathtt{a}}(\omega + 2\pi n)\tilde{\phi}(\omega + 2\pi n) \right\rvert^2 \, d\omega$$

$$= \int_0^{2\pi} \sum_{n\in\mathbb{Z}} \left\lvert \breve{\mathtt{a}}(\omega + 2\pi n) \right\rvert^2 \left\lvert \tilde{\phi}(\omega + 2\pi n) \right\rvert^2 \, d\omega$$

$$= \int_0^{2\pi} \sum_{n\in\mathbb{Z}} \left\lvert \breve{\mathtt{a}}(\omega) \right\rvert^2 \left\lvert \tilde{\phi}(\omega + 2\pi n) \right\rvert^2 \, d\omega \qquad \text{by Proposition 2.39 page 25}$$

$$= \int_0^{2\pi} \left\lvert \breve{\mathtt{a}}(\omega) \right\rvert^2 \frac{1}{2\pi} 2\pi \sum_{n\in\mathbb{Z}} \left\lvert \tilde{\phi}(\omega + 2\pi n) \right\rvert^2 \, d\omega$$

$$= \frac{1}{2\pi} \int_0^{2\pi} \left\lvert \breve{\mathtt{a}}(\omega) \right\rvert^2 \tilde{S}_{\phi\phi}(\omega) \, d\omega \qquad \text{by def. of } \tilde{S}_{\phi\phi}(\omega) \text{ (Theorem 4.28 page 51)}$$

(2) lemma:

$$\frac{1}{2\pi} \int_0^{2\pi} \lvert \breve{\mathtt{a}}(\omega) \rvert^2 \, d\omega = \frac{1}{2\pi} \int_0^{2\pi} \left\lvert \sum_{n\in\mathbb{Z}} a_n e^{-i\omega n} \right\rvert^2 \, d\omega \qquad \text{by def. of } \textit{DTFT} \text{ (Definition 2.38 page 24)}$$

$$= \frac{1}{2\pi} \int_0^{2\pi} \left[ \sum_{n\in\mathbb{Z}} a_n e^{-i\omega n} \right] \left[ \sum_{m\in\mathbb{Z}} a_m e^{-i\omega m} \right]^* \, d\omega$$

$$= \frac{1}{2\pi} \int_0^{2\pi} \left[ \sum_{n\in\mathbb{Z}} a_n e^{-i\omega n} \right] \left[ \sum_{m\in\mathbb{Z}} a_m^* e^{i\omega m} \right] \, d\omega$$

$$= \frac{1}{2\pi} \sum_{n\in\mathbb{Z}} \sum_{m\in\mathbb{Z}} a_n a_m^* \int_0^{2\pi} e^{-i\omega(n-m)} \, d\omega$$

$$= \frac{1}{2\pi} \sum_{n\in\mathbb{Z}} \sum_{m\in\mathbb{Z}} a_n a_m^* 2\pi \breve{\delta}_{nm}$$

$$= \sum_{n\in\mathbb{Z}} \lvert a_n \rvert^2 \qquad \text{by def. of } \breve{\delta} \text{ (Definition 3.12 page 35)}$$

(3) Proof for ( $\Longleftarrow$ ) case:

$$\boxed{A \sum_{n\in\mathbb{Z}} \lvert a_n \rvert^2} = \frac{A}{2\pi} \int_0^{2\pi} \lvert \breve{\mathtt{a}}(\omega) \rvert^2 \, d\omega \qquad \text{by lemma 2 page 53}$$

$$= \frac{1}{2\pi} \int_0^{2\pi} \lvert \breve{\mathtt{a}}(\omega) \rvert^2 A \, d\omega$$

$$\boxed{\leq} \frac{1}{2\pi} \int_0^{2\pi} |\breve{a}(\omega)|^2 \tilde{S}_{\phi\phi}(\omega) \, d\omega \qquad \text{by right hypothesis}$$

$$= \boxed{\left\| \sum_{n \in \mathbb{Z}} a_n \mathbf{T}^n \phi(x) \right\|^2} \qquad \text{by lemma 1 page 53}$$

$$= \frac{1}{2\pi} \int_0^{2\pi} |\breve{a}(\omega)|^2 \tilde{S}_{\phi\phi}(\omega) \, d\omega \qquad \text{by lemma 1 page 53}$$

$$\boxed{\leq} \frac{1}{2\pi} \int_0^{2\pi} |\breve{a}(\omega)|^2 B \, d\omega \qquad \text{by right hypothesis}$$

$$= \frac{B}{2\pi} \int_0^{2\pi} |\breve{a}(\omega)|^2 \, d\omega$$

$$= \boxed{B \sum_{n \in \mathbb{Z}} |a_n|^2} \qquad \text{by lemma 2 page 53}$$

(4) Proof for ( $\implies$ ) case:

    (a) Let $Y \triangleq \left\{ \omega \in [0, 2\pi] \,\middle|\, \tilde{S}_{\phi\phi}(\omega) > \alpha \right\}$
        and $X \triangleq \left\{ \omega \in [0, 2\pi] \,\middle|\, \tilde{S}_{\phi\phi}(\omega) < \alpha \right\}$

    (b) Let $\mathbb{1}_A(x)$ be the *set indicator* (Definition 2.3 page 15) of a set $A$.
        Let $(b_n)_{n \in \mathbb{Z}}$ be the *inverse DTFT* (Theorem 2.43 page 27) of $\mathbb{1}_Y(\omega)$ such that
        $\mathbb{1}_Y(\omega) \triangleq \sum_{n \in \mathbb{N}} b_n e^{-i\omega n} \triangleq \tilde{b}(\omega)$.
        Let $(a_n)_{n \in \mathbb{Z}}$ be the *inverse DTFT* (Theorem 2.43 page 27) of $\mathbb{1}_X(\omega)$ such that
        $\mathbb{1}_X(\omega) \triangleq \sum_{n \in \mathbb{N}} a_n e^{-i\omega n} \triangleq \breve{a}(\omega)$.

    (c) Proof that $\alpha \leq B$:
        Let $\mu(A)$ be the *measure* of a set $A$.

$$\boxed{B} \sum_{n \in \mathbb{Z}} |b_n|^2 \geq \left\| \sum_{n \in \mathbb{Z}} b_n \mathbf{T}^n \phi(x) \right\|^2 \qquad \text{by left hypothesis}$$

$$= \frac{1}{2\pi} \int_0^{2\pi} |\tilde{b}(\omega)|^2 \tilde{S}_{\phi\phi}(\omega) \, d\omega \qquad \text{by lemma 1 page 53}$$

$$= \frac{1}{2\pi} \int_0^{2\pi} |\mathbb{1}_Y(\omega)|^2 \tilde{S}_{\phi\phi}(\omega) \, d\omega \qquad \text{by definition of } \mathbb{1}_Y(\omega) \text{ (item 4b page 54)}$$

$$= \frac{1}{2\pi} \int_Y |1|^2 \tilde{S}_{\phi\phi}(\omega) \, d\omega \qquad \text{by definition of } \mathbb{1}_Y(\omega) \text{ (item 4b page 54)}$$

$$\boxed{\geq} \frac{\alpha}{2\pi} \mu(Y) \qquad \text{by definition of } Y \text{ (item 4a page 54)}$$

$$= \int_0^{2\pi} |\mathbb{1}_Y(\omega)|^2 \, d\omega \qquad \text{by definition of } \mathbb{1}_Y(\omega) \text{ (item 4b page 54)}$$

$$= \int_0^{2\pi} \left| \sum_{n \in \mathbb{Z}} b_n e^{-i\omega n} \right|^2 \, d\omega \qquad \text{by definition of } (b_n) \text{ (item 4b page 54)}$$

$$= \int_0^{2\pi} \left| \breve{\mathsf{b}}(\omega) \right|^2 \, \mathrm{d}\omega \qquad\qquad \text{by definition of } \breve{\mathsf{b}}(\omega) \text{ (item 4b page 54)}$$

$$= \boxed{\alpha} \sum_{n \in \mathbb{Z}} |b_n|^2 \qquad\qquad \text{by lemma 2 page 53}$$

(d)  Proof that $\tilde{S}_{\phi\phi}(\omega) \leq B$:

    (i)  $\tilde{S}_{\phi\phi}(\omega) > \alpha$ whenever $\omega \in Y$ (item 4a page 54).

    (ii)  But even then, $\alpha \leq B$ (item 4c page 54).

    (iii)  So, $\tilde{S}_{\phi\phi}(\omega) \leq B$.

(e)  Proof that $A \leq \alpha$:

Let $\mu(A)$ be the *measure* of a set $A$.

$$\boxed{A} \sum_{n \in \mathbb{Z}} |a_n|^2 \leq \left\| \sum_{n \in \mathbb{Z}} a_n \mathbf{T}^n \phi(x) \right\|^2 \qquad\qquad \text{by left hypothesis}$$

$$= \frac{1}{2\pi} \int_0^{2\pi} |\breve{\mathsf{a}}(\omega)|^2 \tilde{S}_{\phi\phi}(\omega) \, \mathrm{d}\omega \qquad\qquad \text{by lemma 1 page 53}$$

$$= \frac{1}{2\pi} \int_0^{2\pi} \left| \mathbb{1}_X(\omega) \right|^2 \tilde{S}_{\phi\phi}(\omega) \, \mathrm{d}\omega \qquad\qquad \text{by definition of } \mathbb{1}_X(\omega) \text{ (Definition 2.3 page 15)}$$

$$= \frac{1}{2\pi} \int_X |1|^2 \tilde{S}_{\phi\phi}(\omega) \, \mathrm{d}\omega \qquad\qquad \text{by definition of } \mathbb{1}_X(\omega) \text{ (Definition 2.3 page 15)}$$

$$\boxed{\leq} \frac{\alpha}{2\pi} \mu(X) \qquad\qquad \text{by definition of } X \text{ (item 4a page 54)}$$

$$= \int_0^{2\pi} \left| \mathbb{1}_X(\omega) \right|^2 \, \mathrm{d}\omega \qquad\qquad \text{by definition of } \mathbb{1}_X(\omega) \text{ (Definition 2.3 page 15)}$$

$$= \int_0^{2\pi} \left| \sum_{n \in \mathbb{Z}} a_n e^{-i\omega n} \right|^2 \, \mathrm{d}\omega \qquad\qquad \text{by definition of } \langle\!\langle a_n \rangle\!\rangle \text{ (lemma 2 page 53)}$$

$$= \int_0^{2\pi} |\breve{\mathsf{a}}(\omega)|^2 \, \mathrm{d}\omega \qquad\qquad \text{by definition of } \breve{\mathsf{a}}(\omega) \text{ (lemma 2 page 53)}$$

$$= \boxed{\alpha} \sum_{n \in \mathbb{Z}} |a_n|^2 \qquad\qquad \text{by lemma 2 page 53}$$

(f)  Proof that $A \leq \tilde{S}_{\phi\phi}(\omega)$:

    (i)  $\tilde{S}_{\phi\phi}(\omega) < \alpha$ whenever $\omega \in X$ (item 4a page 54).

    (ii)  But even then, $A \leq \alpha$ (item 4e page 55).

    (iii)  So, $A \leq \tilde{S}_{\phi\phi}(\omega)$.

In the case that $\mathsf{f}$ and $\mathsf{g}$ are *orthonormal*, the spectral density relations simplify considerably (next).

**Theorem 4.30** [120] *Let $\tilde{S}_{\mathrm{ff}}$ and $\tilde{S}_{\mathrm{fg}}$ be* SPECTRAL DENSITY FUNCTIONS *(Definition 4.27 page 51).*

$$\langle \mathsf{f}(x) \mid \mathbf{T}^n \mathsf{f}(x) \rangle = \bar{\delta}_n \quad \text{(ORTHONORMAL)} \quad \iff \quad \tilde{S}_{\mathrm{ff}}(\omega) = 1 \quad \forall \mathsf{f} \in L_{\mathbb{F}}^2$$

$$\langle \mathsf{f}(x) \mid \mathbf{T}^n \mathsf{g}(x) \rangle = 0 \quad \iff \quad \tilde{S}_{\mathrm{fg}}(\omega) = 0 \quad \forall \mathsf{f}, \mathsf{g} \in L_{\mathbb{F}}^2$$

✎ PROOF:

(1) Proof that $\langle \mathsf{f}(x) \mid \mathsf{f}(x-n) \rangle = \bar{\delta}_n \iff \tilde{S}_{\mathrm{ff}}(\omega) = 1$: This follows directly from Theorem 4.29 (page 52) with $A = B = 1$ (by *Parseval's Identity* Theorem 3.16 page 36 since $\{\mathbf{T}^n \mathsf{f}\}$ is *orthonormal*)

(2) Alternate proof that $\langle \mathsf{f}(x) \mid \mathsf{f}(x-n) \rangle = \bar{\delta}_n \implies \tilde{S}_{\mathrm{ff}}(\omega) = 1$:

$$\tilde{S}_{\mathrm{ff}}(\omega) = \sum_{n \in \mathbb{Z}} R_{\mathrm{ff}}(n) e^{-i\omega n} \qquad \text{by definition of } \tilde{S}_{\mathrm{fg}} \text{ (Definition 4.27 page 51)}$$

$$= \sum_{n \in \mathbb{Z}} \langle \mathsf{f}(x) \mid \mathsf{f}(x-n) \rangle e^{-i\omega n} \qquad \text{by definition of } R_{\mathrm{ff}} \text{ (Definition 4.25 page 51)}$$

$$= \sum_{n \in \mathbb{Z}} \bar{\delta}_n e^{-i\omega n} \qquad \text{by left hypothesis}$$

$$= 1 \qquad \text{by definition of } \bar{\delta} \text{ (Definition 3.12 page 35)}$$

(3) Alternate proof that $\langle \mathsf{f}(x) \mid \mathsf{f}(x-n) \rangle = \bar{\delta}_n \impliedby \tilde{S}_{\mathrm{ff}}(\omega) = 1$:

$$\langle \mathsf{f}(x) \mid \mathsf{f}(x-n) \rangle$$

$$= \langle \bar{\mathbf{F}} \mathsf{f}(x) \mid \bar{\mathbf{F}} \mathsf{f}(x-n) \rangle \qquad \text{by } unitary \text{ prop. of } \bar{\mathbf{F}} \text{ (Theorem 2.27 page 21)}$$

$$= \langle \tilde{\mathsf{f}}(\omega) \mid e^{-i\omega n} \tilde{\mathsf{f}}(\omega) \rangle \qquad \text{by } shift \text{ property of } \bar{\mathbf{F}} \text{ (Theorem 2.28 page 22)}$$

$$= \int_{\mathbb{R}} \tilde{\mathsf{f}}(\omega) e^{i\omega n} \tilde{\mathsf{f}}^*(\omega) \, \mathrm{d}\omega \qquad \text{by def. of } \langle \triangle \mid \triangledown \rangle \text{ (Definition 2.2 page 15)}$$

$$= \int_{\mathbb{R}} |\tilde{\mathsf{f}}(\omega)|^2 e^{i\omega n} \, \mathrm{d}\omega$$

$$= \sum_{n \in \mathbb{Z}} \int_{2\pi n}^{2\pi(n+1)} |\tilde{\mathsf{f}}(\omega)|^2 e^{i\omega n} \, \mathrm{d}\omega$$

$$= \sum_{n \in \mathbb{Z}} \int_0^{2\pi} |\tilde{\mathsf{f}}(u+2\pi n)|^2 e^{i(u+2\pi n)n} \, \mathrm{d}u \qquad \text{where } u \triangleq \omega - 2\pi n \implies \omega = u + 2\pi n$$

$$= \frac{1}{2\pi} \int_0^{2\pi} \left[ 2\pi \sum_{n \in \mathbb{Z}} |\tilde{\mathsf{f}}(u+2\pi n)|^2 \right] e^{iun} e^{i2\pi n n} \, \mathrm{d}u \qquad \overset{1}{\nearrow}$$

$$= \frac{1}{2\pi} \int_0^{2\pi} \tilde{S}_{\mathrm{ff}}(\omega) e^{iun} \, \mathrm{d}u \qquad \text{by Theorem 4.28 page 51}$$

$$= \frac{1}{2\pi} \int_0^{2\pi} e^{iun} \, \mathrm{d}u \qquad \text{by right hypothesis}$$

$$= \bar{\delta}_n \qquad \text{by definition of } \bar{\delta} \text{ (Definition 3.12 page 35)}$$

---

[120] ☞ [64], page 50, ⟨PROPOSITION 2.1.11⟩, ☞ [137], PAGE 23, ⟨COROLLARY 2.9⟩, ☞ [68], PAGES 214–215, ⟨LEMMA 9.2⟩, ☞ [106], PAGE 306, ⟨COROLLARY 6.4.9⟩

(4) Proof that $\langle \mathsf{f}(x) \mid \mathsf{g}(x-n) \rangle = 0 \implies \tilde{S}_{\mathsf{fg}}(\omega) = 0$:

$$\tilde{S}_{\mathsf{fg}}(\omega) = \sum_{n \in \mathbb{Z}} R_{\mathsf{fg}}(n) e^{-i\omega n} \qquad \text{by definition of } \tilde{S}_{\mathsf{fg}} \text{ (Definition 4.27 page 51)}$$

$$= \sum_{n \in \mathbb{Z}} \langle \mathsf{f}(x) \mid \mathsf{g}(x-n) \rangle e^{-i\omega n} \qquad \text{by definition of } R_{\mathsf{fg}} \text{ (Definition 4.25 page 51)}$$

$$= \sum_{n \in \mathbb{Z}} 0 e^{-i\omega n} \qquad \text{by left hypothesis}$$

$$= 0$$

(5) Proof that $\langle \mathsf{f}(x) \mid \mathsf{g}(x-n) \rangle = 0 \impliedby \tilde{S}_{\mathsf{fg}}(\omega) = 0$:

$$\langle \mathsf{f}(x) \mid \mathsf{g}(x-n) \rangle$$

$$= \langle \tilde{\mathbf{F}} \mathsf{f}(x) \mid \tilde{\mathbf{F}} \mathsf{g}(x-n) \rangle \qquad \text{by } \textit{unitary} \text{ property of } \tilde{\mathbf{F}} \text{ (Theorem 2.27 page 21)}$$

$$= \langle \tilde{\mathsf{f}}(\omega) \mid e^{-i\omega n} \tilde{\mathsf{g}}(\omega) \rangle \qquad \text{by } \textit{shift property} \text{ of } \tilde{\mathbf{F}} \text{ (Theorem 2.28 page 22)}$$

$$= \int_{\mathbb{R}} \tilde{\mathsf{f}}(\omega) e^{i\omega n} \tilde{\mathsf{g}}^*(\omega) \, \mathrm{d}\omega \qquad \text{by definition of } \langle \triangle \mid \triangledown \rangle \text{ (Definition 2.2 page 15)}$$

$$= \int_{\mathbb{R}} \tilde{\mathsf{f}}(\omega) \tilde{\mathsf{g}}^*(\omega) e^{i\omega n} \, \mathrm{d}\omega$$

$$= \sum_{n \in \mathbb{Z}} \int_{2\pi n}^{2\pi(n+1)} \tilde{\mathsf{f}}(\omega) \tilde{\mathsf{g}}^*(\omega) e^{i\omega n} \, \mathrm{d}\omega$$

$$= \sum_{n \in \mathbb{Z}} \int_{0}^{2\pi} \tilde{\mathsf{f}}(u+2\pi n) \tilde{\mathsf{g}}^*(u+2\pi n) e^{i(u+2\pi n)n} \, \mathrm{d}u \qquad \text{where } u \triangleq \omega - 2\pi n \implies \omega = u + 2\pi n$$

$$= \frac{1}{2\pi} \int_{0}^{2\pi} \left[ 2\pi \sum_{n \in \mathbb{Z}} \tilde{\mathsf{f}}(u+2\pi n) \tilde{\mathsf{g}}^*(u+2\pi n) \right] e^{iun} \underbrace{e^{i2\pi n n}}_{1} \, \mathrm{d}u$$

$$= \frac{1}{2\pi} \int_{0}^{2\pi} \tilde{S}_{\mathsf{fg}}(u) e^{iun} \, \mathrm{d}u \qquad \text{by Theorem 4.28 page 51}$$

$$= \frac{1}{2\pi} \int_{0}^{2\pi} 0 \cdot e^{iun} \, \mathrm{d}u \qquad \text{by right hypothesis}$$

$$= 0$$

**Lemma 4.31** [121] *Let $\boldsymbol{\Omega} \triangleq (X, +, \cdot, (\mathbb{F}, \dotplus, \dottimes), \boldsymbol{T})$ be a topological linear space. Let* $\operatorname{span} A$ *be the* SPAN *of a set $A$ (Definition 3.2 page 32). Let $\tilde{\mathsf{f}}(\omega)$ and $\tilde{\mathsf{g}}(\omega)$ be the* FOURIER TRANSFORM*s (Definition 2.22 page 20) of the functions $\mathsf{f}(x)$ and $\mathsf{g}(x)$, respectively, in $L_{\mathbb{R}}^2$ (Definition 2.2 page 15). Let $\mathring{\mathsf{a}}(\omega)$ be the DTFT (Definition 2.38 page 24) of a sequence $(a_n)_{n \in \mathbb{Z}}$ in $\ell_{\mathbb{R}}^2$ (Definition 2.33 page 23).*

$$\left\{ \begin{array}{l} \textit{(1).} \quad \left\{ \boldsymbol{T}^n \mathsf{f} \mid_{n \in \mathbb{Z}} \right\} \textit{ is a } \text{SCHAUDER BASIS} \textit{ for } \boldsymbol{\Omega} \quad \text{and} \\ \textit{(2).} \quad \left\{ \boldsymbol{T}^n \mathsf{g} \mid_{n \in \mathbb{Z}} \right\} \textit{ is a } \text{SCHAUDER BASIS} \textit{ for } \boldsymbol{\Omega} \end{array} \right\} \implies \left\{ \begin{array}{l} \exists (a_n)_{n \in \mathbb{Z}} \quad \text{such that} \\ \tilde{\mathsf{f}}(\omega) = \mathring{\mathsf{a}}(\omega) \tilde{\mathsf{g}}(\omega) \end{array} \right\}$$

---

[121] [30], page 140

✎PROOF:   Let $V_0'$ be the space spanned by $\{\mathbf{T}^n\phi|_{n\in\mathbb{Z}}\}$.

$$\tilde{f}(\omega) \triangleq \tilde{\mathbf{F}}f \qquad\qquad\qquad \text{by Definition 2.22 page 20}$$

$$= \tilde{\mathbf{F}} \sum_{n\in\mathbb{Z}} a_n\mathbf{T}g \qquad\qquad \text{by (2)}$$

$$= \sum_{n\in\mathbb{Z}} a_n\tilde{\mathbf{F}}\mathbf{T}g$$

$$= \underbrace{\sum_{n\in\mathbb{Z}} a_n e^{-i\omega n}}_{\breve{a}(\omega)} \underbrace{\tilde{\mathbf{F}}g}_{\tilde{g}(\omega)} \qquad\qquad \text{by Corollary 4.18 page 47}$$

$$= \breve{a}(\omega)\tilde{g}(\omega) \qquad\qquad \text{by Definition 2.38 (page 24) and Definition 2.22 page 20}$$

$$V_0 \triangleq \left\{ f(x) \,\middle|\, f(x) = \sum_{n\in\mathbb{Z}} b_n\mathbf{T}^n g(x) \right\} \qquad \text{by } \textit{Riesz basis} \text{ hypothesis}$$

$$= \left\{ f(x) \,\middle|\, \tilde{\mathbf{F}}f(x) = \tilde{\mathbf{F}} \sum_{n\in\mathbb{Z}} b_n\mathbf{T}^n g(x) \right\}$$

$$= \left\{ f(x) \,\middle|\, \tilde{f}(\omega) = \tilde{b}(\omega)\tilde{g}(\omega) \right\}$$

$$= \left\{ f(x) \,\middle|\, \tilde{f}(\omega) = \tilde{b}(\omega)\breve{a}(\omega)\tilde{f}(\omega) \right\}$$

$$= \left\{ f(x) \,\middle|\, \tilde{f}(\omega) = \tilde{c}(\omega)\tilde{f}(\omega) \right\} \qquad \text{where } \tilde{c}(\omega) \triangleq \tilde{b}(\omega)\breve{a}(\omega)$$

$$= \left\{ f(x) \,\middle|\, f(x) = \sum_{n\in\mathbb{Z}} c_n f(x-n) \right\}$$

$$\triangleq V_0'$$

☞

**Theorem 4.32**  (Battle-Lemarié orthogonalization)  [122] *Let* $\tilde{f}(\omega)$ *be the* FOURIER TRANS-FORM *(Definition 2.22 page 20) of a function* $f \in L^2_{\mathbb{R}}$.

$$\left\{ \begin{array}{l} \textit{1.} \quad \left\{ \mathbf{T}^n g|_{n\in\mathbb{Z}} \right\} \textit{ is a } \text{RIESZ BASIS } \textit{for } L^2_{\mathbb{R}} \quad \textit{and} \\[2ex] \textit{2.} \quad \tilde{f}(\omega) \triangleq \dfrac{\tilde{g}(\omega)}{\sqrt{2\pi \sum\limits_{n\in\mathbb{Z}} |\tilde{g}(\omega+2\pi n)|^2}} \end{array} \right\} \implies \left\{ \begin{array}{l} \left\{ \mathbf{T}^n f|_{n\in\mathbb{Z}} \right\} \\ \textit{is an } \text{ORTHONORMAL BASIS } \textit{for } L^2_{\mathbb{R}} \end{array} \right\}$$

✎PROOF:

---

[122] ☞ [137], page 25, ⟨Remark 2.4⟩, ☞ [132], page 71, ☞ [94], page 72, ☞ [95], page 225, ☞ [30], page 140, ⟨(5.3.3)⟩

(1) Proof that $\left\{\mathbf{T}^n \mathrm{f}\,|_{n\in\mathbb{Z}}\right\}$ is orthonormal:

$$\bar{S}_{\phi\phi}(\omega) = 2\pi \sum_{n\in\mathbb{Z}} \left|\bar{\mathrm{f}}(\omega+2\pi n)\right|^2 \qquad\qquad \text{by Theorem } 4.28 \text{ page } 51$$

$$= 2\pi \sum_{n\in\mathbb{Z}} \left| \frac{\tilde{\mathrm{g}}(\omega+2\pi n)}{\sqrt{2\pi \sum_{m\in\mathbb{Z}} |\tilde{\mathrm{g}}(\omega+2\pi(n-m))|^2}} \right|^2 \qquad \text{by left hypothesis}$$

$$= \sum_{n\in\mathbb{Z}} \left| \frac{\tilde{\mathrm{g}}(\omega+2\pi n)}{\sqrt{\sum_{m\in\mathbb{Z}} |\tilde{\mathrm{g}}(\omega+2\pi m)|^2}} \right|^2$$

$$= \sum_{n\in\mathbb{Z}} \left| \frac{1}{\sqrt{\sum_{m\in\mathbb{Z}} |\tilde{\mathrm{g}}(\omega+2\pi m)|^2}} \right|^2 |\tilde{\mathrm{g}}(\omega+2\pi n)|^2$$

$$= \frac{1}{\sum_{m\in\mathbb{Z}} |\tilde{\mathrm{g}}(\omega+2\pi m)|^2} \sum_{n\in\mathbb{Z}} |\tilde{\mathrm{g}}(\omega+2\pi n)|^2$$

$$= 1$$

$$\implies \left\{\mathrm{f}_n\,|_{n\in\mathbb{Z}}\right\} \text{ is orthonormal} \qquad\qquad \text{by Theorem } 4.30 \text{ page } 56$$

(2) Proof that $\left\{\mathbf{T}^n \mathrm{f}\,|_{n\in\mathbb{Z}}\right\}$ is a basis for $V_0$: by Lemma 4.31 page 57.

(3) Proof that $\left\{\mathbf{D}^m\mathbf{T}^n \mathrm{f}\,|_{n\in\mathbb{Z}}\right\}$ is a basis for $V_m$: by Theorem 5.12 page 71.

## 4.3 Examples

**Example 4.33** (linear functions) [123] Let $\mathbf{T}$ be the *translation operator* (Definition 4.1 page 38). Let $\mathcal{L}(\mathbb{C},\mathbb{C})$ be the set of all *linear* functions in $L^2_{\mathbb{R}}$.

    1. $\{x, \mathbf{T}x\}$ is a *basis* for $\mathcal{L}(\mathbb{C},\mathbb{C})$    and

    2. $\mathrm{f}(x) = \mathrm{f}(1)x - \mathrm{f}(0)\mathbf{T}x$      $\forall \mathrm{f} \in \mathcal{L}(\mathbb{C},\mathbb{C})$

---

[123] ✎ [65], page 2

✎PROOF:   By left hypothesis, f is *linear*; so let $f(x) \triangleq ax + b$

$$
\begin{aligned}
f(1)x - f(0)\mathbf{T}x &= f(1)x - f(0)(x-1) & \text{by Definition 4.1 page 38} \\
&= (ax+b)|_{x=1}\, x - (ax+b)|_{x=0}\,(x-1) & \text{by left hypothesis and definition of f} \\
&= (a+b)x - b(x-1) \\
&= ax + bx - bx + b \\
&= ax + b \\
&= f(x) & \text{by left hypothesis and definition of f}
\end{aligned}
$$

✌

**Example 4.34** (Cardinal Series)   Let **T** be the *translation operator* (Definition 4.1 page 38). The *Paley-Wiener* class of functions $PW_\sigma^2$ are those functions which are "*bandlimited*" with respect to their Fourier transform. The cardinal series forms an orthogonal basis for such a space. The *Fourier coefficients* for a projection of a function f onto the Cardinal series basis elements is particularly simple—these coefficients are samples of $f(x)$ taken at regular intervals. In fact, one could represent the coefficients using inner product notation with the *Dirac delta distribution* $\delta$ as follows:

$$\langle f(x) \mid \mathbf{T}^n \delta(x) \rangle \triangleq \int_{\mathbb{R}} f(x)\delta(x-n)\,\mathrm{d}x = f(n)$$

1. $\left\{ \mathbf{T}^n \dfrac{\sin(\pi x)}{\pi x} \Big|_{n \in \mathbb{N}} \right\}$ is a *basis* for $PW_\sigma^2$   and

2. $f(x) = \underbrace{\displaystyle\sum_{n=1}^{\infty} f(n)\mathbf{T}^n \dfrac{\sin(\pi x)}{\pi x}}_{Cardinal\ series}$          $\forall f \in PW_\sigma^2,\ \sigma \le \frac{1}{2}$

**Example 4.35** (Fourier Series)
1. $\left\{ \mathbf{D}_n e^{ix} |_{n \in \mathbb{Z}} \right\}$ is a *basis* for $L(0,\ 2\pi)$          and
2. $f(x) = \dfrac{1}{\sqrt{2\pi}} \displaystyle\sum_{n \in \mathbb{Z}} \alpha_n \mathbf{D}_n e^{ix}$          $\forall x \in (0, 2\pi), f \in L(0, 2\pi)$   where
3. $\alpha_n \triangleq \dfrac{1}{\sqrt{2\pi}} \displaystyle\int_0^{2\pi} f(x) \mathbf{D}_n e^{-ix}\,\mathrm{d}x$   $\forall f \in L(0, 2\pi)$

**Example 4.36** (Fourier Transform)
1. $\left\{ \mathbf{D}_\omega e^{ix} |_{\omega \in \mathbb{R}} \right\}$ is a *basis* for $L_{\mathbb{R}}^2$          and
2. $f(x) = \dfrac{1}{\sqrt{2\pi}} \displaystyle\int_{\mathbb{R}} \tilde{f}(\omega) \mathbf{D}_x e^{i\omega}\,\mathrm{d}\omega$   $\forall f \in L_{\mathbb{R}}^2$   where
3. $\tilde{f}(\omega) \triangleq \dfrac{1}{\sqrt{2\pi}} \displaystyle\int_{\mathbb{R}} f(x) \mathbf{D}_\omega e^{-ix}\,\mathrm{d}x$   $\forall f \in L_{\mathbb{R}}^2$

**Example 4.37** (Gabor Transform) [124]

1. $\left\{ \left( \mathbf{T}_\tau e^{-\pi x^2} \right) \left( \mathbf{D}_\omega e^{ix} \right) \Big|_{\tau,\omega \in \mathbb{R}} \right\}$ is a *basis* for $\boldsymbol{L}_\mathbb{R}^2$          and

2. $\mathsf{f}(x) = \int_\mathbb{R} \mathsf{G}(\tau,\omega) \, \mathbf{D}_x e^{i\omega} \, \mathrm{d}\omega$          $\forall x \in \mathbb{R}, \, \mathsf{f} \in \boldsymbol{L}_\mathbb{R}^2$   where

3. $\mathsf{G}(\tau,\omega) \triangleq \int_\mathbb{R} \mathsf{f}(x) \left( \mathbf{T}_\tau e^{-\pi x^2} \right) \left( \mathbf{D}_\omega e^{-ix} \right) \, \mathrm{d}x$   $\forall x \in \mathbb{R}, \, \mathsf{f} \in \boldsymbol{L}_\mathbb{R}^2$

**Example 4.38** (wavelets)  Let $\psi(x)$ be a *mother wavelet*.

1. $\left\{ \mathbf{D}^k \mathbf{T}^n \psi(x) \big|_{k,n \in \mathbb{Z}} \right\}$ is a *basis* for $\boldsymbol{L}_\mathbb{R}^2$          and

2. $\mathsf{f}(x) = \sum_{k \in \mathbb{Z}} \sum_{n \in \mathbb{Z}} \alpha_{k,n} \mathbf{D}^k \mathbf{T}^n \psi(x)$   $\forall \mathsf{f} \in \boldsymbol{L}_\mathbb{R}^2$   where

3. $\alpha_n \triangleq \int_\mathbb{R} \mathsf{f}(x) \mathbf{D}^k \mathbf{T}^n \psi^*(x) \, \mathrm{d}x$   $\forall \mathsf{f} \in \boldsymbol{L}_\mathbb{R}^2$

**Definition 4.39** [125]

1. $\mathbf{D}_\alpha \mathsf{f}(x) \triangleq \sqrt{\alpha} \mathsf{f}(\alpha x)$          $\forall \mathsf{f} \in \mathbb{R}^\mathbb{R}$          (**generalized dilation operator**)

2. $\mathbf{E}_\alpha \mathsf{f}(x) \triangleq e^{i2\pi\alpha x} \mathsf{f}(x)$          $\forall \mathsf{f} \in \mathbb{R}^\mathbb{R}$          (**modulation operator**)

**Example 4.40** (Poisson Summation Formula)  Let $\mathbf{E}$ be the *modulation operator* and $\mathbf{T}$ the *translation operator* (Definition 4.1 page 38). Let $\tilde{\mathsf{f}}(\omega)$ be the Fourier transform of a function $\mathsf{f}(x)$.

$$\underbrace{\sum_{n \in \mathbb{Z}} \mathbf{E}_n \tilde{\mathsf{f}}(2\pi n)}_{\text{modulated summation in "frequency"}} = \underbrace{\frac{1}{\sqrt{2\pi}} \sum_{n \in \mathbb{Z}} \mathbf{T}^n \mathsf{f}(x)}_{\text{summation in "time"}} \qquad \forall \mathsf{f} \in \boldsymbol{L}_\mathbb{R}^2$$

✎ PROOF:   See Theorem 4.21 page 48.   ✏

**Example 4.41** (B-splines)  Let $\mathbf{M}$ be the *multiplication operator* and $\mathbf{T}$ the *translation operator* (Definition 4.1 page 38). Let $\mathsf{N}_n(x)$ be the *$n$th order cardinal B-spline*.

$$\mathsf{N}_n(x) = \frac{1}{n} x \mathsf{N}_{n-1}(x) - \frac{1}{n} x \mathbf{T} \mathsf{N}_{n-1}(x) + \frac{n+1}{n} \mathbf{T} \mathsf{N}_{n-1}(x) \qquad \forall n \in \mathbb{N} \setminus \{1\}, \forall x \in \mathbb{R}$$

**Example 4.42** (Fourier Series analysis)  Let $\mathbf{D}_\alpha$ be the *dilation operator* and $\mathbf{E}$ the *modulation operator* (Definition 4.39 page 61). Let $\hat{\mathbf{F}}$ be the Fourier Series operator.

The *inverse Fourier Series* operator $\hat{\mathbf{F}}^{-1}$ is given by

$$\mathbf{f}(n) = \sum_{n\in\mathbb{Z}} \tilde{\mathbf{f}}(n)\frac{1}{\sqrt{n}}\mathbf{D}_n e^{-i2\pi t} = \sum_{n\in\mathbb{Z}} \mathbf{E}_n \tilde{\mathbf{f}}(n) \qquad \forall \mathbf{f}\in L_{\mathbb{R}}^2$$

where

$$\tilde{\mathbf{f}}(\omega) \triangleq \int_0^1 \mathbf{E}_{-n}\mathbf{f}(x)\,\mathrm{d}x$$

✎Proof:    See Theorem 2.17 page 18.

# 5    Demonstrations/Applications

This section demontrates some applications of the transversal operators in some well known analytic systems. In particular, the usefulness of the operators in the proofs of some traditional results is demonstrated.

## 5.1    Multiresolution analysis

### 5.1.1    Definition

A multiresolution analysis provides "coarse" approximations of a function in a linear space $L_{\mathbb{R}}^2$ at multiple "scales" or "resolutions". Key to this process is a sequence of *scaling functions*. Most traditional transforms feature a single *scaling function* $\phi(x)$ set equal to one ($\phi(x) = 1$). This allows for convenient representation of the most basic functions, such as constants.[126] A multiresolution system, on the other hand, uses a generalized form of the scaling concept:[127]

(1) Instead of the scaling function simply being set *equal to unity* ($\phi(x) = 1$), a multiresolution analysis (Definition 5.1 page 63) is often constructed in such a way that the scaling function $\phi(x)$ forms a *partition of unity* such that $\sum_{n\in\mathbb{Z}} \mathbf{T}^n \phi(x) = 1$.

(2) Instead of there being *just one* scaling function, there is an entire sequence of scaling functions $\big(\mathbf{D}^j \phi(x)\big)_{j\in\mathbb{Z}}$, each corresponding to a different "*resolution*".

---

[126]  [73], page 8

[127] The concept of a scaling space was perhaps first introduced by Taizo Iijima in 1959 in Japan, and later as the *Gaussian Pyramid* by Burt and Adelson in the 1980s in the West. ✎ [94], page 70, ✎ [69], ✎ [20], ✎ [2], ✎ [91], ✎ [4], ✎ [54], ✎ [134], ⟨historical survey⟩

**Definition 5.1** [128] Let $\left(\!\left( V_j \right)\!\right)_{j\in\mathbb{Z}}$ be a sequence of subspaces on $L_\mathbb{R}^2$. Let $A^-$ be the *closure* of a set $A$. The sequence $\left(\!\left( V_j \right)\!\right)_{j\in\mathbb{Z}}$ is a **multiresolution analysis** on $L_\mathbb{R}^2$ if

1. $V_j = V_j^-$                $\forall j\in\mathbb{Z}$        (*closed*)              and
2. $V_j \subset V_{j+1}$              $\forall j\in\mathbb{Z}$        (*linearly ordered*)        and
3. $\left( \bigcup_{j\in\mathbb{Z}} V_j \right)^- = L_\mathbb{R}^2$                    (*dense* in $L_\mathbb{R}^2$)        and
4. $f \in V_j \iff \mathbf{D}f \in V_{j+1}$    $\forall j\in\mathbb{Z}, f\in L_\mathbb{R}^2$    (*self-similar*)        and
5. $\exists \phi$    such that    $\left\{ \mathbf{T}^n\phi|_{n\in\mathbb{Z}} \right\}$ is a *Riesz basis* for $V_0$.

A *multiresolution analysis* is also called an **MRA**. An element $V_j$ of $\left(\!\left( V_j \right)\!\right)_{j\in\mathbb{Z}}$ is a **scaling subspace** of the space $L_\mathbb{R}^2$. The pair $\left( L_\mathbb{R}^2, \left(\!\left( V_j \right)\!\right) \right)$ is a **multiresolution analysis space**, or **MRA space**. The function $\phi$ is the **scaling function** of the *MRA space*.

The traditional definition of the *MRA* also includes the following:

6. $f \in V_j \iff \mathbf{T}^n f \in V_j$    $\forall n,j\in\mathbb{Z}, f\in L_\mathbb{R}^2$    (*translation invariant*)
7. $\bigcap_{j\in\mathbb{Z}} V_j = \{0\}$                    (*greatest lower bound* is $\mathbf{0}$)

However, Proposition 5.2 (next) and Proposition 5.3 (page 64) demonstrate that these follow from the *MRA* as defined in Definition 5.1.

**Proposition 5.2** [129] *Let* MRA *be defined as in Definition 5.1 page 63*.

$$\left\{ \left(\!\left( V_j \right)\!\right)_{j\in\mathbb{Z}} \text{ is an MRA} \right\} \implies \underbrace{\left\{ f \in V_j \iff \mathbf{T}^n f \in V_j \quad \forall n,j\in\mathbb{Z}, f\in L_\mathbb{R}^2 \right\}}_{\text{TRANSLATION INVARIANT}}$$

✎ PROOF:

$\mathbf{T}^n f \in V_j$

$\iff \mathbf{T}^n f \in \operatorname{span}\left\{ \mathbf{D}^j\mathbf{T}^m\phi |_{m\in\mathbb{Z}} \right\}$                by definition of $\{\phi\}$ (Definition 5.1 page 63)

$\iff \exists\, (\alpha_n)_{n\in\mathbb{Z}}$    such that    $\mathbf{T}^n f(x) = \sum_{k\in\mathbb{Z}} \alpha_k \mathbf{D}^j\mathbf{T}^k\phi(x)$        by definition of $\{\phi\}$ (Definition 5.1 page 63)

$\iff \exists\, (\alpha_n)_{n\in\mathbb{Z}}$    such that    $f(x) = \mathbf{T}^{-n}\sum_{k\in\mathbb{Z}} \alpha_k \mathbf{D}^j\mathbf{T}^k\phi(x)$    by definition of $\mathbf{T}$ (Definition 4.1 page 38)

$= \sum_{k\in\mathbb{Z}} \alpha_k \mathbf{T}^{-n}\mathbf{D}^j\mathbf{T}^k\phi(x)$

$= \sum_{k\in\mathbb{Z}} \alpha_k \mathbf{D}^j\mathbf{T}^{k-2n}\phi(x)$        by Proposition 4.8 page 41

---

$$= \sum_{\ell \in \mathbb{Z}} \alpha_{\ell + 2n} \mathbf{D}^j \mathbf{T}^\ell \phi(x) \qquad \text{where } \ell \triangleq k - 2n \implies k = \ell + 2n$$

$$= \sum_{\ell \in \mathbb{Z}} \beta_\ell \mathbf{D}^j \mathbf{T}^\ell \phi(x) \qquad \text{where } \beta_\ell \triangleq \alpha_{\ell + 2n}$$

$$\iff \quad \mathsf{f} \in V_j \qquad \text{by definition of } \{\mathbf{T}^n \phi\} \text{ (Definition 5.1 page 63)}$$

**Proposition 5.3** [130] *Let* MRA *be defined as in Definition 5.1 page 63.*

$$\Big\{ \, (\!(V_j)\!)_{j \in \mathbb{Z}} \text{ is an MRA} \, \Big\} \qquad \implies \qquad \Big\{ \, \bigcap_{j \in \mathbb{Z}} V_j = \{0\} \qquad \text{(greatest lower bound is 0)} \, \Big\}$$

**Proof:**

(1) Let $\mathbf{P}_j$ be the *projection operator* that generates the scaling subspace $V_j$ such that
$$V_j = \{\mathbf{P}_j \mathsf{f} \mid \mathsf{f} \in L^2_{\mathbb{R}}\}$$

(2) lemma: Functions with *compact support* are *dense* in $L^2_{\mathbb{R}}$. Therefore, we only need to prove that the proposition is true for functions with support in $[-R, R]$, for all $R > 0$.

(3) For some function $\mathsf{f} \in L^2_{\mathbb{R}}$, let $(\!(\mathsf{f}_n)\!)_{n \in \mathbb{Z}}$ be a sequence of functions in $L^2_{\mathbb{R}}$ with *compact support* such that
$$\operatorname{supp} \mathsf{f}_n \subseteq [-R, R] \text{ for some } R > 0 \quad \text{and} \quad \mathsf{f}(x) = \lim_{n \to \infty} (\!(\mathsf{f}_n(x))\!).$$

(4) lemma: $\bigcap V_j = \{0\} \quad \iff \quad \lim_{j \to -\infty} \|\mathbf{P}_j \mathsf{f}\| = 0 \quad \forall \mathsf{f} \in L^2_{\mathbb{R}}$. Proof:

$$\bigcap_{j \in \mathbb{Z}} V_j = \bigcap_{j \in \mathbb{Z}} \{\mathbf{P}_j \mathsf{f} \mid \mathsf{f} \in L^2_{\mathbb{R}}\} \qquad \text{by definition of } V_j \text{ (definition 1 page 64)}$$

$$= \lim_{j \to -\infty} \{\mathbf{P}_j \mathsf{f} \mid \mathsf{f} \in L^2_{\mathbb{R}}\} \qquad \text{by definition of } \cap$$

$$= 0 \iff \lim_{j \to -\infty} \|\mathbf{P}_j \mathsf{f}\| = 0 \qquad \text{by } \textit{nondegenerate} \text{ property of } \|\cdot\| \text{ (Definition 3.6 page 33)}$$

(5) lemma: $\lim_{j \to -\infty} \|\mathbf{P}_j \mathsf{f}\| = 0 \quad \forall \mathsf{f} \in L^2_{\mathbb{R}}$. Proof:

---

[130] ✎ [137], pages 19–28, ⟨Proposition 2.14⟩, ✎ [64], page 45, ⟨Theorem 1.6⟩, ✎ [106], pages 313–314, ⟨Lemma 6.4.28⟩

Let $\mathbb{1}_A(x)$ be the *set indicator function* (Definition 2.3 page 15)

$$\lim_{j \to -\infty} \left\| \mathbf{P}_j f \right\|^2$$

$$= \lim_{j \to -\infty} \left\| \mathbf{P}_j \lim_{n \to \infty} (\!(f_n)\!) \right\|^2 \qquad \text{by definition 3 page 64}$$

$$\leq \lim_{j \to -\infty} B \sum_{n \in \mathbb{Z}} \left| \left\langle \mathbf{P}_j \lim_{n \to \infty} (\!(f_n)\!) \,\middle|\, \mathbf{D}^j \mathbf{T}^n \phi \right\rangle \right|^2 \qquad \text{by \emph{frame} prop. (Proposition 3.24 page 38)}$$

$$= \lim_{j \to -\infty} B \sum_{n \in \mathbb{Z}} \left| \left\langle \lim_{n \to \infty} (\!(f_n)\!) \,\middle|\, \mathbf{D}^j \mathbf{T}^n \phi \right\rangle \right|^2 \qquad \text{by definition of } \mathbf{P}_j \text{ (definition 1 page 64)}$$

$$= \lim_{j \to -\infty} B \sum_{n \in \mathbb{Z}} \left| \left\langle \mathbb{1}_{[-R,\,R]}(x) \lim_{n \to \infty} (\!(f_n)\!) \,\middle|\, \mathbf{D}^j \mathbf{T}^n \phi(x) \right\rangle \right|^2 \qquad \text{by def. of } (\!(f_n)\!) \text{ (definition 3 page 64)}$$

$$= \lim_{j \to -\infty} B \sum_{n \in \mathbb{Z}} \left| \left\langle \lim_{n \to \infty} (\!(f_n)\!) \,\middle|\, \mathbb{1}_{[-R,\,R]}(x) \mathbf{D}^j \mathbf{T}^n \phi(x) \right\rangle \right|^2 \qquad \text{prop. of } \langle \triangle \mid \triangledown \rangle \text{ in } L^2_{\mathbb{R}}$$

$$\leq \lim_{j \to -\infty} B \sum_{n \in \mathbb{Z}} \left\| \lim_{n \to \infty} (\!(f_n)\!) \right\|^2 \left\| \mathbb{1}_{[-R,\,R]}(x) \mathbf{D}^j \mathbf{T}^n \phi(x) \right\|^2 \qquad \text{by \emph{Cauchy-Schwarz Inequality}}$$

$$= \lim_{j \to -\infty} B \sum_{n \in \mathbb{Z}} \|f\|^2 \left\| \mathbb{1}_{[-R,\,R]}(x) \mathbf{D}^j \mathbf{T}^n \phi(x) \right\|^2 \qquad \text{by def. of } (\!(f_n)\!) \text{ (definition 3 page 64)}$$

$$= \lim_{j \to -\infty} B \sum_{n \in \mathbb{Z}} \|f\|^2 \left\| \left[ \underbrace{\mathbf{D}^j \mathbf{D}^{-j}}_{\mathbf{I}} \mathbb{1}_{[-R,\,R]}(x) \right] \left[ \mathbf{D}^j \mathbf{T}^n \phi(x) \right] \right\|^2 \qquad \text{by Proposition 4.3 page 39}$$

$$= \lim_{j \to -\infty} B \sum_{n \in \mathbb{Z}} \|f\|^2 \, 2^{j/2} \mathbf{D}^j \left\{ \left[ \mathbf{D}^{-j} \mathbb{1}_{[-R,\,R]}(x) \right] \left[ \mathbf{T}^n \phi(x) \right] \right\}^2 \qquad \text{by Proposition 4.7 page 41}$$

$$= \lim_{j \to -\infty} B \sum_{n \in \mathbb{Z}} \|f\|^2 \left\| \mathbf{D}^j \left\{ 2^{j/2} 2^{-j/2} \mathbb{1}_{[-R,\,R]}(2^{-j}x) \left[ \mathbf{T}^n \phi(x) \right] \right\} \right\|^2 \qquad \text{by Proposition 4.3 page 39}$$

$$= \lim_{j \to -\infty} B \sum_{n \in \mathbb{Z}} \|f\|^2 \left\| \mathbf{D}^j \left\{ \left[ \underbrace{\mathbf{T}^n \mathbf{T}^{-n}}_{\mathbf{I}} \mathbb{1}_{[-R,\,R]}(2^{-j}x) \right] \left[ \mathbf{T}^n \phi(x) \right] \right\} \right\|^2 \qquad \text{by Proposition 4.3 page 39}$$

$$= \lim_{j \to -\infty} B \sum_{n \in \mathbb{Z}} \|f\|^2 \left\| \mathbf{D}^j \left\{ \left[ \mathbf{T}^n \mathbb{1}_{[-R,\,R]}(2^{-j}x + n) \right] \left[ \mathbf{T}^n \phi(x) \right] \right\} \right\|^2 \qquad \text{by Proposition 4.3 page 39}$$

$$= \lim_{j \to -\infty} B \sum_{n \in \mathbb{Z}} \|f\|^2 \left\| \mathbf{D}^j \mathbf{T}^n \left\{ \mathbb{1}_{[-R,\,R]}(2^{-j}x + n) \phi(x) \right\} \right\|^2 \qquad \text{by Proposition 4.3 page 39}$$

$$= \lim_{j \to -\infty} B \sum_{n \in \mathbb{Z}} \|f\|^2 \left\| \mathbb{1}_{[-R,\,R]}(2^{-j}x + n) \phi(x) \right\|^2 \qquad \text{by Theorem 4.15 page 45}$$

$$= B \|f\|^2 \sum_{n \in \mathbb{Z}} \lim_{j \to -\infty} \left\| \mathbb{1}_{[-2^j R + n,\, 2^j R + n]}(u) \phi(2^{-j}(u - n)) \right\|^2 \qquad u \triangleq 2^j x + n \implies x = 2^{-j}(u - n)$$

$$= B \|f\|^2 \sum_{n \in \mathbb{Z}} \lim_{j \to -\infty} \int_{-2^j R + n}^{2^j R + n} \left| \phi(2^{-j}(u - n)) \right|^2 \, \mathrm{d}u$$

$$= B \, \|f\|^2 \sum_{n \in \mathbb{Z}} \int_n^n |\phi(0)|^2 \, \mathrm{d}u$$

$$= 0$$

(6)  Final step in proof that $\bigcap V_j = \{0\}$: by lemma **4** page **64** and lemma **5** page **65**

☞

The MRA (Definition **5.1** page **63**) is more than just an interesting mathematical toy. Under some very "reasonable" conditions (next proposition), as $j \to \infty$, the *scaling subspace* $V_j$ is *dense* in $L^2_\mathbb{R}$ …meaning that with the MRA we can represent any "reasonable" function to within an arbitrary accuracy.

**Proposition 5.4** [131]

$$\left\{ \begin{array}{ll} (1). & (\mathbf{T}^n\phi) \; is \; a \; \text{Riesz sequence} & and \\ (2). & \tilde{\phi}(\omega) \; is \; \text{continuous} \; at \; 0 & and \\ (3). & \tilde{\phi}(0) \neq 0 \end{array} \right\} \implies \left\{ \left( \bigcup_{j \in \mathbb{Z}} V_j \right)^{-} = L^2_\mathbb{R} \quad (\text{dense in } L^2_\mathbb{R}) \right\}$$

### 5.1.2  Dilation equation

Several functions in mathematics exhibit a kind of *self-similar* or *recursive* property:

☞ If a function $f(x)$ is *linear*, then (Example **4.33** page **59**)
$$f(x) = f(1)x - f(0)\mathbf{T}x.$$

☞ If a function $f(x)$ is sufficiently *bandlimited*, then the *Cardinal series* (Example **4.34** page **60**) demonstrates
$$f(x) = \sum_{n=1}^{\infty} f(n)\mathbf{T}^n \frac{\sin[\pi(x)]}{\pi(x)}.$$

☞ *B-splines* are another example:
$$\mathsf{N}_n(x) = \frac{1}{n} x \mathsf{N}_{n-1}(x) - \frac{1}{n} x \mathbf{T} \mathsf{N}_{n-1}(x) + \frac{n+1}{n} \mathbf{T} \mathsf{N}_{n-1}(x) \qquad \forall n \in \mathbb{N}\setminus\{1\}, \forall x \in \mathbb{R}.$$

The scaling function $\phi(x)$ (Definition **5.1** page **63**) also exhibits a kind of *self-similar* property. By Definition **5.1** page **63**, the dilation $\mathbf{D}f$ of each vector $f$ in $V_0$ is in $V_1$. If $\{\mathbf{T}^n\phi|_{n\in\mathbb{Z}}\}$ is a basis for $V_0$, then $\{\mathbf{DT}^n\phi|_{n\in\mathbb{Z}}\}$ is a basis for $V_1$, $\{\mathbf{D}^2\mathbf{T}^n\phi|_{n\in\mathbb{Z}}\}$ is a basis for $V_2$, …; and in general $\{\mathbf{D}^j\mathbf{T}^m\phi \mid j \in \mathbb{Z}\}$ is a basis for $V_j$. Also, if $\phi$ is in $V_0$, then it is also in $V_1$ (because $V_0 \subset V_1$). And because $\phi$ is in $V_1$ and because $\{\mathbf{DT}^n\phi|_{n\in\mathbb{Z}}\}$ is a basis for $V_1$, $\phi$ is a linear combination of the elements in $\{\mathbf{DT}^n\phi|_{n\in\mathbb{Z}}\}$. That is, $\phi$ can be represented as a linear

---

[131] ✎ [137], pages 28–31, ⟨Proposition 2.15⟩, ✎ [53], pages 35–37, ⟨Proposition 2.3⟩

combination of translated and dilated versions of itself. The resulting equation is called the *dilation equation* (Definition 5.5, next).[132]

**Definition 5.5** [133] Let $\left( L^2_{\mathbb{R}}, \left( V_j \right) \right)$ be a *multiresolution analysis space* with scaling function $\phi$ (Definition 5.1 page 63). Let $\left( h_n \right)_{n \in \mathbb{Z}}$ be a *sequence* (Definition 2.32 page 23) in $\ell^2_{\mathbb{R}}$ (Definition 2.33 page 23). The equation

$$\phi(x) = \sum_{n \in \mathbb{Z}} h_n \mathbf{DT}^n \phi(x) \qquad \forall x \in \mathbb{R}$$

is called the **dilation equation**. It is also called the **refinement equation**, **two-scale difference equation**, and **two-scale relation**.

**Theorem 5.6** (dilation equation) *Let an* MRA SPACE *and* SCALING FUNCTION *be as defined in Definition 5.1 page 63.*

$$\left\{ \begin{array}{l} \left( L^2_{\mathbb{R}}, \left( V_j \right) \right) \text{ is an } \text{MRA SPACE} \\ \text{with } \text{SCALING FUNCTION } \phi \end{array} \right\} \implies \underbrace{\left\{ \begin{array}{l} \exists \left( h_n \right)_{n \in \mathbb{Z}} \text{ such that} \\ \phi(x) = \sum_{n \in \mathbb{Z}} h_n \mathbf{DT}^n \phi(x) \qquad \forall x \in \mathbb{R} \end{array} \right\}}_{\text{DILATION EQUATION IN "TIME"}}$$

✎PROOF:

$$\begin{aligned} \phi &\in V_0 & \text{by Definition 5.1 page 63} \\ &\subseteq V_1 & \text{by Definition 5.1 page 63} \\ &= \operatorname{span}\left\{ \mathbf{DT}^n \phi(x) |_{n \in \mathbb{Z}} \right\} \\ &\implies \exists \left( h_n \right)_{n \in \mathbb{Z}} \quad \text{such that} \quad \phi = \sum_{n \in \mathbb{Z}} h_n \mathbf{DT}^n \phi \end{aligned}$$

✏

**Lemma 5.7** [134] *Let* $\phi(x)$ *be a function in* $L^2_{\mathbb{R}}$ (Definition 2.2 page 15). *Let* $\tilde{\phi}(\omega)$ *be the* FOURIER TRANSFORM *of* $\phi(x)$. *Let* $\breve{h}(\omega)$ *be the* DISCRETE TIME FOURIER TRANSFORM *of a sequence* $\left( h_n \right)_{n \in \mathbb{Z}}$.

$$(A) \quad \phi(x) = \sum_{n \in \mathbb{Z}} h_n \mathbf{DT}^n \phi(x) \quad \forall x \in \mathbb{R} \quad \iff \quad \tilde{\phi}(\omega) = \frac{\sqrt{2}}{2} \breve{h}\left( \frac{\omega}{2} \right) \tilde{\phi}\left( \frac{\omega}{2} \right) \qquad \forall \omega \in \mathbb{R} \qquad (1)$$

$$\iff \quad \tilde{\phi}(\omega) = \tilde{\phi}\left( \frac{\omega}{2^N} \right) \prod_{n=1}^{N} \frac{\sqrt{2}}{2} \breve{h}\left( \frac{\omega}{2^n} \right) \quad \forall n \in \mathbb{N}, \omega \in \mathbb{R} \qquad (2)$$

✎PROOF:

---

[132]The property of *translation invariance* is of particular significance in the theory of *normed linear spaces* (a Hilbert space is a complete normed linear space equipped with an inner product).

[133] ☞ [73], page 7

[134] ☞ [95], page 228

(1) Proof that (A) $\implies$ (1):

$$\tilde{\phi}(\omega) \triangleq \tilde{\mathbf{F}}\phi$$

$$= \tilde{\mathbf{F}} \sum_{n \in \mathbb{Z}} h_n \mathbf{D}\mathbf{T}^n \phi(x) \qquad \text{by (A)}$$

$$= \sum_{n \in \mathbb{Z}} h_n \tilde{\mathbf{F}}\mathbf{D}\mathbf{T}^n \phi(x)$$

$$= \sum_{n \in \mathbb{Z}} h_n \underbrace{\frac{\sqrt{2}}{2} e^{-i\frac{\omega}{2}n} \phi\left(\frac{\omega}{2}\right)}_{\tilde{\mathbf{F}}\mathbf{D}\mathbf{T}^n \phi(x)} \qquad \text{by Proposition 4.19 page 47}$$

$$= \frac{\sqrt{2}}{2} \underbrace{\left[\sum_{n \in \mathbb{Z}} h_n e^{-i\frac{\omega}{2}n}\right]}_{\tilde{h}(\omega/2)} \tilde{\phi}\left(\frac{\omega}{2}\right)$$

$$= \frac{\sqrt{2}}{2} \tilde{h}\left(\frac{\omega}{2}\right) \tilde{\phi}\left(\frac{\omega}{2}\right) \qquad \text{by definition of } DTFT \text{ (Definition 2.38 page 24)}$$

(2) Proof that (A) $\impliedby$ (1):

$$\phi(x) = \tilde{\mathbf{F}}^{-1} \tilde{\phi}(\omega) \qquad \text{by definition of } \tilde{\phi}(\omega)$$

$$= \tilde{\mathbf{F}}^{-1} \frac{\sqrt{2}}{2} \tilde{h}\left(\frac{\omega}{2}\right) \tilde{\phi}\left(\frac{\omega}{2}\right) \qquad \text{by (1)}$$

$$= \tilde{\mathbf{F}}^{-1} \frac{\sqrt{2}}{2} \sum_{n \in \mathbb{Z}} h_n e^{-i\frac{\omega}{2}n} \tilde{\phi}\left(\frac{\omega}{2}\right) \qquad \text{by definition of } DTFT \text{ (Definition 2.38 page 24)}$$

$$= \frac{\sqrt{2}}{2} \sum_{n \in \mathbb{Z}} h_n \tilde{\mathbf{F}}^{-1} e^{-i\frac{\omega}{2}n} \tilde{\phi}\left(\frac{\omega}{2}\right) \qquad \text{by property of linear operators}$$

$$= \frac{\sqrt{2}}{2} \sum_{n \in \mathbb{Z}} h_n \tilde{\mathbf{F}}^{-1} \tilde{\mathbf{F}}\mathbf{D}\mathbf{T}^n \phi \qquad \text{by Proposition 4.19 page 47}$$

$$= \sum_{n \in \mathbb{Z}} h_n \mathbf{D}\mathbf{T}^n \phi(x)$$

(3) Proof that (1) $\implies$ (2):

  (a) Proof for $N = 1$ case:

$$\tilde{\phi}\left(\frac{\omega}{2^N}\right) \prod_{n=1}^{N} \frac{\sqrt{2}}{2} \tilde{h}\left(\frac{\omega}{2^n}\right)\bigg|_{N=1} = \frac{\sqrt{2}}{2} \tilde{h}\left(\frac{\omega}{2}\right) \tilde{\phi}\left(\frac{\omega}{2}\right)$$

$$= \tilde{\phi}(\omega) \qquad \text{by (1)}$$

  (b) Proof that [$N$ case] $\implies$ [$N + 1$ case]:

$$\tilde{\phi}\left(\frac{\omega}{2^{N+1}}\right) \prod_{n=1}^{N+1} \frac{\sqrt{2}}{2} \tilde{h}\left(\frac{\omega}{2^n}\right) = \left[\prod_{n=1}^{N} \frac{\sqrt{2}}{2} \tilde{h}\left(\frac{\omega}{2^n}\right)\right] \underbrace{\frac{\sqrt{2}}{2} \tilde{h}\left(\frac{\omega}{2^{N+1}}\right) \tilde{\phi}\left(\frac{\omega}{2^{N+1}}\right)}_{\tilde{\phi}(\omega/2^N)}$$

$$= \tilde{\phi}(\omega/2^N) \prod_{n=1}^{N} \frac{\sqrt{2}}{2} \breve{h}\left(\frac{\omega}{2^n}\right)$$

$$= \tilde{\phi}(\omega) \hspace{3cm} \text{by } [N \text{ case}] \text{ hypothesis}$$

(4) Proof that (1) $\Longleftarrow$ (2):

$$\tilde{\phi}(\omega) = \tilde{\phi}\left(\frac{\omega}{2^N}\right) \prod_{n=1}^{N} \frac{\sqrt{2}}{2} \breve{h}\left(\frac{\omega}{2^n}\right)\Bigg|_{N=1} \hspace{2cm} \text{by (2)}$$

$$= \tilde{\phi}\left(\frac{\omega}{2}\right) \frac{\sqrt{2}}{2} \breve{h}\left(\frac{\omega}{2}\right)$$

$$= \frac{\sqrt{2}}{2} \breve{h}\left(\frac{\omega}{2}\right) \tilde{\phi}\left(\frac{\omega}{2}\right)$$

**Lemma 5.8**   *Let $\phi(x)$ be a function in $\mathbf{L}_{\mathbb{R}}^2$* *(Definition 2.2 page 15)*. *Let $\tilde{\phi}(\omega)$ be the* Fourier trans-form *of $\phi(x)$. Let $\breve{h}(\omega)$ be the* Discrete time Fourier transform *of $(h_n)$. Let $\prod_{n=1}^{\infty} x_n \triangleq$*

$\lim_{N \to \infty} \prod_{n=1}^{N} x_n$, *with respect to the standard norm in $\mathbf{L}_{\mathbb{R}}^2$.*

$$\left\{ \begin{array}{l} \tilde{\phi}(\omega) = C \prod_{n=1}^{\infty} \frac{\sqrt{2}}{2} \breve{h}\left(\frac{\omega}{2^n}\right) \\ \forall C>0,\, \omega \in \mathbb{R} \hspace{1.5cm} (A) \end{array} \right\} \implies \begin{array}{ll} \phi(x) = \displaystyle\sum_{n \in \mathbb{Z}} h_n \mathbf{D} \mathbf{T}^n \phi(x) & \forall x \in \mathbb{R} \hspace{1cm} (1) \\[2ex] \Longleftrightarrow \quad \tilde{\phi}(\omega) = \frac{\sqrt{2}}{2} \breve{h}\left(\frac{\omega}{2}\right) \tilde{\phi}\left(\frac{\omega}{2}\right) & \forall \omega \in \mathbb{R} \hspace{1cm} (2) \\[2ex] \Longleftrightarrow \quad \tilde{\phi}(\omega) = \tilde{\phi}\left(\frac{\omega}{2^N}\right) \displaystyle\prod_{n=1}^{N} \frac{\sqrt{2}}{2} \breve{h}\left(\frac{\omega}{2^n}\right) & \forall n \in \mathbb{N},\, \omega \in \mathbb{R} \hspace{0.3cm} (3) \end{array}$$

✎ Proof:

(1) Proof that (1) $\Longleftrightarrow$ (2) $\Longleftrightarrow$ (3): by Lemma 5.7 page 67

(2) Proof that (A) $\implies$ (2):

$$\tilde{\phi}(\omega) = C \prod_{n=1}^{\infty} \tfrac{\sqrt{2}}{2} \breve{\mathrm{h}}\left(\frac{\omega}{2^n}\right) \qquad \text{by left hypothesis}$$

$$= C \; \tfrac{\sqrt{2}}{2} \breve{\mathrm{h}}\left(\frac{\omega}{2}\right) \prod_{n=1}^{\infty} \tfrac{\sqrt{2}}{2} \breve{\mathrm{h}}\left(\frac{\omega}{2^{n+1}}\right)$$

$$= C \; \tfrac{\sqrt{2}}{2} \breve{\mathrm{h}}\left(\frac{\omega}{2}\right) \prod_{n=1}^{\infty} \tfrac{\sqrt{2}}{2} \breve{\mathrm{h}}\left(\frac{\omega/2}{2^n}\right)$$

$$= \tfrac{\sqrt{2}}{2} \breve{\mathrm{h}}\left(\frac{\omega}{2}\right) \left[ C \prod_{n=1}^{\infty} \tfrac{\sqrt{2}}{2} \breve{\mathrm{h}}\left(\frac{\omega/2}{2^n}\right) \right]$$

$$= \tfrac{\sqrt{2}}{2} \breve{\mathrm{h}}\left(\frac{\omega}{2}\right) \tilde{\phi}\left(\frac{\omega}{2}\right) \qquad \text{by left hypothesis}$$

$\Longleftrightarrow$

**Proposition 5.9**  *Let $\phi(x)$ be a function in $\boldsymbol{L}_{\mathbb{R}}^2$ (Definition 2.2 page 15). Let $\tilde{\phi}(\omega)$ be the* FOURIER

TRANSFORM *of $\phi(x)$. Let $\breve{\mathrm{h}}(\omega)$ be the* DISCRETE TIME FOURIER TRANSFORM *of $(h_n)$. Let $\displaystyle\prod_{n=1}^{\infty} x_n \triangleq$*

$\displaystyle\lim_{N \to \infty} \prod_{n=1}^{N} x_n$, *with respect to the standard norm in $\boldsymbol{L}_{\mathbb{R}}^2$.*

$$\left\{ \begin{array}{l} \tilde{\phi}(\omega) \text{ is} \\ \text{CONTINUOUS} \\ \text{at } \omega = 0 \end{array} \right\} \implies \left\{ \begin{array}{rcll} \phi(x) & = & \displaystyle\sum_{n \in \mathbb{Z}} h_n \mathbf{DT}^n \phi(x) & \forall x \in \mathbb{R} \quad (1) \\[2ex] \Longleftrightarrow \; \tilde{\phi}(\omega) & = & \tfrac{\sqrt{2}}{2} \breve{\mathrm{h}}\left(\frac{\omega}{2}\right) \tilde{\phi}\left(\frac{\omega}{2}\right) & \forall \omega \in \mathbb{R} \quad (2) \\[2ex] \Longleftrightarrow \; \tilde{\phi}(\omega) & = & \tilde{\phi}\left(\frac{\omega}{2^N}\right) \displaystyle\prod_{n=1}^{N} \tfrac{\sqrt{2}}{2} \breve{\mathrm{h}}\left(\frac{\omega}{2^n}\right) & \forall n \in \mathbb{N}, \omega \in \mathbb{R} \quad (3) \\[2ex] \Longleftrightarrow \; \tilde{\phi}(\omega) & = & \tilde{\phi}(0) \displaystyle\prod_{n=1}^{\infty} \tfrac{\sqrt{2}}{2} \breve{\mathrm{h}}\left(\frac{\omega}{2^n}\right) & \omega \in \mathbb{R} \quad (4) \end{array} \right\}$$

✎ PROOF:

(1) Proof that (1) $\iff$ (2) $\iff$ (3): by Lemma 5.7 page 67

(2) Proof that (3) $\implies$ (4):

$$\tilde{\phi}(0) \prod_{n=1}^{\infty} \tfrac{\sqrt{2}}{2} \breve{\mathrm{h}}\left(\frac{\omega}{2^n}\right) = \lim_{N \to \infty} \tilde{\phi}\left(\frac{\omega}{2^N}\right) \prod_{n=1}^{N} \tfrac{\sqrt{2}}{2} \breve{\mathrm{h}}\left(\frac{\omega}{2^n}\right) \qquad \text{by } \textit{continuity} \text{ and definition of } \prod_{n=1}^{\infty} x_n$$

$$= \tilde{\phi}(\omega) \qquad \text{by (3) and Lemma 5.7 page 67}$$

(3)  Proof that (2) $\Longleftarrow$ (4): by Lemma 5.8 page 69

☞

Definition 5.10 (next) formally defines the coefficients that appear in Theorem 5.6 (page 67).

**Definition 5.10**  Let $\left( L_{\mathbb{R}}^2, \left( V_j \right) \right)$ be a multiresolution analysis space with scaling function $\phi$. Let $\left( h_n \right)_{n \in \mathbb{Z}}$ be a sequence of coefficients such that $\phi = \sum_{n \in \mathbb{Z}} h_n \mathbf{D} \mathbf{T}^n \phi$.  A **multiresolution system** is the tuple $\left( L_{\mathbb{R}}^2, \left( V_j \right), \phi, \left( h_n \right) \right)$. The sequence $\left( h_n \right)_{n \in \mathbb{Z}}$ is the **scaling coefficient sequence**. A multiresolution system is also called an **MRA system**. An *MRA system* is an **orthonormal MRA system** if $\left\{ \mathbf{T}^n \phi |_{n \in \mathbb{Z}} \right\}$ is *orthonormal*.

**Definition 5.11**  Let $\left( L_{\mathbb{R}}^2, \left( V_j \right), \phi, \left( h_n \right) \right)$ be a multiresolution system, and $\mathbf{D}$ the dilation operator.  The **normalization coefficient at resolution** $n$ is the quantity
$$\left\| \mathbf{D}^j \phi \right\|.$$

**Theorem 5.12**  *Let* $\left( L_{\mathbb{R}}^2, \left( V_j \right), \phi, \left( h_n \right) \right)$ *be an* MRA system *(Definition 5.10 page 71). Let* span$A$ *be the* LINEAR SPAN *(Definition 3.2 page 32) of a set* $A$.

$$\underbrace{\operatorname{span}\left\{ \mathbf{T}^n \phi |_{n \in \mathbb{Z}} \right\} = V_0}_{\left\{ \mathbf{T}^n \phi |_{n \in \mathbb{Z}} \right\} \text{ is a BASIS for } V_0} \quad \Longrightarrow \quad \underbrace{\operatorname{span}\left\{ \mathbf{D}^j \mathbf{T}^n \phi |_{n \in \mathbb{Z}} \right\} = V_j}_{\left\{ \mathbf{D}^j \mathbf{T}^n \phi |_{n \in \mathbb{Z}} \right\} \text{ is a BASIS for } V_j} \quad \forall j \in \mathbb{W}$$

✎PROOF:  Proof is by induction:[135]

(1)  induction basis (proof for $j = 0$ case):

$$\operatorname{span}\left\{ \mathbf{D}^j \mathbf{T}^n \phi |_{n \in \mathbb{Z}} \right\} \Big|_{j=0} = \operatorname{span}\left\{ \mathbf{T}^n \phi |_{n \in \mathbb{Z}} \right\}$$
$$= V_0 \qquad\qquad \text{by left hypothesis}$$

---

[135]☞  [121], page 4

(2) induction step (proof that $j$ case $\implies$ $j+1$ case):

$$\mathrm{span}\left\{ \mathbf{D}^{j+1}\mathbf{T}^n\phi\,\middle|\,n\in\mathbb{Z} \right\}$$

$$= \left\{ f\in L^2_{\mathbb{R}} \,\middle|\, \exists\,(\alpha_n) \quad\text{such that}\quad f(x)=\sum_{n\in\mathbb{Z}}\alpha_n\mathbf{D}^{j+1}\mathbf{T}^n\phi \right\} \qquad \text{by def. of } \mathrm{span} \text{ (Definition 3.2 page 32)}$$

$$= \left\{ f\in L^2_{\mathbb{R}} \,\middle|\, \exists\,(\alpha_n) \quad\text{such that}\quad f(x)=\mathbf{D}\sum_{n\in\mathbb{Z}}\alpha_n\mathbf{D}^{j}\mathbf{T}^n\phi \right\}$$

$$= \left\{ f\in L^2_{\mathbb{R}} \,\middle|\, \exists\,(\alpha_n) \quad\text{such that}\quad \mathbf{D}^{-1}f(x)=\sum_{n\in\mathbb{Z}}\alpha_n\mathbf{D}^{j}\mathbf{T}^n\phi \right\}$$

$$= \left\{ [\mathbf{D}f]\in L^2_{\mathbb{R}} \,\middle|\, \exists\,(\alpha_n) \quad\text{such that}\quad \mathbf{D}^{-1}[\mathbf{D}f(x)]=\sum_{n\in\mathbb{Z}}\alpha_n\mathbf{D}^{j}\mathbf{T}^n\phi \right\}$$

$$= \mathbf{D}\left\{ f\in L^2_{\mathbb{R}} \,\middle|\, \exists\,(\alpha_n) \quad\text{such that}\quad f(x)=\sum_{n\in\mathbb{Z}}\alpha_n\mathbf{D}^{j}\mathbf{T}^n\phi \right\}$$

$$= \mathbf{D}\,\mathrm{span}\left\{ \mathbf{D}^{j}\mathbf{T}^n\phi\,\middle|\,n\in\mathbb{Z} \right\} \qquad \text{by def. of } \mathrm{span} \text{ (Definition 3.2 page 32)}$$

$$= \mathbf{D}\,V_j \qquad \text{by induction hypothesis}$$

$$= V_{j+1} \qquad \text{by } \textit{self-similar} \text{ property (Definition 5.1 page 63)}$$

## Example 5.13

In the *Haar* MRA, the scaling function $\phi(x)$ is the *pulse function*

$$\phi(x) = \begin{cases} 1 & \text{for } x\in[0,\,1) \\ 0 & \text{otherwise.} \end{cases}$$

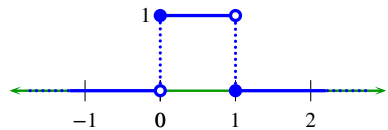

In the subspace $V_j$ $(j\in\mathbb{Z})$ the scaling functions are

$$\mathbf{D}^j\phi(x) = \begin{cases} (2)^{j/2} & \text{for } x\in\bigl[0,\,(2^{-j})\bigr) \\ 0 & \text{otherwise.} \end{cases}$$

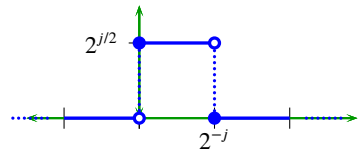

The scaling subspace $V_0$ is the span $V_0 \triangleq \mathrm{span}\left\{ \mathbf{T}^n\phi\,\middle|\,n\in\mathbb{Z} \right\}$. The scaling subspace $V_j$ is the

| subspace | transform | approximation |
|---|---|---|
| $V_0$ | | |
| $V_1$ | | |
| $V_2$ | | |

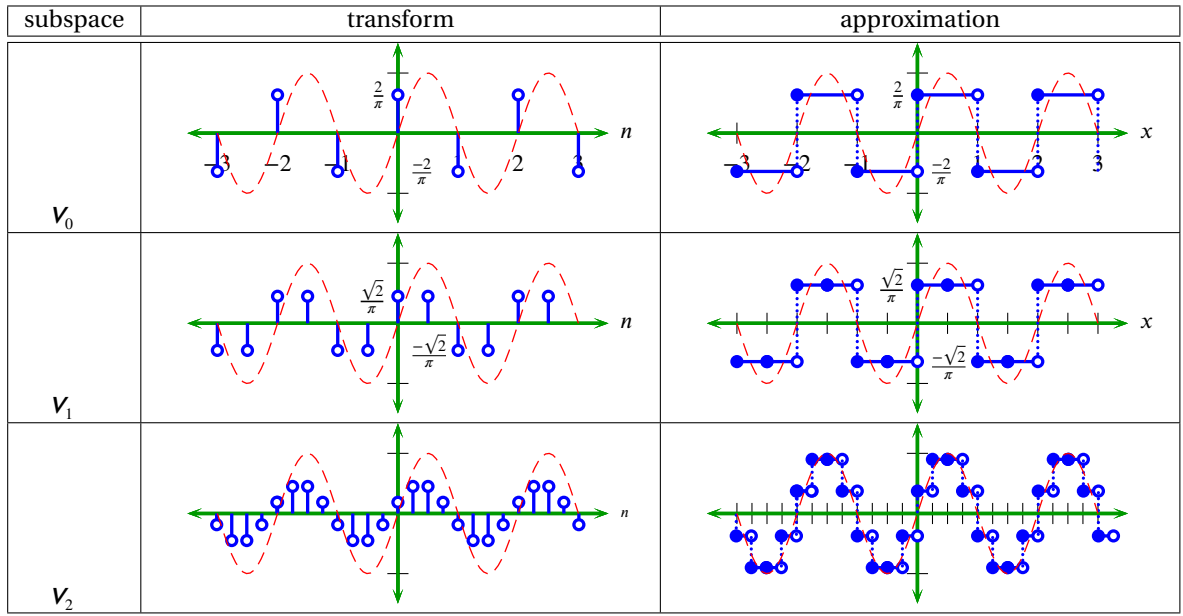

Figure 4: Example approximations of $\sin(\pi x)$ in 3 Haar scaling subspaces (see Example 5.13 page 72)

span $V_j \triangleq \operatorname{span}\left\{\mathbf{D}^j \mathbf{T}^n \phi \mid n \in \mathbb{Z}\right\}$. Note that $\left\|\mathbf{D}^j \mathbf{T}^n \phi\right\|$ for each resolution $j$ and shift $n$ is unity:

$$
\begin{aligned}
\left\|\mathbf{D}^j \mathbf{T}^n \phi\right\|^2 &= \|\phi\|^2 \\
&= \int_0^1 |1|^2 \, \mathrm{d}x && \text{by definition of } \|\cdot\| \text{ on } L_{\mathbb{R}}^2 \\
&= 1
\end{aligned}
$$

Let $f(x) = \sin(\pi x)$. Suppose we want to project $f(x)$ onto the subspaces $V_0$, $V_1$, $V_2$, ....

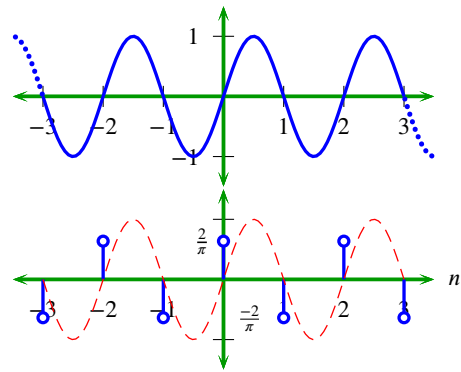

The values of the transform coefficients for the subspace $V_j$ are given by

$$\left[\mathbf{R}_j f(x)\right](n) = \frac{1}{\|\mathbf{D}^j \mathbf{T}^n \phi\|^2} \left\langle f(x) \mid \mathbf{D}^j \mathbf{T}^n \phi \right\rangle$$

$$= \frac{1}{\|\phi\|^2} \left\langle f(x) \mid 2^{j/2} \phi\left(2^j x - n\right) \right\rangle \qquad \text{by Proposition 4.4 page 40}$$

$$= 2^{j/2} \left\langle f(x) \mid \phi\left(2^j x - n\right) \right\rangle$$

$$= 2^{j/2} \int_{2^{-j} n}^{2^{-j}(n+1)} f(x)\, dx$$

$$= 2^{j/2} \int_{2^{-j} n}^{2^{-j}(n+1)} \sin(\pi x)\, dx$$

$$= 2^{j/2} \left(-\frac{1}{\pi}\right) \cos(\pi x)\Big|_{2^{-j} n}^{2^{-j}(n+1)}$$

$$= \frac{2^{j/2}}{\pi} \left[\cos\left(2^{-j} n\pi\right) - \cos\left(2^{-j}(n+1)\pi\right)\right]$$

And the projection $\mathbf{A}_n f(x)$ of the function $f(x)$ onto the subspace $V_j$ is

$$\mathbf{A}_j f(x) = \sum_{n \in \mathbb{Z}} \left\langle f(x) \mid \mathbf{D}^j \mathbf{T}^n \phi \right\rangle \mathbf{D}^j \mathbf{T}^n \phi$$

$$= \frac{2^{j/2}}{\pi} \sum_{n \in \mathbb{Z}} \left[\cos\left(2^{-j} n\pi\right) - \cos\left(2^{-j}(n+1)\pi\right)\right] 2^{j/2} \phi\left(2^j x - n\right)$$

$$= \frac{2^{j}}{\pi} \sum_{n \in \mathbb{Z}} \left[\cos\left(2^{-j} n\pi\right) - \cos\left(2^{-j}(n+1)\pi\right)\right] \phi\left(2^j x - n\right)$$

The transforms into the subspaces $V_0$, $V_1$, and $V_2$, as well as the approximations in those subspaces are as illustrated in Figure 4 (page 73).

### 5.1.3  Necessary Conditions

**Theorem 5.14** (admissibility condition)  *Let* $\hat{h}(z)$ *be the* Z-TRANSFORM *(Definition 2.35 page 24) and* $\check{h}(\omega)$ *the* DISCRETE-TIME FOURIER TRANSFORM *(Definition 2.38 page 24) of a sequence* $(h_n)_{n \in \mathbb{Z}}$.

$$\left\{ \left(L_{\mathbb{R}}^2, \left(V_j\right), \phi, (h_n)\right) \text{ is an MRA SYSTEM } \textit{(Definition 5.10 page 71)} \right\}$$

$$\underset{\Longleftarrow}{\overset{\Longrightarrow}{\phantom{x}}} \quad \underbrace{\left\{ \sum_{n \in \mathbb{Z}} h_n = \sqrt{2} \right\}}_{\textit{(1)} \text{ ADMISSIBILITY } \textit{in "time"}} \quad \Longleftrightarrow \quad \underbrace{\left\{ \hat{h}(z)\Big|_{z=1} = \sqrt{2} \right\}}_{\textit{(2)} \text{ ADMISSIBILITY } \textit{in "z domain"}} \quad \Longleftrightarrow \quad \underbrace{\left\{ \check{h}(\omega)\Big|_{\omega=0} = \sqrt{2} \right\}}_{\textit{(3)} \text{ ADMISSIBILITY } \textit{in "frequency"}}$$

PROOF:

(1) Proof that MRA system $\implies$ (1):

$$\sum_{n\in\mathbb{Z}} h_n = \frac{\int_{\mathbb{R}} \phi(x)\,\mathrm{d}x}{\int_{\mathbb{R}} \phi(x)\,\mathrm{d}x} \sum_{n\in\mathbb{Z}} h_n$$

$$= \frac{1}{\int_{\mathbb{R}} \phi(x)\,\mathrm{d}x} \int_{\mathbb{R}} \sum_{n\in\mathbb{Z}} h_n \phi(x)\,\mathrm{d}x$$

$$= \frac{1}{\int_{\mathbb{R}} \phi(x)\,\mathrm{d}x} \int_{\mathbb{R}} \sum_{n\in\mathbb{Z}} h_n \frac{\sqrt{2}}{\sqrt{2}} \phi(2y-n)2\,\mathrm{d}y \quad \text{let } y \triangleq \frac{x+n}{2} \implies x = 2y - n \implies \mathrm{d}x = 2\,\mathrm{d}y$$

$$= \frac{2}{\sqrt{2}} \frac{1}{\int_{\mathbb{R}} \phi(x)\,\mathrm{d}x} \int_{\mathbb{R}} \sum_{n\in\mathbb{Z}} h_n \mathbf{D}\mathbf{T}^n \phi(y)\,\mathrm{d}y \quad \text{by definitions of } \mathbf{T} \text{ and } \mathbf{D} \text{ (Definition 4.1 page 38)}$$

$$= \sqrt{2} \frac{1}{\int_{\mathbb{R}} \phi(x)\,\mathrm{d}x} \int_{\mathbb{R}} \phi(y)\,\mathrm{d}y \quad \text{by } \textit{dilation equation} \text{ (Theorem 5.6 page 67)}$$

$$= \sqrt{2}$$

(2) Alternate proof that MRA system $\implies$ (1):
Let $f(x) \triangleq 1 \quad \forall x \in \mathbb{R}$.

$$\langle \phi \mid f \rangle = \left\langle \sum_{n\in\mathbb{Z}} h_n \mathbf{D}\mathbf{T}^n \phi \;\middle|\; f \right\rangle \quad \text{by dilation equation (Theorem 5.6 page 67)}$$

$$= \sum_{n\in\mathbb{Z}} h_n \langle \mathbf{D}\mathbf{T}^n \phi \mid f \rangle \quad \text{by linearity of } \langle \triangle \mid \triangledown \rangle$$

$$= \sum_{n\in\mathbb{Z}} h_n \langle \phi \mid (\mathbf{D}\mathbf{T}^n)^* f \rangle \quad \text{by definition of operator adjoint (Theorem 1.35 page 10)}$$

$$= \sum_{n\in\mathbb{Z}} h_n \langle \phi \mid (\mathbf{T}^*)^n \mathbf{D}^* f \rangle \quad \text{by property of operator adjoint (Theorem 1.35 page 10)}$$

$$= \sum_{n\in\mathbb{Z}} h_n \langle \phi \mid (\mathbf{T}^{-1})^n \mathbf{D}^{-1} f \rangle \quad \text{by unitary property of } \mathbf{T} \text{ and } \mathbf{D} \text{ (Proposition 4.10 page 43)}$$

$$= \sum_{n\in\mathbb{Z}} h_n \langle \phi \mid (\mathbf{T}^{-1})^n \tfrac{\sqrt{2}}{2} f \rangle \quad \text{because } f \text{ is a constant and by Proposition 4.3 page 39}$$

$$= \sum_{n\in\mathbb{Z}} h_n \langle \phi \mid \tfrac{\sqrt{2}}{2} f \rangle \quad \text{by } f(x) = 1 \text{ definition}$$

$$= \sum_{n\in\mathbb{Z}} h_n \tfrac{\sqrt{2}}{2} \langle \phi \mid f \rangle \quad \text{by property of } \langle \triangle \mid \triangledown \rangle$$

$$= \tfrac{\sqrt{2}}{2} \langle \phi \mid f \rangle \sum_{n\in\mathbb{Z}} h_n$$

$$\implies \sum_{n\in\mathbb{Z}} h_n = \sqrt{2}$$

(3) Proof that (1) $\iff$ (2) $\iff$ (3): by Proposition 2.40 page 25.

(4)   Proof for $\Longleftarrow$ part: by Counterexample 5.15 page 76.

**Counterexample 5.15** Let $\left(L_{\mathbb{R}}^2, \left(\!\left( V_j \right)\!\right), \phi, \left(\!\left( h_n \right)\!\right)\right)$ be an *MRA system* (Definition 5.10 page 71).

$$\left\{ \left(\!\left( h_n \right)\!\right) \triangleq \sqrt{2}\bar{\delta}_{n-1} \triangleq \left\{ \begin{array}{ll} \sqrt{2} & \text{for } n = 1 \\ 0 & \text{otherwise.} \end{array} \right. \right\} \implies \{\phi(x) = 0\}$$

which means

$$\left\{ \sum_{n\in\mathbb{Z}} h_n = \sqrt{2} \right\} \implies \left\{ \left(L_{\mathbb{R}}^2, \left(\!\left( V_j \right)\!\right), \phi, \left(\!\left( h_n \right)\!\right)\right) \text{ is an MRA system for } L_{\mathbb{R}}^2. \right\}$$

✎ PROOF:

$$\begin{aligned}
\phi(x) &= \sum_{n\in\mathbb{Z}} h_n \mathbf{D}\mathbf{T}^n \phi(x) && \text{by } \textit{dilation equation} \text{ (Theorem 5.6 page 67)} \\
&= \sum_{n\in\mathbb{Z}} h_n \phi(2x - n) && \text{by definitions of } \mathbf{D} \text{ and } \mathbf{T} \text{ (Definition 4.1 page 38)} \\
&= \sum_{n\in\mathbb{Z}} \underbrace{\sqrt{2}\bar{\delta}_{n-1}}_{\left(\!\left( h_n \right)\!\right)} \phi(2x - n) && \text{by definitions of } \left(\!\left( h_n \right)\!\right) \\
&= \sqrt{2}\phi(2x - 1) && \text{by definition of } \phi(x) \\
\implies \phi(x) &= 0
\end{aligned}$$

This implies $\phi(x) = 0$, which implies that $\left(L_{\mathbb{R}}^2, \left(\!\left( V_j \right)\!\right), \phi, \left(\!\left( h_n \right)\!\right)\right)$ is *not* an *MRA system* for $L_{\mathbb{R}}^2$ because

$$\left( \bigcup_{j\in\mathbb{Z}} V_j \right)^{-} = \left( \bigcup_{j\in\mathbb{Z}} \text{span}\left\{ \mathbf{D}^j \mathbf{T}^n \phi \,|_{n\in\mathbb{Z}} \right\} \right)^{-} \neq L_{\mathbb{R}}^2$$

(the *least upper bound* is *not* $L_{\mathbb{R}}^2$).

**Theorem 5.16** (Quadrature condition in "time") *Let* $\left(L_{\mathbb{R}}^2, \left(\!\left( V_j \right)\!\right), \phi, \left(\!\left( h_n \right)\!\right)\right)$ *be an* MRA SYS-TEM *(Definition 5.10 page 71).*

$$\sum_{m\in\mathbb{Z}} h_m \sum_{k\in\mathbb{Z}} h_k^* \left\langle \phi \,|\, \mathbf{T}^{2n-m+k}\phi \right\rangle = \left\langle \phi \,|\, \mathbf{T}^n\phi \right\rangle \qquad \forall n\in\mathbb{Z}$$

✎ PROOF:

$$
\begin{aligned}
\langle \phi \mid \mathbf{T}^n \phi \rangle
&= \left\langle \sum_{m\in\mathbb{Z}} h_m \mathbf{D}\mathbf{T}^m \phi \;\middle|\; \mathbf{T}^n \sum_{k\in\mathbb{Z}} h_k \mathbf{D}\mathbf{T}^k \phi \right\rangle && \text{by dilation equation (Theorem 5.6 page 67)} \\
&= \sum_{m\in\mathbb{Z}} h_m \sum_{k\in\mathbb{Z}} h_k^* \left\langle \mathbf{D}\mathbf{T}^m \phi \mid \mathbf{T}^n \mathbf{D}\mathbf{T}^k \phi \right\rangle && \text{by properties of } \langle \triangle \mid \triangledown \rangle \\
&= \sum_{m\in\mathbb{Z}} h_m \sum_{k\in\mathbb{Z}} h_k^* \left\langle \phi \mid (\mathbf{D}\mathbf{T}^m)^* \, \mathbf{T}^n \mathbf{D}\mathbf{T}^k \phi \right\rangle && \text{by def. of operator adjoint (Proposition 1.33 page 10)} \\
&= \sum_{m\in\mathbb{Z}} h_m \sum_{k\in\mathbb{Z}} h_k^* \left\langle \phi \mid (\mathbf{D}\mathbf{T}^m)^* \, \mathbf{D}\mathbf{T}^{2n} \mathbf{T}^k \phi \right\rangle && \text{by Proposition 4.8 page 41} \\
&= \sum_{m\in\mathbb{Z}} h_m \sum_{k\in\mathbb{Z}} h_k^* \left\langle \phi \mid \mathbf{T}^{*m} \mathbf{D}^* \mathbf{D}\mathbf{T}^{2n} \mathbf{T}^k \phi \right\rangle && \text{by operator star-algebra prop. (Theorem 1.35 page 10)} \\
&= \sum_{m\in\mathbb{Z}} h_m \sum_{k\in\mathbb{Z}} h_k^* \left\langle \phi \mid \mathbf{T}^{-m} \mathbf{D}^{-1} \mathbf{D}\mathbf{T}^{2n} \mathbf{T}^k \phi \right\rangle && \text{by Proposition 4.10 page 43} \\
&= \sum_{m\in\mathbb{Z}} h_m \sum_{k\in\mathbb{Z}} h_k^* \left\langle \phi \mid \mathbf{T}^{2n-m+k} \phi \right\rangle
\end{aligned}
$$

         ✍

Theorem 5.24 (next) presents the *quadrature* necessary conditions of a wavelet system. These relations simplify dramatically in the special case of an *orthonormal wavelet system* (Theorem 2.44 page 27).

**Theorem 5.17** (Quadrature condition in "frequency") [136]
*Let* $\left( L^2_\mathbb{R}, \left( V_j \right), \phi, \left( h_n \right) \right)$ *be an* MRA SYSTEM *(Definition 5.10 page 71). Let* $\tilde{x}(\omega)$ *be the* DISCRETE TIME FOURIER TRANSFORM *for a sequence* $\left( x_n \right)_{n\in\mathbb{Z}}$ *in* $\ell^2_\mathbb{R}$. *Let* $\tilde{S}_{\phi\phi}(\omega)$ *be the* AUTO-POWER SPECTRUM *(Definition 4.27 page 51) of* $\phi$.

$$
\left| \breve{\mathrm{h}}(\omega) \right|^2 \tilde{S}_{\phi\phi}(\omega) + \left| \breve{\mathrm{h}}(\omega+\pi) \right|^2 \tilde{S}_{\phi\phi}(\omega+\pi) = 2\tilde{S}_{\phi\phi}(2\omega)
\qquad
\left( \begin{array}{l} \text{Note: } \tilde{S}_{\phi\phi}(\omega)=1 \\ \text{for ORTHONORMAL } \textit{MRA} \end{array} \right)
$$

✎ PROOF:

$2\tilde{S}_{\phi\phi}(2\omega)$

$$
= 2(2\pi) \sum_{n\in\mathbb{Z}} \left| \tilde{\phi}(2\omega+2\pi n) \right|^2 \qquad \text{by Theorem 4.28 page 51}
$$

$$
= 2(2\pi) \sum_{n\in\mathbb{Z}} \left| \tfrac{\sqrt{2}}{2} \breve{\mathrm{h}}\left( \frac{2\omega+2\pi n}{2} \right) \tilde{\phi}\left( \frac{2\omega+2\pi n}{2} \right) \right|^2 \qquad \text{by Lemma 5.7 page 67}
$$

$$
= 2\pi \sum_{n\in\mathbb{Z}_e} \left| \breve{\mathrm{h}}\left( \frac{2\omega+2\pi n}{2} \right) \right|^2 \left| \tilde{\phi}\left( \frac{2\omega+2\pi n}{2} \right) \right|^2 + 2\pi \sum_{n\in\mathbb{Z}_o} \left| \breve{\mathrm{h}}\left( \frac{2\omega+2\pi n}{2} \right) \right|^2 \left| \tilde{\phi}\left( \frac{2\omega+2\pi n}{2} \right) \right|^2
$$

$$
= 2\pi \sum_{n\in\mathbb{Z}} \left| \breve{\mathrm{h}}\left( \omega+2\pi n \right) \right|^2 \left| \tilde{\phi}\left( \omega+2\pi n \right) \right|^2 + 2\pi \sum_{n\in\mathbb{Z}} \left| \breve{\mathrm{h}}\left( \omega+2\pi n+\pi \right) \right|^2 \left| \tilde{\phi}\left( \omega+2\pi n+\pi \right) \right|^2
$$

---
[136] ✎ [27], page 135, ✎ [52], page 110

$$= 2\pi \sum_{n\in\mathbb{Z}} |\breve{\mathrm{h}}(\omega)|^2 \left|\tilde{\phi}(\omega+2\pi n)\right|^2 + 2\pi \sum_{n\in\mathbb{Z}} |\breve{\mathrm{h}}(\omega+\pi)|^2 \left|\tilde{\phi}(\omega+2\pi n+\pi)\right|^2 \qquad \text{by Proposition 2.39 page 25}$$

$$= |\breve{\mathrm{h}}(\omega)|^2 \left(2\pi \sum_{n\in\mathbb{Z}} |\tilde{\phi}(\omega+2\pi n)|^2\right) + |\breve{\mathrm{h}}(\omega+\pi)|^2 \left(2\pi \sum_{n\in\mathbb{Z}} |\tilde{\phi}(\omega+\pi+2\pi n)|^2\right)$$

$$= |\breve{\mathrm{h}}(\omega)|^2 \, \bar{S}_{\phi\phi}(\omega) + |\breve{\mathrm{h}}(\omega+\pi)|^2 \, \bar{S}_{\phi\phi}(\omega+\pi) \qquad \text{by Theorem 4.28 page 51}$$

### 5.1.4　Sufficient conditions

Theorem 5.18 (next) gives a set of *sufficient* conditions on the *scaling function* (Definition 5.1 page 63) $\phi$ to generate an *MRA*.

**Theorem 5.18** [137] *Let an* MRA *be defined as in Definition 5.1 page 63. Let a* Riesz sequence *be defined as in Definition 3.20 page 37. Let* $V_j \triangleq \overline{\mathrm{span}}\{\,\mathbf{T}^n\phi(x)|_{n\in\mathbb{Z}}\,\}$.

$$\left\{\begin{array}{ll} (1). & (\mathbf{T}^n\phi) \text{ is a Riesz sequence} \qquad \text{and} \\ (2). & \exists\,(h_n) \quad \text{such that} \quad \phi(x) = \sum_{n\in\mathbb{Z}} h_n \mathbf{DT}^n\phi(x) \quad \text{and} \\ (3). & \tilde{\phi}(\omega) \text{ is continuous } at\,0 \qquad\qquad\qquad \text{and} \\ (4). & \tilde{\phi}(0) \neq 0 \end{array}\right\} \implies \left\{\,\big(\!\big(V_j\big)\!\big)_{j\in\mathbb{Z}} \text{ is an MRA}\right\}$$

✎Proof:　For this to be true, each of the conditions in the definition of an *MRA* (Definition 5.1 page 63) must be satisfied:

(1)　Proof that each $V_j$ is *closed*: by definition of $\overline{\mathrm{span}}$

(2)　Proof that $\big(\!\big(V_j\big)\!\big)$ is *linearly ordered*:

$$V_j \subseteq V_{j+1} \iff \overline{\mathrm{span}}\{\mathbf{D}^j\mathbf{T}^n\phi\} \subseteq \overline{\mathrm{span}}\{\mathbf{D}^{j+1}\mathbf{T}^n\phi\} \qquad\qquad \iff (2)$$

(3)　Proof that $\bigcup_{j\in\mathbb{Z}} V_j$ is *dense* in $\pmb{L}^2_{\mathbb{R}}$: by Proposition 5.4 page 66

(4)　Proof of *self-similar* property:

$$\left\{\mathrm{f}\in V_j \iff \mathbf{D}\mathrm{f}\in V_{j+1}\right\} \iff \mathrm{f}\in\overline{\mathrm{span}}\{\mathbf{T}^n\phi\} \iff \mathbf{D}\mathrm{f}\in\overline{\mathrm{span}}\{\mathbf{DT}^n\phi\} \qquad \iff (2)$$

(5)　Proof for *Riesz basis*: by (1) and Proposition 5.4 page 66.

---

## 5.2 Wavelet analysis

### 5.2.1 Definition

The term "wavelet" comes from the French word "*ondelette*", meaning "small wave". And in essence, wavelets are "small waves" (as opposed to the "long waves" of Fourier analysis) that form a basis for the Hilbert space $L^2_{\mathbb{R}}$.[138]

**Definition 5.19** [139] Let $\mathbf{T}$ and $\mathbf{D}$ be as defined in Definition 4.1 page 38. A function $\psi(x)$ in $L^2_{\mathbb{R}}$ is a **wavelet function** for $L^2_{\mathbb{R}}$ if

$\left\{ \mathbf{D}^j \mathbf{T}^n \psi \,|_{j,n\in\mathbb{Z}} \right\}$ is a *Riesz basis* for $L^2_{\mathbb{R}}$.

In this case, $\psi$ is also called the **mother wavelet** of the basis $\left\{ \mathbf{D}^j \mathbf{T}^n \psi \,|_{j,n\in\mathbb{Z}} \right\}$. The sequence of subspaces $\left(\!\!\left( W_j \right)\!\!\right)_{j\in\mathbb{Z}}$ is the **wavelet analysis** induced by $\psi$, where each subspace $W_j$ is defined as

$W_j \triangleq \operatorname{span}\left\{ \mathbf{D}^j \mathbf{T}^n \psi \,|_{n\in\mathbb{Z}} \right\}$ .

A *wavelet analysis* $\left(\!\!\left( W_j \right)\!\!\right)$ is often constructed from a *multiresolution anaysis* (Definition 5.1 page 63) $\left(\!\!\left( V_j \right)\!\!\right)$ under the relationship

$V_{j+1} = V_j \dotplus W_j$,     where $\dotplus$ is subspace addition (*Minkowski addition*).

By this relationship alone, $\left(\!\!\left( W_j \right)\!\!\right)$ is in no way uniquely defined in terms of a multiresolution analysis $\left(\!\!\left( V_j \right)\!\!\right)$. In general there are many possible complements of a subspace $V_j$. To uniquely define such a wavelet subspace, one or more additional constraints are required. One of the most common additional constraints is *orthogonality*, such that $V_j$ and $W_j$ are orthogonal to each other.

### 5.2.2 Dilation equation

Suppose $\left(\!\!\left( \mathbf{T}^n \psi \right)\!\!\right)_{n\in\mathbb{Z}}$ is a basis for $W_0$. By Definition 5.19 page 79, the wavelet subspace $W_0$ is contained in the scaling subspace $V_1$. By Definition 5.1 page 63, the sequence $\left(\!\!\left( \mathbf{D}\mathbf{T}^n \phi \right)\!\!\right)_{n\in\mathbb{Z}}$ is a basis for $V_1$. Because $W_0$ is contained in $V_1$, the sequence $\left(\!\!\left( \mathbf{D}\mathbf{T}^n \phi \right)\!\!\right)_{n\in\mathbb{Z}}$ is also a basis for $W_0$.

**Theorem 5.20** *Let* $\left( L^2_{\mathbb{R}}, \left(\!\!\left( V_j \right)\!\!\right), \phi, \left(\!\!\left( h_n \right)\!\!\right) \right)$ *be a multiresolution system and* $\left(\!\!\left( W_j \right)\!\!\right)_{j\in\mathbb{Z}}$ *a wavelet analysis with respect to* $\left( L^2_{\mathbb{R}}, \left(\!\!\left( V_j \right)\!\!\right), \phi, \left(\!\!\left( h_n \right)\!\!\right) \right)$ *and with wavelet function* $\psi$ .

$\exists \left(\!\!\left( g_n \right)\!\!\right)_{n\in\mathbb{Z}}$    such that    $\psi = \displaystyle\sum_{n\in\mathbb{Z}} g_n \mathbf{D}\mathbf{T}^n \phi$

---

✎ PROOF:

$$\psi \in W_0 \qquad \text{by Definition 5.19 page 79}$$
$$\subseteq V_1 \qquad \text{by Definition 5.19 page 79}$$
$$= \operatorname{span}(\mathbf{DT}^n \phi(x))_{n \in \mathbb{Z}} \qquad \text{by Definition 5.1 page 63 (MRA)}$$
$$\implies \exists (g_n)_{n \in \mathbb{Z}} \quad \text{such that} \quad \psi = \sum_{n \in \mathbb{Z}} g_n \mathbf{DT}^n \phi$$

⇦

A *wavelet system* (next definition) consists of two subspace sequences:

- ✿ A **multiresolution analysis** $(V_j)$ (Definition 5.1 page 63) provides "coarse" approximations of a function in $L^2_{\mathbb{R}}$ at different "scales" or resolutions.
- ✿ A **wavelet analysis** $(W_j)$ provides the "detail" of the function missing from the approximation provided by a given scaling subspace (Definition 5.19 page 79).

**Definition 5.21** Let $\left(L^2_{\mathbb{R}}, (V_j), \phi, (h_n)\right)$ be a multiresolution system (Definition 5.1 page 63) and $(W_j)_{j \in \mathbb{Z}}$ a wavelet analysis (Definition 5.19 page 79) with respect to $(V_j)$. Let $(g_n)_{n \in \mathbb{Z}}$ be a sequence of coefficients such that $\psi = \sum_{n \in \mathbb{Z}} g_n \mathbf{DT}^n \phi$. A **wavelet system** is the tuple $\left(L^2_{\mathbb{R}}, (V_j), (W_j), \phi, \psi, (h_n), (g_n)\right)$ and the sequence $(g_n)_{n \in \mathbb{Z}}$ is the **wavelet coefficient sequence**.

### 5.2.3 Necessary conditions

**Theorem 5.22** (quadrature conditions in "time") *Let* $\left(L^2_{\mathbb{R}}, (V_j), (W_j), \phi, \psi, (h_n), (g_n)\right)$ *be a wavelet system* (Definition 5.21 page 80).

1. $\displaystyle \sum_{m \in \mathbb{Z}} h_m \sum_{k \in \mathbb{Z}} h_k^* \left\langle \phi \,\middle|\, \mathbf{T}^{2n-m+k}\phi \right\rangle = \left\langle \phi \,\middle|\, \mathbf{T}^n \phi \right\rangle \quad \forall n \in \mathbb{Z}$

2. $\displaystyle \sum_{m \in \mathbb{Z}} g_m \sum_{k \in \mathbb{Z}} g_k^* \left\langle \phi \,\middle|\, \mathbf{T}^{2n-m+k}\phi \right\rangle = \left\langle \psi \,\middle|\, \mathbf{T}^n \psi \right\rangle \quad \forall n \in \mathbb{Z}$

3. $\displaystyle \sum_{m \in \mathbb{Z}} h_m \sum_{k \in \mathbb{Z}} g_k^* \left\langle \phi \,\middle|\, \mathbf{T}^{2n-m+k}\phi \right\rangle = \left\langle \phi \,\middle|\, \mathbf{T}^n \psi \right\rangle \quad \forall n \in \mathbb{Z}$

✎ PROOF:

(1)  Proof for (1): by Theorem 5.16 page 76.

(2) Proof for (2):

$\langle \psi \mid \mathbf{T}^n \psi \rangle$

$\displaystyle = \left\langle \sum_{m\in\mathbb{Z}} g_m \mathbf{DT}^m \phi \,\middle|\, \mathbf{T}^n \sum_{k\in\mathbb{Z}} g_k \mathbf{DT}^k \phi \right\rangle$    by Theorem **5.20** page **79**

$\displaystyle = \sum_{m\in\mathbb{Z}} g_m \sum_{k\in\mathbb{Z}} g_k^* \langle \mathbf{DT}^m \phi \mid \mathbf{T}^n \mathbf{DT}^k \phi \rangle$    by properties of $\langle \triangle \mid \triangledown \rangle$

$\displaystyle = \sum_{m\in\mathbb{Z}} g_m \sum_{k\in\mathbb{Z}} g_k^* \langle \phi \mid (\mathbf{DT}^m)^* \mathbf{T}^n \mathbf{DT}^k \phi \rangle$    by definition of operator adjoint (Proposition **1.33** page **10**)

$\displaystyle = \sum_{m\in\mathbb{Z}} g_m \sum_{k\in\mathbb{Z}} g_k^* \langle \phi \mid (\mathbf{DT}^m)^* \mathbf{DT}^{2n} \mathbf{T}^k \phi \rangle$    by Proposition **4.8** page **41**

$\displaystyle = \sum_{m\in\mathbb{Z}} g_m \sum_{k\in\mathbb{Z}} g_k^* \langle \phi \mid \mathbf{T}^{*m} \mathbf{D}^* \mathbf{DT}^{2n} \mathbf{T}^k \phi \rangle$    by operator star-algebra properties (Theorem **1.35** page **10**)

$\displaystyle = \sum_{m\in\mathbb{Z}} g_m \sum_{k\in\mathbb{Z}} g_k^* \langle \phi \mid \mathbf{T}^{-m} \mathbf{D}^{-1} \mathbf{DT}^{2n} \mathbf{T}^k \phi \rangle$    by Proposition **4.10** page **43**

$\displaystyle = \sum_{m\in\mathbb{Z}} g_m \sum_{k\in\mathbb{Z}} g_k^* \langle \phi \mid \mathbf{T}^{2n-m+k} \phi \rangle$

(3) Proof for (3):

$\langle \phi \mid \mathbf{T}^n \psi \rangle$

$\displaystyle = \left\langle \sum_{m\in\mathbb{Z}} h_m \mathbf{DT}^m \phi \,\middle|\, \mathbf{T}^n \sum_{k\in\mathbb{Z}} g_k \mathbf{DT}^k \phi \right\rangle$    by Theorem **5.6** page **67** and Theorem **5.20** page **79**

$\displaystyle = \sum_{m\in\mathbb{Z}} h_m \sum_{k\in\mathbb{Z}} g_k^* \langle \mathbf{DT}^m \phi \mid \mathbf{T}^n \mathbf{DT}^k \phi \rangle$    by properties of $\langle \triangle \mid \triangledown \rangle$

$\displaystyle = \sum_{m\in\mathbb{Z}} h_m \sum_{k\in\mathbb{Z}} g_k^* \langle \phi \mid (\mathbf{DT}^m)^* \mathbf{T}^n \mathbf{DT}^k \phi \rangle$    by definition of operator adjoint (Proposition **1.33** page **10**)

$\displaystyle = \sum_{m\in\mathbb{Z}} h_m \sum_{k\in\mathbb{Z}} g_k^* \langle \phi \mid (\mathbf{DT}^m)^* \mathbf{DT}^{2n} \mathbf{T}^k \phi \rangle$    by Proposition **4.8** page **41**

$\displaystyle = \sum_{m\in\mathbb{Z}} h_m \sum_{k\in\mathbb{Z}} g_k^* \langle \phi \mid \mathbf{T}^{*m} \mathbf{D}^* \mathbf{DT}^{2n} \mathbf{T}^k \phi \rangle$    by operator star-algebra properties (Theorem **1.35** page **10**)

$\displaystyle = \sum_{m\in\mathbb{Z}} h_m \sum_{k\in\mathbb{Z}} g_k^* \langle \phi \mid \mathbf{T}^{-m} \mathbf{D}^{-1} \mathbf{DT}^{2n} \mathbf{T}^k \phi \rangle$    by Proposition **4.10** page **43**

$\displaystyle = \sum_{m\in\mathbb{Z}} h_m \sum_{k\in\mathbb{Z}} g_k^* \langle \phi \mid \mathbf{T}^{2n-m+k} \phi \rangle$

**Proposition 5.23** *Let* $\left( L^2_{\mathbb{R}}, \left( V_j \right), \left( W_j \right), \phi, \psi, \left( h_n \right), \left( g_n \right) \right)$ *be a wavelet system. Let* $\tilde{\phi}(\omega)$ *and* $\tilde{\psi}(\omega)$ *be the* FOURIER TRANSFORM*s of* $\phi(x)$ *and* $\psi(x)$, *respectively. Let* $\breve{g}(\omega)$ *be the* DISCRETE TIME FOURIER TRANSFORM *of* $\left( g_n \right)$.

$$\tilde{\psi}(\omega) = \frac{\sqrt{2}}{2} \breve{g}\left( \frac{\omega}{2} \right) \tilde{\phi}\left( \frac{\omega}{2} \right)$$

✎ PROOF:

$$
\begin{aligned}
\tilde{\psi}(\omega) \triangleq \tilde{\mathbf{F}}\psi \\
&= \tilde{\mathbf{F}}\sum_{n\in\mathbb{Z}} g_n \mathbf{D}\mathbf{T}^n\phi && \text{by Theorem 5.20 page 79} \\
&= \sum_{n\in\mathbb{Z}} g_n \tilde{\mathbf{F}}\mathbf{D}\mathbf{T}^n\phi \\
&= \sum_{n\in\mathbb{Z}} g_n \mathbf{D}^{-1}\tilde{\mathbf{F}}\mathbf{T}^n\phi && \text{by Corollary 4.18 page 47} \\
&= \sum_{n\in\mathbb{Z}} g_n \mathbf{D}^{-1} e^{-i\omega n}\tilde{\mathbf{F}}\phi && \text{by Corollary 4.18 page 47} \\
&= \sum_{n\in\mathbb{Z}} g_n \sqrt{2}\left(\mathbf{D}^{-1} e^{-i\omega n}\right)\left(\mathbf{D}^{-1}\tilde{\mathbf{F}}\phi\right) && \text{by Proposition 4.7 page 41} \\
&= \sqrt{2}\left(\mathbf{D}^{-1}\sum_{n\in\mathbb{Z}} g_n e^{-i\omega n}\right)\left(\mathbf{D}^{-1}\tilde{\mathbf{F}}\phi\right) \\
&= \sqrt{2}\left(\mathbf{D}^{-1}\breve{\mathbf{F}}(g_n)\right)\left(\mathbf{D}^{-1}\tilde{\mathbf{F}}\phi\right) && \text{by definition of } \breve{\mathbf{F}} \\
&= \sqrt{2}\frac{\sqrt{2}}{2}\breve{g}\left(\frac{\omega}{2}\right)\frac{\sqrt{2}}{2}\tilde{\phi}\left(\frac{\omega}{2}\right) && \text{by Proposition 4.3 page 39} \\
&= \frac{\sqrt{2}}{2}\breve{g}\left(\frac{\omega}{2}\right)\tilde{\phi}\left(\frac{\omega}{2}\right)
\end{aligned}
$$

✎

Theorem 5.24 (next) presents the *quadrature* necessary conditions of a wavelet system. These relations simplify dramatically in the special case of an *orthonormal wavelet system* (Theorem 2.44 page 27).

**Theorem 5.24** (Quadrature conditions in "frequency")   [140]   *Let* $\left(\boldsymbol{L}^2_{\mathbb{R}},\,(\!(\,V_j\,)\!)\,\right)$, $(\!(\,W_j\,)\!)$, $\phi$, $\psi$, $(\!(h_n)\!)$, $(\!(g_n)\!)$ *be a wavelet system. Let* $\breve{x}(\omega)$ *be the* DISCRETE TIME FOURIER TRANSFORM *for a sequence* $(x_n)_{n\in\mathbb{Z}}$ *in* $\boldsymbol{\ell}^2_{\mathbb{R}}$. *Let* $\tilde{S}_{\phi\phi}(\omega)$ *be the* AUTO-POWER SPECTRUM *(Definition 4.27 page 51) of* $\phi$, $\tilde{S}_{\psi\psi}(\omega)$ *be the* AUTO-POWER SPECTRUM *of* $\psi$, *and* $\tilde{S}_{\phi\psi}(\omega)$ *be the* CROSS-POWER SPECTRUM *of* $\phi$ *and* $\psi$.

1. $\left|\breve{h}(\omega)\right|^2 \tilde{S}_{\phi\phi}(\omega) + \left|\breve{h}(\omega+\pi)\right|^2 \tilde{S}_{\phi\phi}(\omega+\pi) \quad = \quad 2\tilde{S}_{\phi\phi}(2\omega)$
2. $\left|\breve{g}(\omega)\right|^2 \tilde{S}_{\phi\phi}(\omega) + \left|\breve{g}(\omega+\pi)\right|^2 \tilde{S}_{\phi\phi}(\omega+\pi) \quad = \quad 2\tilde{S}_{\psi\psi}(2\omega)$
3. $\breve{h}(\omega)\breve{g}^*(\omega)\tilde{S}_{\phi\phi}(\omega) + \breve{h}(\omega+\pi)\breve{g}^*(\omega+\pi)\tilde{S}_{\phi\phi}(\omega+\pi) \quad = \quad 2\tilde{S}_{\phi\psi}(2\omega)$

✎ PROOF:

(1)   Proof for (1): by Theorem 5.17 page 77.

---

[140] ✎ [27], page 135, ✎ [52], page 110

(2)  Proof for (2):

$$2\tilde{S}_{\psi\psi}(2\omega) \triangleq 2(2\pi)\sum_{n\in\mathbb{Z}}|\breve{\psi}(2\omega+2\pi n)|^2$$

$$= 2(2\pi)\sum_{n\in\mathbb{Z}}\left|\frac{\sqrt{2}}{2}\breve{g}\left(\frac{2\omega+2\pi n}{2}\right)\tilde{\phi}\left(\frac{2\omega+2\pi n}{2}\right)\right|^2 \qquad \text{by Lemma 5.7 page 67}$$

$$= 2\pi\sum_{n\in\mathbb{Z}_e}\left|\breve{g}\left(\frac{2\omega+2\pi n}{2}\right)\right|^2\left|\tilde{\phi}\left(\frac{2\omega+2\pi n}{2}\right)\right|^2 +$$

$$\qquad 2\pi\sum_{n\in\mathbb{Z}_o}\left|\breve{g}\left(\frac{2\omega+2\pi n}{2}\right)\right|^2\left|\tilde{\phi}\left(\frac{2\omega+2\pi n}{2}\right)\right|^2$$

$$= 2\pi\sum_{n\in\mathbb{Z}}|\breve{g}(\omega+2\pi n)|^2\left|\tilde{\phi}(\omega+2\pi n)\right|^2 + 2\pi\sum_{n\in\mathbb{Z}}|\breve{g}(\omega+2\pi n+\pi)|^2\left|\tilde{\phi}(\omega+2\pi n+\pi)\right|^2$$

$$= 2\pi\sum_{n\in\mathbb{Z}}|\breve{g}(\omega)|^2\left|\tilde{\phi}(\omega+2\pi n)\right|^2 + 2\pi\sum_{n\in\mathbb{Z}}|\breve{g}(\omega+\pi)|^2\left|\tilde{\phi}(\omega+2\pi n+\pi)\right|^2$$

$$= |\breve{g}(\omega)|^2\left(2\pi\sum_{n\in\mathbb{Z}}\left|\tilde{\phi}(\omega+2\pi n)\right|^2 +\right)|\breve{g}(\omega+\pi)|^2\left(2\pi\sum_{n\in\mathbb{Z}}\left|\tilde{\phi}(\omega+\pi+2\pi n)\right|^2\right)$$

$$= |\breve{g}(\omega)|^2\tilde{S}_{\phi\phi}(\omega) + |\breve{g}(\omega+\pi)|^2\tilde{S}_{\phi\phi}(\omega+\pi) \qquad \text{by Theorem 4.28 page 51}$$

(3)  Proof for (3):

$$2\tilde{S}_{\phi\psi}(2\omega) = 2(2\pi)\sum_{n\in\mathbb{Z}}\tilde{\phi}(2\omega+2\pi n)\breve{\psi}^*(2\omega+2\pi n)$$

$$= 2(2\pi)\sum_{n\in\mathbb{Z}}\frac{\sqrt{2}}{2}\breve{h}(\omega+\pi n)\tilde{\phi}(\omega+\pi n)\frac{\sqrt{2}}{2}\breve{g}^*(\omega+\pi n)\tilde{\phi}^*(\omega+\pi n) \qquad \text{by Lemma 5.7 page 67}$$

$$= 2\pi\sum_{n\in\mathbb{Z}}\breve{h}(\omega+\pi n)\breve{g}^*(\omega+\pi n)\left|\tilde{\phi}(\omega+\pi n)\right|^2$$

$$= 2\pi\sum_{n\in\mathbb{Z}_o}\breve{h}(\omega+\pi n)\breve{g}^*(\omega+\pi n)\left|\tilde{\phi}(\omega+\pi n)\right|^2$$

$$\qquad + 2\pi\sum_{n\in\mathbb{Z}_e}\breve{h}(\omega+\pi n)\breve{g}^*(\omega+\pi n)\left|\tilde{\phi}(\omega+\pi n)\right|^2$$

$$= 2\pi\sum_{n\in\mathbb{Z}}\breve{h}(\omega+2\pi n+\pi)\breve{g}^*(\omega+2\pi n+\pi)\left|\tilde{\phi}(\omega+2\pi n+\pi)\right|^2$$

$$\qquad + 2\pi\sum_{n\in\mathbb{Z}}\breve{h}(\omega+2\pi n)\breve{g}^*(\omega+2\pi n)\left|\tilde{\phi}(\omega+2\pi n)\right|^2$$

$$= 2\pi\sum_{n\in\mathbb{Z}}\breve{h}(\omega+\pi)\breve{g}^*(\omega+\pi)\left|\tilde{\phi}(\omega+2\pi n+\pi)\right|^2 + 2\pi\sum_{n\in\mathbb{Z}}\breve{h}(\omega)\breve{g}^*(\omega)\left|\tilde{\phi}(\omega+2\pi n)\right|^2$$

$$= \breve{h}(\omega)\breve{g}^*(\omega)\left(2\pi\sum_{n\in\mathbb{Z}}\left|\tilde{\phi}(\omega+2\pi n)\right|^2\right)$$

$$\qquad + \breve{h}(\omega+\pi)\breve{g}^*(\omega+\pi)\left(2\pi\sum_{n\in\mathbb{Z}}\left|\tilde{\phi}(\omega+\pi+2\pi n)\right|^2\right)$$

$$= \breve{h}(\omega)\breve{g}^*(\omega)\left(2\pi \sum_{n\in\mathbb{Z}} |\tilde{\phi}(\omega + 2\pi n)|^2\right) + \breve{h}(\omega + \pi)\breve{g}^*(\omega + \pi)\left(2\pi \sum_{n\in\mathbb{Z}} |\tilde{\phi}(\omega + \pi + 2\pi n)|^2\right)$$

$$= \breve{h}(\omega)\breve{g}^*(\omega)\tilde{S}_{\phi\phi}(\omega) + \breve{h}(\omega + \pi)\breve{g}^*(\omega + \pi)\tilde{S}_{\phi\phi}(\omega + \pi) \qquad \text{by Theorem 4.28 page 51}$$

### 5.2.4 Sufficient condition

In this text, an often used sufficient condition for designing the *wavelet coefficient sequence* $(g_n)$ (Definition 5.21 page 80) is the *conjugate quadrature filter condition*. It expresses the sequence $(g_n)$ in terms of the *scaling coefficient sequence* (Definition 5.10 page 71) and a "shift" integer $N$ as $g_n = \pm(-1)^n h_{N-n}^*$.

**Theorem 5.25** *Let* $\left(L_\mathbb{R}^2, (V_j), (W_j), \phi, \psi, (h_n), (g_n)\right)$ *be a* WAVELET SYSTEM *(Definition 5.21 page 80). Let* $\breve{g}(\omega)$ *be the* DTFT *(Definition 2.38 page 24) and* $\hat{g}(z)$ *the* Z-TRANSFORM *(Definition 2.35 page 24) of* $(g_n)$.

$$\underbrace{g_n = \pm(-1)^n h_{N-n}^*, \; N\in\mathbb{Z}}_{\text{CONJUGATE QUADRATURE FILTER}} \quad \Longleftrightarrow \quad \breve{g}(\omega) = \pm(-1)^N e^{-i\omega N} \breve{h}^*(\omega + \pi)\Big|_{\omega=\pi} \qquad (1)$$

$$\Longrightarrow \quad \sum_{n\in\mathbb{Z}} (-1)^n g_n = \sqrt{2} \qquad (2)$$

$$\Longleftrightarrow \quad \hat{g}(z)\Big|_{z=-1} = \sqrt{2} \qquad (3)$$

$$\Longleftrightarrow \quad \breve{g}(\omega)\Big|_{\omega=\pi} = \sqrt{2} \qquad (4)$$

✎PROOF:

(1) Proof that CQF $\Longleftrightarrow$ (1): by Theorem 2.46 page 30

(2) Proof that CQF $\Longrightarrow$ (4):

$$\breve{g}(\pi) = \breve{g}(\omega)\Big|_{\omega=\pi}$$

$$= \pm(-1)^N e^{-i\omega N} \breve{h}^*(\omega + \pi)\Big|_{\omega=\pi} \qquad \text{by Theorem 2.46 page 30}$$

$$= \pm(-1)^N e^{-i\pi N} \breve{h}^*(2\pi)$$

$$= \pm(-1)^N (-1)^N \breve{h}^*(0) \qquad \text{by Proposition 2.39 page 25}$$

$$= \sqrt{2} \qquad \text{by } \textit{admissibility condition (Theorem 5.14 page 74)}$$

(3) Proof that (2) $\Longleftrightarrow$ (3) $\Longleftrightarrow$ (4): by Proposition 2.41 page 26

## 5.3 Partition of unity systems

### 5.3.1 Motivation

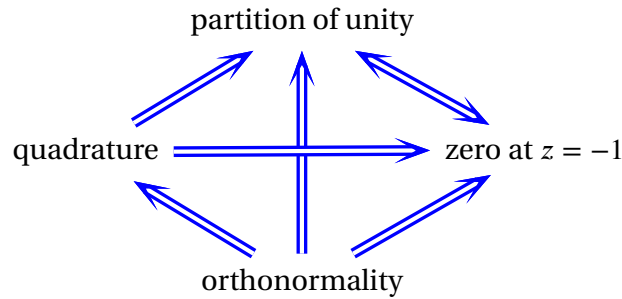

Figure 5: Implications of scaling function properties

A very common property of scaling functions (Definition 5.1 page 63) is the *partition of unity* property (Definition 5.26 page 86). The partition of unity is a kind of generalization of *orthonormality*; that is, *all* orthonormal scaling functions form a partition of unity. But the partition of unity property is not just a consequence of orthonormality, but also a generalization of orthonormality, in that if you remove the orthonormality constraint, the partition of unity is still a reasonable constraint in and of itself.

There are two reasons why the partition of unity property is a reasonable constraint on its own:

   ☞ Without a partition of unity, it is difficult to represent a function as simple as a constant.[141]

   ☞ For a multiresolution system $\left( L_{\mathbb{R}}^2, \left( V_j \right), \phi, \left( h_n \right) \right)$, the partition of unity property is equivalent to $\sum_{n \in \mathbb{Z}} (-1)^n h_n = 0$ (Theorem 5.32 page 88). As viewed from the perspective of discrete time signal processing, this implies that the scaling coefficients form a "*lowpass filter*"; lowpass filters provide a kind of "coarse approximation" of a function. And that is what the scaling function is "supposed" to do—to provide a coarse approximation at some resolution or "scale" (Definition 5.1 page 63).

---

[141] ☞ [73], page 8

### 5.3.2   Definition and results

**Definition 5.26** [142]
A function $f \in \mathbb{R}^{\mathbb{R}}$ forms a **partition of unity** if
$$\sum_{n \in \mathbb{Z}} \mathbf{T}^n f(x) = 1 \qquad \forall x \in \mathbb{R}.$$

**Theorem 5.27** [143] *Let* $\left(\mathbf{L}^2_{\mathbb{R}}, \left(V_j\right), \phi, \left(h_n\right)\right)$ *be a multiresolution system (Definition 5.10 page 71).* *Let* $\tilde{\mathbf{F}}f(\omega)$ *be the* FOURIER TRANSFORM *(Definition 2.22 page 20) of a function* $f \in \mathbf{L}^2_{\mathbb{R}}$. *Let* $\bar{\delta}_n$ *be the* KRONECKER DELTA FUNCTION.

$$\underbrace{\sum_{n \in \mathbb{Z}} \mathbf{T}^n f = c}_{\text{PARTITION OF UNITY } in \text{ "time"}} \qquad \Longleftrightarrow \qquad \underbrace{\left[\tilde{\mathbf{F}}f\right](2\pi n) = \bar{\delta}_n}_{\text{PARTITION OF UNITY } in \text{ "frequency"}}$$

✎PROOF:   Let $\mathbb{Z}_e$ be the set of even integers and $\mathbb{Z}_o$ the set of odd integers.

(1)   Proof for ( $\Longrightarrow$ ) case:

$$c = \sum_{m \in \mathbb{Z}} \mathbf{T}^m f(x) \qquad\qquad\qquad \text{by left hypothesis}$$

$$= \sum_{m \in \mathbb{Z}} f(x - m) \qquad\qquad\qquad \text{by definition of } \mathbf{T} \text{ (Definition 4.1 page 38)}$$

$$= \sqrt{2\pi} \sum_{m \in \mathbb{Z}} \tilde{f}(2\pi m) e^{i 2\pi m x} \qquad\qquad \text{by } PSF \text{ (Theorem 4.21 page 48)}$$

$$= \underbrace{\sqrt{2\pi} \tilde{f}(2\pi n) e^{i 2\pi n x}}_{\text{real and constant for } n = 0} + \sqrt{2\pi} \sum_{m \in \mathbb{Z}_n} \tilde{f}(2\pi m) e^{i 2\pi m x}}_{\text{complex and non-constant}}$$

$$\Longrightarrow \sqrt{2\pi} \tilde{f}(2\pi n) = c \bar{\delta}_n \qquad\qquad \text{because } c \text{ is real and constant for all } t$$

(2) Proof for ( $\Longleftarrow$ ) case:

$$\sum_{n\in\mathbb{Z}} \mathbf{T}^n f(x) = \sum_{n\in\mathbb{Z}} f(x-n) \qquad \text{by definition of } \mathbf{T} \; \text{(Definition 4.1 page 38)}$$

$$= \sqrt{2\pi} \sum_{n\in\mathbb{Z}} \tilde{f}(2\pi n) e^{-i2\pi nx} \qquad \text{by } PSF \; \text{(Theorem 4.21 page 48)}$$

$$= \sqrt{2\pi} \sum_{n\in\mathbb{Z}} \frac{c}{\sqrt{2\pi}} \bar{\delta}_n e^{-i2\pi nx} \qquad \text{by right hypothesis}$$

$$= \sqrt{2\pi} \frac{c}{\sqrt{2\pi}} e^{-i2\pi 0x} \qquad \text{by definition of } \bar{\delta}_n \; \text{(Definition 3.12 page 35)}$$

$$= c$$

**Corollary 5.28**

$$\left\{ \begin{array}{l} \exists \mathbf{g} \in L^2_{\mathbb{R}} \; \textit{such that} \\ f(x) = \mathbb{1}_{[-1,\,1)}(x) \star \mathbf{g}(x) \end{array} \right\} \implies \left\{ \begin{array}{l} f(x) \; \textit{generates} \\ \textit{a} \; \text{PARTITION OF UNITY} \end{array} \right\}$$

**Example 5.29** All *B-splines* form a partition of unity. All B-splines of order $n = 1$ or greater can be generated by convolution with a *pulse* function, similar to that specified in Corollary 5.28 (page 87).

**Example 5.30** Let a function $f$ be defined in terms of the cosine function (Definition 2.5 page 16) as follows:

$$f(x) \triangleq \left\{ \begin{array}{ll} \cos^2\left(\frac{\pi}{2}x\right) & \text{for } |x| \leq 1 \\ 0 & \text{otherwise} \end{array} \right.$$

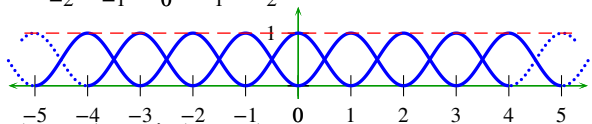

Then $f$ forms a *partition of unity*:

Note that $\tilde{f}(\omega) = \dfrac{1}{2\sqrt{2\pi}} \Bigg[ \underbrace{\dfrac{2\sin\omega}{\omega}}_{2\operatorname{sinc}(\omega)} + \underbrace{\dfrac{\sin(\omega-\pi)}{(\omega-\pi)}}_{\operatorname{sinc}(\omega-\pi)} + \underbrace{\dfrac{\sin(\omega+\pi)}{(\omega+\pi)}}_{\operatorname{sinc}(\omega+\pi)} \Bigg]$

and so $\tilde{f}(2\pi n) = \frac{1}{\sqrt{2\pi}} \bar{\delta}_n$:

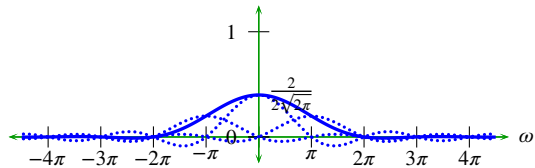

**Example 5.31** (raised cosine) [144] Let a function f be defined in terms of the cosine function (Definition 2.5 page 16) as follows:

$$\text{Let } f(x) \triangleq \begin{cases} 1 & \text{for } 0 \le |x| < \dfrac{1-\beta}{2} \\ \dfrac{1}{2}\left\{ 1 + \cos\left[ \dfrac{\pi}{\beta}\left( |x| - \dfrac{1-\beta}{2} \right) \right] \right\} & \text{for } \dfrac{1-\beta}{2} \le |x| < \dfrac{1+\beta}{2} \\ 0 & \text{otherwise} \end{cases}$$

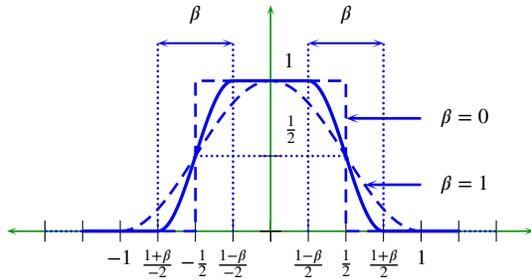

Then f forms a *partition of unity*:

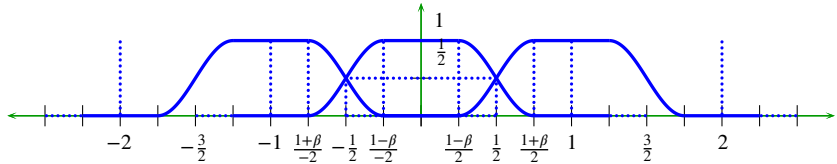

### 5.3.3 Scaling functions with partition of unity

The $Z$ transform (Definition 2.35 page 24) of a sequence $(h_n)$ with sum $\sum_{n \in \mathbb{Z}} (-1)^n h_n = 0$ has a zero at $z = -1$. Somewhat surprisingly, the *partition of unity* and *zero at $z = -1$* properties are actually equivalent (next theorem).

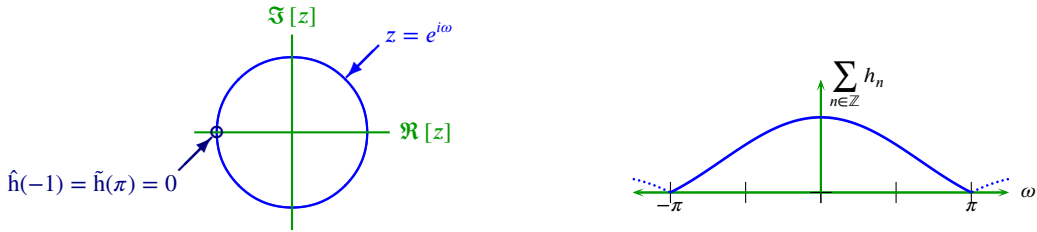

**Theorem 5.32** [145] *Let $\left( L_{\mathbb{R}}^2, \left( V_j \right), \phi, (h_n) \right)$ be a multiresolution system (Definition 5.10 page 71). Let $\tilde{F}f(\omega)$ be the FOURIER TRANSFORM (Definition 2.22 page 20) of a function $f \in L_{\mathbb{R}}^2$. Let $\bar{\delta}_n$ be the*

---

KRONECKER DELTA FUNCTION.

$$\underbrace{\sum_{n\in\mathbb{Z}}\mathbf{T}^n\phi = c \quad \textit{for some } c\in\mathbb{R}\setminus 0}_{\textit{(1) PARTITION OF UNITY}} \iff \underbrace{\sum_{n\in\mathbb{Z}}(-1)^n h_n = 0}_{\textit{(2) ZERO AT } z=-1} \iff \underbrace{\sum_{n\in\mathbb{Z}}h_{2n} = \sum_{n\in\mathbb{Z}}h_{2n+1} = \frac{\sqrt{2}}{2}}_{\textit{(3) sum of even = sum of odd} = \frac{\sqrt{2}}{2}}$$

✎PROOF:   Let $\mathbb{Z}_e$ be the set of even integers and $\mathbb{Z}_o$ the set of odd integers.

(1)   Proof that (1) $\impliedby$ (2):

$$\sum_{n\in\mathbb{Z}}\mathbf{T}^n\phi = \sum_{n\in\mathbb{Z}}\mathbf{T}^n\left[\sum_{m\in\mathbb{Z}}h_m\mathbf{D}\mathbf{T}^m\phi\right] \qquad \text{by dilation equation (Theorem 5.6 page 67)}$$

$$= \sum_{m\in\mathbb{Z}}h_m\sum_{n\in\mathbb{Z}}\mathbf{T}^n\mathbf{D}\mathbf{T}^m\phi$$

$$= \sum_{m\in\mathbb{Z}}h_m\sum_{n\in\mathbb{Z}}\mathbf{D}\mathbf{T}^{2n}\mathbf{T}^m\phi \qquad \text{by Proposition 4.8 page 41}$$

$$= \mathbf{D}\sum_{m\in\mathbb{Z}}h_m\sum_{n\in\mathbb{Z}}\mathbf{T}^{2n}\mathbf{T}^m\phi$$

$$= \mathbf{D}\sum_{m\in\mathbb{Z}}h_m\left[\sqrt{\frac{2\pi}{2}}\hat{\mathbf{F}}^{-1}\mathbf{S}_2\tilde{\mathbf{F}}(\mathbf{T}^m\phi)\right] \qquad \text{by PSF (Theorem 4.21 page 48)}$$

$$= \sqrt{\pi}\mathbf{D}\sum_{m\in\mathbb{Z}}h_m\hat{\mathbf{F}}^{-1}\mathbf{S}_2 e^{-i\omega m}\tilde{\mathbf{F}}\phi \qquad \text{by Corollary 4.18 page 47}$$

$$= \sqrt{\pi}\mathbf{D}\sum_{m\in\mathbb{Z}}h_m\hat{\mathbf{F}}^{-1}e^{-i\frac{2\pi}{2}km}\mathbf{S}_2\tilde{\mathbf{F}}\phi \qquad \text{by definition of } \mathbf{S} \text{ (Theorem 4.21 page 48)}$$

$$= \sqrt{\pi}\mathbf{D}\sum_{m\in\mathbb{Z}}h_m\hat{\mathbf{F}}^{-1}(-1)^{km}\mathbf{S}_2\tilde{\mathbf{F}}\phi$$

$$= \sqrt{\pi}\mathbf{D}\sum_{m\in\mathbb{Z}}h_m\left[\frac{\sqrt{2}}{2}\sum_{k\in\mathbb{Z}}(-1)^{km}\left(\mathbf{S}_2\tilde{\mathbf{F}}\phi\right)e^{i\frac{2\pi}{2}kx}\right] \qquad \text{by def. of } \hat{\mathbf{F}}^{-1} \text{ (Theorem 2.17 page 18)}$$

$$= \frac{\sqrt{2\pi}}{2}\mathbf{D}\sum_{k\in\mathbb{Z}}\left(\mathbf{S}_2\tilde{\mathbf{F}}\phi\right)e^{i\pi kx}\sum_{m\in\mathbb{Z}}(-1)^{km}h_m$$

$$= \frac{\sqrt{2\pi}}{2}\mathbf{D}\sum_{k\in\mathbb{Z}_e}\left(\mathbf{S}_2\tilde{\mathbf{F}}\phi\right)e^{i\pi kx}\sum_{m\in\mathbb{Z}}(-1)^{km}h_m$$

$$\qquad + \frac{\sqrt{2\pi}}{2}\mathbf{D}\sum_{k\in\mathbb{Z}_o}\left(\mathbf{S}_2\tilde{\mathbf{F}}\phi\right)e^{i\pi kx}\sum_{m\in\mathbb{Z}}(-1)^{km}h_m$$

$$= \frac{\sqrt{2\pi}}{2}\mathbf{D}\sum_{k\in\mathbb{Z}_e}\left(\mathbf{S}_2\tilde{\mathbf{F}}\phi\right)e^{i\pi kx}\underbrace{\sum_{m\in\mathbb{Z}}h_m}_{\sqrt{2}}$$

$$+ \frac{\sqrt{2\pi}}{2} \mathbf{D} \sum_{k \in \mathbb{Z}_o} (\mathbf{S}_2 \tilde{\mathbf{F}} \phi) e^{i\pi kx} \underbrace{\sum_{m \in \mathbb{Z}} (-1)^m h_m}_{0}$$

$$= \sqrt{\pi} \mathbf{D} \sum_{k \in \mathbb{Z}_e} (\mathbf{S}_2 \tilde{\mathbf{F}} \phi) e^{i\pi kx} \qquad \text{by Theorem 5.14 (page 74) and right hyp.}$$

$$= \sqrt{\pi} \mathbf{D} \sum_{k \in \mathbb{Z}_e} \tilde{\phi}\left(\frac{2\pi}{2}k\right) e^{i\pi kx} \qquad \text{by definitions of } \tilde{\mathbf{F}} \text{ and } \mathbf{S}_2$$

$$= \sqrt{\pi} \mathbf{D} \sum_{k \in \mathbb{Z}} \tilde{\phi}(2\pi k) e^{i2\pi kx} \qquad \text{by definition of } \mathbb{Z}_e$$

$$= \frac{1}{2} \mathbf{D} \left\{ \sqrt{2\pi} \sum_{k \in \mathbb{Z}} \tilde{\phi}(2\pi k) e^{i2\pi kx} \right\}$$

$$= \frac{1}{2} \mathbf{D} \sum_{n \in \mathbb{Z}} \phi(x + n) \qquad \text{by PSF (Theorem 4.21 page 48)}$$

$$= \frac{1}{2} \mathbf{D} \sum_{n} \mathbf{T}^n \phi \qquad \text{by definition of } \mathbf{T} \text{ (Definition 4.1 page 38)}$$

The above equation sequence demonstrates that

$$\mathbf{D} \sum_n \mathbf{T}^n \phi = \sqrt{2} \sum_n \mathbf{T}^n \phi$$

(essentially that $\sum_n \mathbf{T}^n \phi$ is equal to it's own dilation). This implies that $\sum_n \mathbf{T}^n \phi$ is a constant (Proposition 4.11 page 43).

(2)  Proof that (1) $\implies$ (2):

$$c = \sum_{n \in \mathbb{Z}} \mathbf{T}^n \phi \qquad \text{by left hypothesis}$$

$$= \sqrt{2\pi} \, \hat{\mathbf{F}}^{-1} \mathbf{S} \tilde{\mathbf{F}} \phi \qquad \text{by PSF (Theorem 4.21 page 48)}$$

$$= \sqrt{2\pi} \, \hat{\mathbf{F}}^{-1} \mathbf{S} \sqrt{2} \underbrace{\left( \mathbf{D}^{-1} \sum_{n \in \mathbb{Z}} h_n e^{-i\omega n} \right) \left( \mathbf{D}^{-1} \tilde{\mathbf{F}} \phi \right)}_{\tilde{\mathbf{F}} \phi} \qquad \text{by Lemma 5.7 page 67}$$

$$= 2\sqrt{\pi} \, \hat{\mathbf{F}}^{-1} \left( \mathbf{S} \mathbf{D}^{-1} \sum_{n \in \mathbb{Z}} h_n e^{-i\omega n} \right) (\mathbf{S}\tilde{\mathbf{F}}\mathbf{D}\phi) \qquad \text{by Corollary 4.18 page 47}$$

$$= 2\sqrt{\pi} \, \hat{\mathbf{F}}^{-1} \left( \mathbf{S} \frac{1}{\sqrt{2}} \sum_{n \in \mathbb{Z}} h_n e^{-i\frac{\omega}{2}n} \right) (\mathbf{S}\tilde{\mathbf{F}}\mathbf{D}\phi) \qquad \text{by Proposition 4.3 page 39}$$

$$= \sqrt{2\pi} \, \hat{\mathbf{F}}^{-1} \left( \sum_{n \in \mathbb{Z}} h_n e^{-i\frac{2\pi k}{2}n} \right) (\mathbf{S}\tilde{\mathbf{F}}\mathbf{D}\phi) \qquad \text{by def. of } \mathbf{S} \text{ (Theorem 4.21 page 48)}$$

$$= \sqrt{2\pi}\,\hat{\mathbf{F}}^{-1}\left(\sum_{n\in\mathbb{Z}} h_n(-1)^{kn}\right)(\mathbf{SD}^{-1}\mathbf{F}\phi)$$

$$= \sqrt{2\pi}\,\hat{\mathbf{F}}^{-1}\left(\sum_{n\in\mathbb{Z}} h_n(-1)^{kn}\right)\left(\mathbf{S}\frac{1}{\sqrt{2}}\tilde{\phi}\left(\frac{\omega}{2}\right)\right) \qquad \text{by def. of } \mathbf{S} \text{ (Theorem 4.21 page 48)}$$

$$= \sqrt{2\pi}\,\hat{\mathbf{F}}^{-1}\left(\sum_{n\in\mathbb{Z}} h_n(-1)^{kn}\right)\left(\frac{1}{\sqrt{2}}\tilde{\phi}\left(\frac{2\pi k}{2}\right)\right)$$

$$= \sqrt{\pi}\,\sum_{k\in\mathbb{Z}}\sum_{n\in\mathbb{Z}} h_n(-1)^{kn}\,\tilde{\phi}(\pi k)e^{i2\pi kx} \qquad \text{by def. of } \hat{\mathbf{F}}^{-1} \text{ (Theorem 2.17 page 18)}$$

$$= \sqrt{\pi}\,\sum_{k\text{ even}}\sum_{n\in\mathbb{Z}} h_n(-1)^{kn}\,\tilde{\phi}(\pi k)e^{i2\pi kx}$$
$$\quad + \sqrt{\pi}\,\sum_{k\text{ odd}}\sum_{n\in\mathbb{Z}} h_n(-1)^{kn}\,\tilde{\phi}(\pi k)e^{i2\pi kx}$$

$$= \sqrt{\pi}\,\sum_{k\text{ even}}\left(\underset{\textcolor{red}{\sqrt{2}}}{\sum_{n\in\mathbb{Z}} h_n}\right)\tilde{\phi}(\pi k)e^{i2\pi kx}$$
$$\quad + \sqrt{\pi}\,\sum_{k\text{ odd}}\left(\sum_{n\in\mathbb{Z}} h_n(-1)^{n}\right)\tilde{\phi}(\pi k)e^{i2\pi kx}$$

$$= \sqrt{\pi}\,\sum_{k\in\mathbb{Z}}\sqrt{2}\,\tilde{\phi}(\pi 2k)e^{i2\pi 2kx}$$
$$\quad + \sqrt{\pi}\,\sum_{k\in\mathbb{Z}}\left(\sum_{n\in\mathbb{Z}} h_n(-1)^{n}\right)\tilde{\phi}(\pi[2k+1])e^{i2\pi[2k+1]x} \qquad \text{by Theorem 5.14 page 74}$$

$$= \frac{\sqrt{2\pi}}{\sqrt{2\pi}}\tilde{\phi}(0) + \sqrt{\pi}e^{i2\pi x}\sum_{n\in\mathbb{Z}} h_n(-1)^{n}\sum_{k\in\mathbb{Z}}\tilde{\phi}(\pi[2k+1])e^{i4\pi kx} \qquad \text{by left hyp. and Theorem 5.27 page 86}$$

$$\implies \left(\sum_{n\in\mathbb{Z}} h_n(-1)^{n}\right) = 0 \qquad \text{because the right side must equal } c$$

(3)   Proof that (2) $\implies$ (3):

$$\sum_{n\in\mathbb{Z}_{\mathrm{e}}} h_n = \sum_{n\in\mathbb{Z}_{\mathrm{o}}} h_n = \frac{1}{2}\sum_{n\in\mathbb{Z}} h_n \qquad \text{by (2) and Proposition 2.41 page 26}$$
$$= \frac{\sqrt{2}}{2} \qquad \text{by } \textit{admissibility condition} \text{ (Theorem 5.14 page 74)}$$

(4)   Proof that (2) $\impliedby$ (3):

$$\frac{\sqrt{2}}{2} = \underbrace{\sum_{n\in\mathbb{Z}_{\mathrm{e}}}(-1)^{n}h_n}_{\text{even terms}} + \underbrace{\sum_{n\in\mathbb{Z}_{\mathrm{o}}}(-1)^{n}h_n}_{\text{odd terms}} \qquad \text{by (3)}$$

$$\implies \sum_{n \in \mathbb{Z}} (-1)^n h_n = 0 \qquad\qquad \text{by Proposition } 2.41 \text{ page } 26$$

**Proposition 5.33**

$\phi(x)$ *generates a* PARTITION OF UNITY $\implies$ $\phi(x)$ *generates an* MRA *system.*

*Telecommunications Engineering Department, National Chiao-Tung University, Hsinchu, Taiwan*

dgreenhoe@gmail.com

---

[146]This document was typeset using X∃LATEX, which is part of the TEX family of typesetting engines, which is arguably the greatest development since the Gutenberg Press. Graphics were rendered using the pstricks and related packages and LATEX graphics support. The main roman, *italic*, and **bold** font typefaces used are all from the Heuristica family of typefaces (based on the Utopia typeface, released by Adobe Systems Incorporated). The math font is XITS from the XITS font project. The font used for the title in the footers is Adventor (similar to Avant-Garde) from the TEX-Gyre Project. The font used for the version number in the footers is LIQUID CRYSTAL (Liquid Crystal) from FontLab Studio. This handwriting font is Lavi from the Free Software Foundation.